\documentclass{amsart}

\usepackage{amsmath,amsthm,amssymb,amsfonts,enumerate,color}
\usepackage{mathrsfs}
\usepackage{amstext,amsxtra}
\usepackage[margin=1.5in]{geometry}
\usepackage{graphicx}
\usepackage{float}

\usepackage{tikz}
\usetikzlibrary{matrix,arrows,calc,intersections,fit}
\usetikzlibrary{decorations.markings}
\usepackage{tikz-cd}
\usepackage{textcomp}
\usepgflibrary{shapes} 
\usepgflibrary[shapes] 
\usetikzlibrary{shapes} 
\usetikzlibrary[shapes] 

\usepackage[colorlinks,urlcolor=black,linkcolor=blue,citecolor=blue,hypertexnames=false]{hyperref}
\allowdisplaybreaks
\usepackage{pgf,tikz}
\usepackage{stmaryrd}
\usepackage{bm}

{\theoremstyle{plain}
\newtheorem{theorem}{Theorem}[section]
\newtheorem{proposition}[theorem]{Proposition}
\newtheorem{corollary}[theorem]{Corollary}
\newtheorem{lemma}[theorem]{Lemma}}
{\theoremstyle{definition}
\newtheorem{example}[theorem]{Example}
\newtheorem{definition}[theorem]{Definition}
\newtheorem{construction}[theorem]{Construction}}
{\theoremstyle{remark}
\newtheorem{remark}[theorem]{Remark}}

\newcommand{\mult}{\operatorname{mult}}
\newcommand{\aut}{\operatorname{Aut}}
\newcommand{\conj}{\operatorname{conj}}
\newcommand{\id}{\operatorname{id}}

\newcommand{\sym}{\operatorname{Sym}}

\newcommand{\hl}{\mathcal{H}}
\newcommand{\sal}{\mathcal{S}}

\newcommand{\cl}{\mathcal{C}}

\newcommand{\ml}{\mathcal{M}}

\newcommand{\fl}{\mathcal{F}}
\newcommand{\zl}{\mathcal{Z}}

\newcommand{\rb}{\mathbb{R}}

\newcommand{~}{\quad}
\newcommand{\cb}{\mathbb{C}}

\newcommand{\pb}{\mathbb{P}}

\newcommand{\sk}{\mathfrak{s}}

\newcommand{\undl}{\underline}

\definecolor{cardinalred}{RGB}{140,21,21}
\definecolor{coolgray}{RGB}{77,79,83}
\definecolor{black}{RGB}{0,0,0}
\definecolor{beige}{RGB}{210,194,149}
\definecolor{darkbeige}{RGB}{179,153,93}
\definecolor{darkcardinal}{RGB}{94,48,50}
\definecolor{lightcardinal}{RGB}{141,60,30}
\definecolor{darkpurple}{RGB}{83,40,79}
\definecolor{darkcyan}{RGB}{0,124,146}
\definecolor{skyblue}{RGB}{0,152,219}
\definecolor{seablue}{RGB}{10,100,180}
\definecolor{darkblue}{RGB}{20,80,150}
\definecolor{treegreen}{RGB}{0,155,118}
\definecolor{darkorange}{RGB}{168,101,12}
\definecolor{beigegray}{RGB}{95,87,79}
\definecolor{boxgray}{RGB}{238,235,233}
\definecolor{footergray}{RGB}{199,209,197}

\usepgflibrary{plotmarks} 
\usepgflibrary[plotmarks] 
\usetikzlibrary{plotmarks} 
\usetikzlibrary[plotmarks]

\begin{document}

\title{Asymptotics for real monotone double Hurwitz numbers}

\author{Yanqiao Ding}
\author{Qinhao He}

\address{School of Mathematics and Statistics, Zhengzhou University, Zhengzhou, 450001, China}

\email{yqding@zzu.edu.cn}

\address{School of Mathematics and Statistics, Zhengzhou University, Zhengzhou, 450001, China}

\email{1141786482@qq.com}



\date{\today}

\begin{abstract}
In recent years, monotone double Hurwitz numbers were introduced as a naturally combinatorial modification of double Hurwitz numbers.
Monotone double Hurwitz numbers share many structural properties with their classical counterparts, such as piecewise polynomaility, while the quantitative properties of these two numbers are quite different.
We consider real analogues of monotone double Hurwitz numbers and study the asymptotics for these real analogues.
The key ingredient is an interpretation of real tropical covers with arbitrary splittings as factorizations in the symmetric group which generalizes the result from Guay-Paquet, Markwig, and Rau (Int. Math. Res. Not. IMRN, 2016(1):258-293, 2016).
By using the above interpretation, we consider three types of real analogues of monotone double Hurwitz numbers: real monotone double Hurwitz numbers relative to simple splittings, relative to arbitrary splittings and real mixed double Hurwitz numbers. Under certain conditions,
we find lower bounds for these real analogues,
and obtain logarithmic asymptotics for real monotone double Hurwitz numbers relative to arbitrary splittings and real mixed double Hurwitz numbers.
In particular, under given conditions real mixed double Hurwitz numbers are logarithmically equivalent to complex double Hurwitz numbers.
We construct a family of real tropical covers and use them to show that real monotone double Hurwitz numbers relative to simple splittings are logarithmically equivalent to monotone double Hurwitz numbers with specific conditions.
This is consistent with the logarithmic equivalence of real double Hurwitz numbers and complex double Hurwitz numbers.
\end{abstract}

\maketitle

\tableofcontents

\section{Introduction}
Hurwitz numbers are important geometric invariants in enumerative geometry.
The study on Hurwitz numbers is related to many fields of mathematics, such as the moduli space of algebraic curves, integrable system, representation theory, combinatorics, tropical geometry and mathematical physics \cite{bbm-2011,cjm-2010,cm-2016,dyz-2017,elsv-2001,lzz-2000,gjv-2005,okounkov-2000,op-2006}.
Roughly speaking, Hurwitz numbers count the number of ramified covers of a Riemann surface by
Riemann surfaces with specified
ramification profiles over a fixed set of points.
An equivalent way to define Hurwitz numbers is to enumerate factorizations of the identity into a product of permutations of given cycle types in the symmetric group \cite{cm-2016,hurwitz-1891}.

The double Hurwitz number $H^\cb_g(\lambda,\mu)$ is of particular interest.
It enumerates ramified covers of $\cb P^1$ by genus $g$ surfaces with
ramification profiles $\lambda$, $\mu$ over $0$, $\infty$
and simple ramification over other branch points,
where $\lambda$ and $\mu$ are two partitions of an
integer $d$.
Many interesting results about the double Hurwitz number were obtained in recent decades.
For example, the polynomiality of
the generating function of double Hurwitz numbers and the wall-crossing
formulas \cite{cjm-2011,dl-2022,gjv-2005,johnson-2015,ssv-2008}.
In particular, methods from tropical geometry were applied to study double Hurwitz numbers \cite{bbm-2011,cjm-2010}.

In recent years, a modification of the double Hurwitz number was introduced by I. P. Goulden, M. Guay-Paquet, and J. Novak in \cite{ggpn-2014} which is called the monotone double Hurwitz number.
The monotone double Hurwitz number $\vec H^\cb_g(\lambda,\mu)$ is a combinatorial interpretation of the asymptotic expansion of the Harish-Chandra-Itzykson-Zuber (HCIZ) random matrix model.
Once the $d$ sheets of a ramified cover are labelled by integers $1,2,\ldots,d$ and a unramified point in $\cb P^1$ is chosen, one can construct a monodromy representation for the covering map.
By using the monodromy representation,
a simple ramification over a branch point corresponds to a transposition $(a_i,b_i)$ with $a_i<b_i$.
We denote by $l(\lambda)$ the number of parts of $\lambda$, and call $l(\lambda)$ the length of $\lambda$.
The monotone double Hurwitz number $\vec H^\cb_g(\lambda,\mu)$ counts the same covers which are counted by the double Hurwitz number $H^\cb_g(\lambda,\mu)$ with an additional condition:
\begin{equation}\label{eq:intro}
b_i\leq b_{i+1} \text{ for all } 1\leq i\leq r-1,
\end{equation}
where $r=l(\lambda)+l(\mu)+2g-2$ (See equation $(\ref{eq:monot-Hurwitz})$ for the precise definition).
If the condition in equation $(\ref{eq:intro})$ is modified to $b_i\leq b_{i+1}$ for all $1\leq i\leq k-1$,
where $k\leq r$, the corresponding count of ramified covers is called the mixed double Hurwitz number.
Since the introducing of monotone Hurwitz numbers, abundant researches on monotone Hurwitz numbers were carried out.
Goulden, Guay-Paquet, and Novak calculated single monotone Hurwitz numbers and derived the polynomiality of single monotone Hurwitz numbers in \cite{ggpn-2013,ggpn-2013a}.
There is an important quantitative difference between the two numbers $H^\cb_g(\lambda,\mu)$ and $\vec H^\cb_g(\lambda,\mu)$.
Based on explicit formulas for Hurwitz numbers and monotone Hurwitz numbers \cite{dyz-2017,ggpn-2013,ggpn-2013a}, the Hurwitz number grows superexponentially in the degree $d$, while its monotone analogue exhibits only exponential growth.
The convergence of monotone Hurwitz generating functions was further considered in \cite{ggpn-2017}.
A tropical approach to monotone Hurwitz numbers was introduced in \cite{dk-2017,hahn-2019} via the monodromy graph.
The polynomiality of monotone Hurwitz numbers was proved in \cite{hahn-2019,kls-2019}.
Another tropical interpretation of the monotone double Hurwitz number was introduced by Hahn and Lewa\'nski in \cite{hl-2022} which expresses monotone double Hurwitz numbers in terms of tropical covers weighted by Gromov-Witten invariants.
By applying their tropical approach to the monotone double Hurwitz number,
Hahn and Lewa\'nski derived wall-crossing and recursion formulas for tropical Jucys covers in \cite{hl-2020}.

We are interested in real analogues of double Hurwitz numbers.
The real double Hurwitz number $H^\rb_g(\lambda,\mu;s)$ counts
real ramified covers of $\cb P^1$ by genus $g$ surfaces with
ramification profiles $\lambda$, $\mu$ over $0$, $\infty$
and simple ramification over other branch points (See equation $(\ref{eq:def-real-Hurwitz})$ for the precise definition).
The number $s$ stands for the number of positive real branch points,
and the real double Hurwitz number $H^\rb_g(\lambda,\mu;s)$ depends on the number $s$.
Actually, in real enumerative geometry the number of real solutions for
a real enumerative problem usually depends on the positions of the point constraints \cite{iks2003,wel2005a,wel2005b}.
In order to solve real enumerative problems,
it is important to find lower bounds
for the numbers of real solutions.
For example, the signed counts of real algebraic curves in real surfaces or threefolds which are called the Welschinger invariants \cite{iks2013b,ks-2015,wel2005a,wel2005b}, and the real Gromov-Witten invariants of real symplectic manifolds \cite{gz2018}.
Itenberg and Zvonkine
\cite{iz-2018} found that the signed count method works in the study of real polynomials.
The signed count of real polynomials invented by Itenberg and Zvonkine is an invariant and
is logarithmically equivalent to the count of
complex polynomials under certain parity conditions.
El Hilany and Rau \cite{er-2019} generalized the construction of Itenberg and Zvonkine to the enumerative problem of real simple rational
functions $\frac{f(x)}{x-p}$, $f(x)\in\rb[x]$, $p\in\rb$.
Tropical geometry played an important role in the study of real algebraic geometry \cite{mikhalkin-2005}.
Guay-Paquet, Markwig, and Rau \cite{gpmr-2015} introduced a tropical approach to study real double Hurwitz numbers with positive real branch points based on factorizations in the symmetric group.
Markwig and Rau \cite{mr-2015} established the general theory of tropical real Hurwitz numbers and expressed real double Hurwitz numbers with arbitrary number of positive real branch points in terms of tropical covers weighted by multiplicities.
By using the tropical calculation of real double Hurwitz numbers derived in \cite{mr-2015},
Rau \cite{rau2019} analysed the combinatorial structure of tropical covers,
and obtained a lower bound, which is called the zigzag number,
of real double Hurwitz numbers.
Rau also proved the logarithmic equivalence of real double Hurwitz numbers
and complex double Hurwitz numbers under certain parity conditions.
The parity condition in \cite{rau2019} was removed by the first author in \cite{d-2020}.

In this paper, we consider real analogues of monotone double Hurwitz numbers and study the asymptotic growths of these real analogues when the degree goes to infinity and only
simple ramification points are added.
Our first problem is how to define the real counterparts of monotone double Hurwitz numbers.
The combinatorial nature of monotone double Hurwitz numbers makes it is not as obvious as in the case of real double Hurwtiz numbers, where we have a geometric definition by considering ramified covers over $\cb P^1$.
To solve this problem, we generalize the monodromy graph approach to real double Hurwitz numbers with real positive branch points introduced in \cite{gpmr-2015},
and interpret the tropical computation of real double Hurwitz numbers derived by Markwig and Rau in \cite{mr-2015} via factorizations in the symmetric group.
Then we construct real tropical covers with given real factorizations (c.f. Construction $\ref{const2}$).
With the help of Construction $\ref{const2}$, we give a combinatorial proof of Markwig and Rau's correspondence theorem \cite[Corollary $5.9$]{mr-2015}.
By using our factorization interpretation of real double Hurwitz numbers with arbitrary number of positive real branch points,
we introduce the real monotone double Hurwitz number $\vec{H}^\rb_g(\lambda,\mu;\undl S(s))$ (resp. the real $k$-mixed double Hurwitz number $\vec{H}^\rb_g(\lambda,\mu;\undl S(s),k)$) with $s$ positive real branch points under the splitting $\undl S(s)$
(See equation $(\ref{eq:signed-mixed-Hurwitz})$ and equation $(\ref{eq:mixed-Hurwitz})$ for the precise definition),
where $\undl S(s)=\{\sk_1,\ldots,\sk_r\}$ is a sequence of signs with $s$ positive entries (See Definition \ref{def-seq-signs}).
If the first $s$ entries in the sequence of signs $\undl S(s)$ are positive and all the remaining $r-s$ entries are negative,
we use $\vec H^\rb_g(\lambda,\mu;s)$ as the abbreviation of $\vec{H}^\rb_g(\lambda,\mu;\undl S(s))$.
We call the numbers
$$
\vec H^\rb_g(\lambda,\mu):=\inf_{0\leq s\leq r}\vec H^\rb_g(\lambda,\mu;s),\text{ and }
\vec\hl^{\rb}_g(\lambda,\mu):=\inf_{\undl S(s)}\vec H^\rb_g(\lambda,\mu;\undl S(s))
$$
the \textit{real monotone double Hurwitz number} relative to simple splittings and arbitrary splittings, respectively.
The number $\vec\hl^{\rb}_g(\lambda,\mu;k):=\inf_{\undl S(s)}\vec{H}^\rb_g(\lambda,\mu;\undl S(s),k)$ is called \textit{the real $k$-mixed double Hurwitz number}.
Our second goal is to study the asymptotic growths of $\vec H^\rb_g(\lambda,\mu)$, $\vec\hl^{\rb}_g(\lambda,\mu)$ and $\vec\hl^{\rb}_g(\lambda,\mu;k)$ as the degree goes to infinity and only
simple ramification points are added.
\subsection{Results}
We state our results about the asymptotics for real monotone (resp. mixed) double Hurwitz numbers as follows, and refer the readers to Corollary $\ref{cor:asym-simple-splitting}$, Corollary $\ref{cor:asymp-arbitrary-splitting}$ and Corollary $\ref{cor:asym-mix-zigzag1}$ for the precise statements.
\begin{theorem}
\label{thm:main}
With specific conditions, we have
\begin{enumerate}
    \item[$(1)$] real monotone double Hurwitz numbers relative to simple splittings are logarithmically equivalent to monotone double Hurwitz numbers;
    \item[$(2)$] the logarithmic asymptotics for real monotone double Hurwitz numbers relative to arbitrary splittings is at least $m\log m$ as the degree $d+2m$ goes to infinity;
    \item[$(3)$] real $k$-mixed double Hurwitz numbers are logarithmically equivalent to complex double Hurwitz numbers.
\end{enumerate}
\end{theorem}
The result $(1)$ in Theorem $\ref{thm:main}$ is consistent with other results about the logarithmic equivalence of real enumerative invariants and complex enumerative invariants.
For example, the logarithmic equivalence of Welschinger invariants and Gromov-Witten invariants in \cite{iks2004,iks2007,shustin2015},
the logarithmic equivalence of the signed counts of real polynomials (or simple rational functions) and the corresponding complex counts in \cite{er-2019,iz-2018}, and the logarithmic equivalence of real double Hurwitz numbers and complex double Hurwitz numbers in \cite{d-2020,rau2019}.


It follows from Proposition $\ref{prop:optimal-asymp}$ that
in the case when $\lambda=\mu=(1^{2m+1})$, the asymptotic growth $m\log m$ in Theorem $\ref{thm:main}(2)$ is the optimal estimate of real monotone double Hurwitz numbers relative to arbitrary splittings as the degree goes to infinity.
The difference between the asymptotic growths of real monotone double Hurwitz numbers relative to simple splittings and relative to arbitrary splittings shows the tremendous influence of the number of sign changes on the lower bounds of real monotone double Hurwitz numbers.

Our idea to prove the main theorem is as follows.
Zigzag covers proposed by Rau \cite{rau2019} have special combinatorial properties, \textit{i.e.} any zigzag cover has a unique real structure for any splitting of the branch points.
Since we want to study the asymptotics for real monotone double Hurwitz numbers, we only consider real monotone factorizations associated to zigzag covers according to Construction $\ref{const2}$.
Note that not every zigzag cover can be associated with a monotone factorization,
so we refine Rau's construction to the monotone setting.
We construct (universally) monotone zigzag covers such that for any real structure of a (universally) monotone zigzag cover there is at least one real monotone factorization associated to it.
Since the number of real monotone factorizations of certain type with a fixed starting permutation depends on the chosen permutation (see Example $\ref{exa:real-double-monotone}$),
we use the lower bounds of the numbers of real monotone factorizations associated to all real structures of a monotone zigzag cover as the multiplicity of the monotone zigzag cover.
We achieve the asymptotic growth of real monotone double Hurwitz numbers relative to arbitrary splittings by using universally monotone zigzag covers.
That is the optimal estimate of the real monotone double Hurwitz number relative to arbitrary splittings in the sense of Proposition $\ref{prop:optimal-asymp}$.
In order to get the optimal asymptotic growth of real monotone double Hurwitz numbers relative to simple splittings,
we first construct a family of monotone tropical covers,
then we glue a monotone zigzag cover with this family of monotone tropical covers.
It turns out that for any $s$ with a fixed parity, there are abundant tropical covers obtained in this way which have the required real structure.
The abundance of such tropical covers implies the asymptotic growth of real monotone double Hurwitz numbers relative to simple splittings is $2m\log m$ as the degree $d+2m$ goes to infinity.
The asymptotics for real $k$-mixed double Hurwitz numbers is obtained by applying the construction of universally monotone zigzag covers and Rau's estimate of zigzag numbers \cite{rau2019}.


\subsection{Organization of the paper}
In section $\ref{sec:2}$, we interpret real double Hurwitz numbers with arbitrary number of positive branch points as real factorizations of certain type and give a construction of real tropical covers with given real factorizations.
We construct (universally) monotone zigzag covers which provide lower bounds of real monotone double Hurwitz numbers in the next section.
In section $\ref{sec:asym-real-mono}$ and section $\ref{sec:5}$, we study the asymptotics for real monotone and mixed double Hurwitz numbers.

\subsection*{Acknowledgements}
The work on this text was initiated when the first author visited at Institut de Math\'ematiques
de Jussieu-Paris Rive Gauche in 2020.
The first author would like to thank IMJ-PRG for their
hospitality and excellent working conditions,
and he is also deeply grateful to Ilia Itenberg for
valuable discussions and suggestions.
The authors would like to thank Jianfeng Wu for a careful reading of the preliminary manuscript.
This work was supported by the Natural Science Foundation of Henan (No. 212300410287) and the Natural Science Foundation of China (No. 12101565).

\section{Correspondence theorem for real and tropical real double Hurwitz numbers revisited}
\label{sec:2}
In this section, we first describe real double Hurwitz numbers with arbitrary number of positive real branch points via monodromy representations following the approach in \cite{cadoret-2005,gpmr-2015} and \cite[Chapter $7$]{cm-2016},
then we build a relationship between the monodromy representations and tropical real double Hurwitz numbers with arbitrary number of positive real branch points which generalizes the result in \cite{gpmr-2015}.
Via this relationship, we give another proof of the correspondence theorem for real double Hurwitz numbers and their tropical counterparts established by Markwig and Rau \cite{mr-2015}.
Our another purpose of introducing the combinatorial properties of factorizations in the symmetric group in this section is to find possible ways to consider the real counterparts of monotone double Hurwitz numbers which will be discussed in Section $\ref{sec:real-mono}$.

\subsection{Preliminary on real double Hurwitz numbers}
\label{subsec:real-Hurwitz}
In this subsection, we review some basic facts about real double Hurwitz numbers from \cite{cadoret-2005,cm-2016,gpmr-2015,mr-2015,rau2019}.

Let $d\geq1$, $g\geq0$ be two fixed integers,
and let $\lambda$, $\mu$ be two partitions of $d$.
We assume that the set $\undl p=\{p_1,\ldots,p_r\}$
consists of $r=l(\lambda)+l(\mu)+2g-2$ points in $\cb P^1\setminus\{0,\infty\}$.
\begin{definition}\label{cpx-ramified-covering}
A \textit{Hurwitz cover of type $(g,\lambda,\mu,\undl p)$}
is a degree $d$ holomorphic map $\pi:C\to\cb P^1$ satisfying:
\begin{itemize}
    \item $C$ is a connected Riemann surface of genus $g$;
    \item $\pi$ ramifies with profiles $\lambda$
    and $\mu$ over $0$ and $\infty$ respectively;
    \item $\pi$ ramifies at $\undl p$ with simple branch points and is unramified everywhere else.
\end{itemize}
\end{definition}
An isomorphism of Riemann surfaces
$\psi:C_1\to C_2$ is called an isomorphism of two Hurwitz covers $\pi_1:C_1\to\cb P^1$ and $\pi_2:C_2\to\cb P^1$, if $\pi_1=\pi_2\circ\psi$.
Denote by $\aut^\cb(\pi)$ the space of isomorphic Hurwitz covers of the Hurwitz cover $\pi:C\to\cb P^1$.
The double Hurwitz number is the following weighted sum:
\begin{equation}
\label{eq:def-cpx-Hurwitz}
H^\cb_g(\lambda,\mu)=\sum_{[\pi]}
\frac{d!}{|\aut^\cb(\pi)|},
\end{equation}
where we sum over all isomorphism classes of
Hurwitz covers of type $(g,\lambda,\mu,\undl p)$.
The weighted sum $H^\cb_g(\lambda,\mu)$ does not
depend on the positions of the points in $\undl p$ (c.f. \cite{cm-2016}).

\begin{remark}\label{rem:def-Hurwitz-number}
Double Hurwitz numbers are usually defined to be
$\frac{1}{d!}H^\cb_g(\lambda,\mu)$ in the literature.
The number $\frac{1}{d!}H^\cb_g(\lambda,\mu)$ enumerates branched covers with unlabelled sheets.
Since the monotonicity condition, which is used to define monotone Hurwitz numbers,
depends on a total ordering of the sheets of a branched cover, the labelling matters in the study on monotone Hurwitz numbers.
For the sake of consistency, we consider Hurwitz number as an invariant counting the number of branched covers with labelled sheets and use the notation in equation $(\ref{eq:def-cpx-Hurwitz})$.
We refer the readers to see \cite[Remark $1.2$]{ggpn-2013} for more details.
\end{remark}

In the following, we assume further that the set $\undl p=\{p_1,\ldots,p_r\}\subset\rb P^1\setminus\{0,\infty\}$ and $p_1<p_2<\cdots<p_r$.
\begin{definition}\label{real-fct1}
A \textit{real Hurwitz cover of type $(g,\lambda,\mu,\undl p)$}
is a Hurwitz cover together with an anti-holomorphic involution $\tau:C\to C$
such that the cover $\pi:C\to\cb P^1$ satisfies $\pi\circ\tau=\conj\circ\pi$, where $\conj$ is the standard complex conjugation.
\end{definition}
An isomorphism of two real Hurwitz covers
$(\pi_1:C_1\to\cb P^1,\tau_1)$
and $(\pi_2:C_2\to\cb P^1,\tau_2)$
is an isomorphism of Riemann surfaces
$\varphi:C_1\to C_2$
such that $\pi_1=\pi_2\circ\varphi$ and
$\varphi\circ\tau_1=\tau_2\circ\varphi$.
We denote by $\aut^\rb(\pi,\tau)$ the space of isomorphic real Hurwitz covers of the real Hurwitz cover $(\pi,\tau)$.
Assume that $|\undl p\cap\rb^+|=s$,
then the real double Hurwitz number is defined as the weighted sum:
\begin{equation}
\label{eq:def-real-Hurwitz}
H^\rb_g(\lambda,\mu;s)=\sum_{[(\pi,\tau)]}
\frac{d!}{|\aut^\rb(\pi,\tau)|},
\end{equation}
where we sum over all isomorphism classes of
real Hurwitz covers of type $(g,\lambda,\mu,\undl p)$.
The weighted sum $H^\rb_g(\lambda,\mu;s)$
depends on the number $s$ of real positive points in $\undl p$.

\subsection{Real double Hurwitz numbers via factorizations}
An equivalent way to define the Hurwitz number
is to count the number of particular tuples in the symmetric group. In this subsection,
we match any real Hurwitz cover with a class of tuples in the symmetric group following the approach in \cite{cadoret-2005,gpmr-2015},
and describe real double Hurwitz numbers with arbitrary number of positive real branch points via the numbers of certain factorizations.

Let $\sal_d$ denote the symmetric group of order $d$.
The tuples in the symmetric group that we count are the images of generators of the fundamental group
$\pi_1(\cb P^1\setminus\{0,\infty,p_1,\ldots,p_r\},x_0)$ under the monodromy representation
$$
\Psi:\pi_1(\cb P^1\setminus\{0,\infty,p_1,\ldots,p_r\},x_0)\to\sal_d,
$$
where $x_0$ is an unramified point of $\pi$.
In the real case,
when we choose different base
points and different generators of the fundamental group
$\pi_1(\cb P^1\setminus\{0,\infty,p_1,\ldots,p_r\},x_0)$,
the monodromy map $\Psi$ induces collections of
tuples in the symmetric group with different combinatorial properties.
The reason for this phenomenon is that the combinatorial properties of the complex conjugate action on different generators of
$\pi_1(\cb P^1\setminus\{0,\infty,p_1,\ldots,p_r\},x_0)$ are different.
We give three examples to illustrate it in more detail.

\begin{example}
\label{exa:real-loops1}
Suppose that $x_0$ is a base point in $\rb P^1\setminus\{0,\infty,p_1,\ldots,p_r\}$
which is strictly smaller than all the points in $\undl p\cup\{0\}$.
Let $l_0,\ldots,l_r$ be $r+1$ loops around the holes of $\cb\setminus\{0,p_1,\ldots,p_r\}$ depicted in Figure $\ref{fig:real-loops1}$.
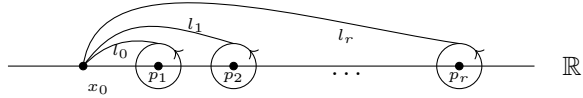
\begin{figure}[H]
    \centering
    \begin{tikzpicture}
    \draw (0,0)--(7,0);
    \foreach \Point in {(2,0),(3,0),(6,0)}
    \draw[decoration={markings, mark=at position 0.125 with {\arrow{>}}},
        postaction={decorate}
        ] \Point circle (0.3);
    \foreach \Point in {(1,0), (2,0),(3,0),(6,0)}
    \draw[fill=black] \Point circle (0.05);
    \draw (1.2,-0.3) node{\tiny{$x_0$}};
    \draw (4.5,-0.1) node{$\ldots$};
    \draw (7.5,0) node{$\rb$};
    \draw[bend left,-]  (1,0) to (2,0.3);
    \draw (1,0) .. controls (1.5,0.7) and (2,0.6) .. (3,0.3);
    \draw (1,0) .. controls (1.2,1.5) and (4,0.6) .. (6,0.3);
    \draw (1.5,0.2) node{\tiny{$l_0$}} (2.5,0.55) node{\tiny{$l_1$}} (4.5,0.4) node{\tiny{$l_r$}} (2,-0.15) node{\tiny{$p_1$}} (3,-0.15) node{\tiny{$p_2$}} (6,-0.15) node{\tiny{$p_r$}};
    \end{tikzpicture}
    \caption{Generators of $\pi_1(\cb P^1\setminus\{0,\infty,p_1,\ldots,p_r\},x_0)$ in Example \ref{exa:real-loops1}.}
    \label{fig:real-loops1}
\end{figure}
In this case, it is well known that the action of the complex
conjugation on the fundamental group $\pi_1(\cb P^1\setminus\{0,\infty,p_1,\ldots,p_r\},x_0)$
is determined by
\begin{equation}
\label{eq:exa1}
    \conj\circ(l_i\cdots l_0)=(l_i\cdots l_0)^{-1}.
\end{equation}
Equation $(\ref{eq:exa1})$ determines the combinatorical properties of the tuples in $\sal_d$ that we want to count.
\end{example}

\begin{example}\label{exa:real-loops2}
Suppose that the points in $\undl p$ satisfy
$p_1<\ldots<p_{r-s}<0<p_{r-s+1}<\ldots<p_r$.
Assume that the base point $x_0$ is chosen such that $p_{r-s}<x_0<0$.
Let $l_0,l_1,\ldots,l_r$ be $r+1$ loops depicted in
Figure $\ref{fig:real-loops2}$.
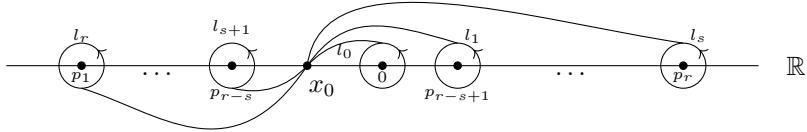
\begin{figure}[H]
    \centering
    \begin{tikzpicture}
    \draw (-3,0)--(7,0);
    \foreach\Point in {(-2,0), (0,0)}
    \draw[decoration={markings, mark=at position 0.125 with {\arrow{>}}},
        postaction={decorate}
        ] \Point circle (0.3);
    \foreach \Point in {(2,0),(3,0),(6,0)}
    \draw[decoration={markings, mark=at position 0.125 with {\arrow{>}}},
        postaction={decorate}
        ] \Point circle (0.3);
    \foreach \Point in {(-2,0),(0,0),(1,0), (2,0),(3,0),(6,0)}
    \draw[fill=black] \Point circle (0.05);
    \draw (-2,-0.15) node{\tiny{$p_1$}};
    \draw (0,-0.4) node{\tiny{$p_{r-s}$}};
    \draw (1.2,-0.3) node{$x_0$};
    \draw (2,-0.15) node{\tiny{$0$}};
    \draw (3,-0.4) node{\tiny{$p_{r-s+1}$}};
    \draw (6,-0.15) node{\tiny{$p_r$}};
    \draw (-1,-0.1) node{$\ldots$};
    \draw (4.5,-0.1) node{$\ldots$};
    \draw (7.5,0) node{$\rb$};
    \draw[bend left,-]  (1,0) to (2,0.3) (1.5,0.2) node{\tiny{$l_0$}};
    \draw (1,0) .. controls (1.5,0.7) and (2,0.6) .. (3,0.3);
    \draw (1,0) .. controls (1.2,1.5) and (4,0.6) .. (6,0.3);
    \draw[bend left,-]  (1,0) to (0,-0.3) (0,0.5)node{\tiny{$l_{s+1}$}};
    \draw (1,0) .. controls (0,-1.5) and (-1,-0.6) .. (-2,-0.3);
    \draw (3.2,0.4) node{\tiny{$l_1$}} (6.2,0.4) node{\tiny{$l_s$}} (-2,0.4) node{\tiny{$l_r$}};
    \end{tikzpicture}
    \caption{Generators of $\pi_1(\cb P^1\setminus\{0,\infty,p_1,\ldots,p_r\},x_0)$ in Example \ref{exa:real-loops2}.}
    \label{fig:real-loops2}
\end{figure}
The loops $l_0,l_1,\dots,l_r$ are generators of the
fundamental group $\pi_1(\cb P^1\setminus\{0,\infty,p_1,\ldots,p_r\},x_0)$.
The complex conjugate action on
$\pi_1(\cb P^1\setminus\{0,\infty,p_1,\ldots,p_r\},p_0)$
is determined by:
\begin{equation}\label{eq:exa2}
\begin{aligned}
    &\conj\circ(l_i\cdots l_0)=(l_i\cdots l_0)^{-1},
    \text{ for }0\leq i\leq s;    \\
    &\conj\circ(l_j\cdots l_{s+1})=(l_j\cdots l_{s+1})^{-1},
    \text{for }s+1\leq j\leq r.
\end{aligned}
\end{equation}
\end{example}

\begin{example}
\label{exa:real-loops3}
Let $\undl p$, $x_0$ and $l_0$, $l_1,\ldots,l_r$ be the same as in Example $\ref{exa:real-loops2}$.
Now we choose $r+1$ different loops $\bar l_0,\bar l_1,\ldots,\bar l_r$
around the holes of $\cb\setminus\{0,p_1,\ldots,p_r\}$ as follows:
first, we choose $r+1$ non-intersecting paths $c_0,c_1,\ldots,c_r$ which rotate counterclockwise from $x_0$ to the set $\{0,p_1,\ldots,p_r\}$.
We require that the end points of $c_0$, $c_1,\ldots,c_r$ are $0$, $p_{r-s}$, $p_{r-s+1},\ldots,p_i$, $p_{2r-2s-i+1},\ldots,p_1,p_{2r-2s}$, $p_{2r-2s+1},\ldots,p_r$, respectively, if $r\leq2s+1$,
and the end points of $c_0,c_1,\ldots,c_r$ are
$0$, $p_{r-s}$, $p_{r-s+1},\ldots,p_i$, $p_{2r-2s-i+1},\ldots,p_{r-2s+1},p_{r}$, $p_{r-2s},\ldots,p_1$, respectively, if $r>2s+1$.
Let $\theta_{i}$ be the angle between the
positive real axis and the curve $c_{i}$
at $x_0$.
We also require that the paths $c_1,\ldots,c_r$ are chosen such that $\theta_{i+1}>\theta_i$, $i=1,\ldots,r-1$.
Then we add small positive oriented loops around the end points of the paths $c_0,c_1,\ldots,c_r$ respectively.
Finally, we obtain the loops $\bar l_0,\bar l_1,\ldots,\bar l_r$ as depicted in Figure $\ref{fig:real-loops3}$.
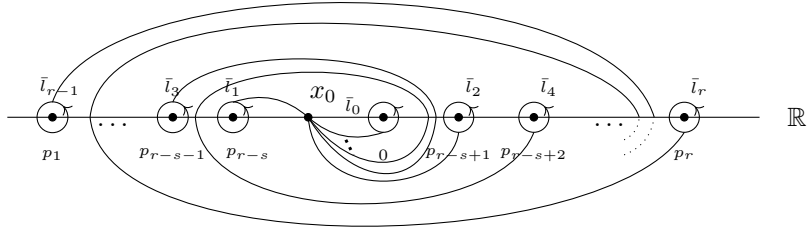
\begin{figure}[ht]
    \centering
    \begin{tikzpicture}
    \draw (-3,0)--(7,0);
    \foreach\Point in {(-2.4,0),(-0.8,0), (0,0)}
    \draw[decoration={markings, mark=at position 0.125 with {\arrow{>}}},
        postaction={decorate}
        ] \Point circle (0.205);
    \foreach \Point in {(2,0),(3,0),(4,0),(6,0)}
    \draw[decoration={markings, mark=at position 0.125 with {\arrow{>}}},
        postaction={decorate}
        ] \Point circle (0.2);
    \foreach \Point in {(-2.4,0),(-0.8,0),(0,0),(1,0), (2,0),(3,0),(4,0),(6,0)}
    \draw[fill=black] \Point circle (0.05);
    \draw[bend right,-]  (1,0) to (2,-0.2) (1.6,0.15) node{\tiny{$\bar l_0$}};
    \draw[bend right,-]  (1,0) to (0,0.2);
    \draw (1,0) .. controls (1.1,-1.2) and (2.9,-1) .. (3,-0.2);
    \draw[bend right,dotted]  (5.2,-0.3) to (5.4,0);
    \draw (5.4,0) .. controls (5,1.3) and (-1.7,2) .. (-1.9,0);
    \draw[bend right,dotted]  (5.2,-0.5) to (5.6,0);
    \draw (5.6,0) .. controls (5,2) and (-1.8,2) .. (-2.4,0.2);
    \draw (-1.9,0) .. controls (-1.7,-2) and (5,-1.8) .. (6,-0.2);
    \draw (1,0) .. controls (1.6,-1) and (2.6,-1) .. (2.7,0);
    \draw (2.7,0) .. controls (2.6,1) and (-0.5,1) .. (-0.8,0.2);
    \draw (1,0) .. controls (1.6,-0.8) and (2.5,-0.8) .. (2.6,0);
    \draw (2.6,0) .. controls (2.5,0.8) and (-0.4,0.8) .. (-0.5,0);
    \draw (4,-0.2) .. controls (3.5,-1.5) and (-0.4,-1.5) .. (-0.5,0);
    \draw[line width=0.4mm,dotted]  (1.5,-0.3) to (1.6,-0.5);
    \draw (-2.4,-0.5) node{\tiny{$p_1$}};
    \draw (0.2,-0.5) node{\tiny{$p_{r-s}$}};
    \draw (-0.8,-0.5) node{\tiny{$p_{r-s-1}$}};
    \draw (1.2,0.3) node{$x_0$};
    \draw (2,-0.5) node{\tiny{$0$}};
    \draw (3,-0.5) node{\tiny{$p_{r-s+1}$}};
    \draw (4,-0.5) node{\tiny{$p_{r-s+2}$}};
    \draw (6,-0.5) node{\tiny{$p_r$}} (0,0.4) node{\tiny{$\bar l_1$}} (-0.8,0.4) node{\tiny{$\bar l_3$}}(3.2,0.4) node{\tiny{$\bar l_2$}} (4.2,0.4) node{\tiny{$\bar l_4$}} (-2.3,0.4) node{\tiny{$\bar l_{r-1}$}} (6.2,0.4) node{\tiny{$\bar l_r$}};
    \draw (-1.6,-0.1) node{$\ldots$};
    \draw (5,-0.1) node{$\ldots$};
    \draw (7.5,0) node{$\rb$};
    \end{tikzpicture}
    \caption{Generators of $\pi_1(\cb P^1\setminus\{0,\infty,p_1,\ldots,p_r\},x_0)$, when $2s=r$, in Example \ref{exa:real-loops3}.}
    \label{fig:real-loops3}
\end{figure}
The loops $\bar l_0,\bar l_1,\dots,\bar l_r$ are generators of the
fundamental group $\pi_1(\cb P^1\setminus\{0,\infty,p_1,\ldots,p_r\},x_0)$.
In $\pi_1(\cb P^1\setminus\{0,\infty,p_1,\ldots,p_r\},x_0)$,
we have
$$
\begin{aligned}
    \bar l_0&\sim l_0,\\
    \bar l_1&\sim l_{s+1},\\
    \bar l_2&\sim \bar l_0^{-1}l_1\bar l_0,\\
    &\ldots,\\
    \bar l_{2i-1}&\sim (\bar L_1\bar L_3\cdots\bar L_{2i-5}\bar L_{2i-3})^{-1}l_{s+i} (\bar L_1\bar L_3\cdots\bar L_{2i-5}\bar L_{2i-3}), i\geq2\\
    \bar l_{2i}&\sim (\bar L_0\bar L_{2}\cdots\bar L_{2i-4}\bar L_{2i-2})^{-1} l_{i} (\bar L_0\bar L_{2}\cdots\bar L_{2i-4}\bar L_{2i-2}), i\geq2\\
    &\ldots,
\end{aligned}
$$
where $\bar L_k=\bar l_k\bar l_{k-1}\cdots\bar l_1\bar l_0$, $k\geq0$.
The relations in equation $(\ref{eq:exa2})$ imply:
\begin{equation}
\label{eq:exa3}
    \conj\circ(\bar L_0\bar L_1\cdots\bar L_{i-1}\bar L_i)=(\bar L_0\bar L_1\cdots\bar L_{i-1}\bar L_i)^{-1}, \text{ where } 0\leq i\leq r.
\end{equation}
\end{example}

Now we consider the constructions of loops around the holes
of $\cb\setminus\{0,p_1,\ldots,p_r\}$ in general.
\begin{definition}
\label{def-seq-signs}
A \textit{sequence $\undl S(s)$ of signs with $s$ positive entries} is a sequence of integers $\undl S(s)=\{\sk_1,\sk_2,\ldots,\sk_r\}$ such that
$\sk_i=\pm1$ and $|\rb^+\cap\undl S(s)|=s$.
A sequence $\undl S(s)$ of signs with $s$ positive entries is called a \textit{simple sequence},
if $\sk_1=\sk_2=\cdots=\sk_s=+1$ and $\sk_{s+1}=\cdots=\sk_r=-1$.
\end{definition}
It is easy to see that a simple sequence $\undl S(s)$ of signs with $s$ positive entries is determined by the number $s$.

\begin{construction}
\label{constr1}
Let $\undl S(s)=\{\sk_1,\sk_2,\ldots,\sk_r\}$ be any sequence of signs with $s$ positive entries.
Assume that $\undl p$, $x_0$ and $l_0$, $l_1,\ldots,l_r$
are the data given in Example $\ref{exa:real-loops2}$.
We choose additional $r+1$ loops $l_0',l_1',\ldots,l_r'$
around the $r+1$ holes in $\cb\setminus\{0,p_1,\ldots,p_r\}$.
The hole rounded by $l_i'$ is determined by the sequence of signs $\undl S(s)$, where $i=1,\ldots,r$.
\begin{itemize}
    \item $l_0'$ is always a loop around $0$ with base point $x_0$;
    \item $l_1'$ is a loop around the smallest positive number in $\undl p$, that is $p_{r-s+1}$,
    if $\sk_1=+1$. Otherwise, $l_1'$ is a loop around the largest negative number in $\undl p$, that is $p_{r-s}$;
    \item Assume that the loop $l_i'$ is already chosen to around a number $q\in\undl p$,
    then $l_{i+1}'$ is a loop around the number $q'\in\undl p$
    which is adjacent to $q$ and $|q'|>|q|$, if $\sk_{i+1}$
    has the same sign as $\sk_i$.
    Otherwise, $l_{i+1}'$ is a loop around the number
    $q''\in\undl p$, where $q''$ is the number with the minimal absolute value in the subset of $\undl p$
    consisting of numbers which have different signs from $q$
    and are not rounded by the loops $l_1',\ldots, l_i'$.
\end{itemize}
Loops $l_0', l_1',\ldots, l_r'$ are constructed
by fixing paths from $x_0$ to their end points and adding small positive oriented loops around their end points.
The paths from $x_0$ to their end points $\{0\}\cup\undl p$ are chosen
in the same way as the paths in Example \ref{exa:real-loops3}.
Note that the loops in Example \ref{exa:real-loops3}
are constructed according to the sequence of signs $\{\sk_1,\ldots,\sk_r\}$, where $\sk_i=+1$ for even $i$ and $\sk_i=-1$ for odd $i$, where $i=1,\ldots,r$.
\end{construction}

Obviously,
the $r+1$ loops $l_0', l_1',\dots, l_r'$ are the
generators of the fundamental group
$\pi_1(\cb P^1\setminus\{0,\infty,p_1,\ldots,p_r\},x_0)$.
For any sequence of signs
$\undl S(s)=\{\sk_1,\ldots,\sk_r\}$, there are integers
$i_1<i_2<\ldots<i_k$ dividing the
sequence $\sk_1,\ldots,\sk_r$ into consecutive
subsequences $\sk_1,\ldots,\sk_{i_1}$; $\sk_{i_1+1},\ldots,\sk_{i_2}$;
$\ldots$; $\sk_{i_k+1},\ldots,\sk_r$ such that all the entries in a
subsequence have same signs, while entries in
two adjacent subsequences have different signs.
The number $k$ of the integers $i_1,i_2,\ldots,i_k$ is exactly
the number of sign changes in the sequence
$\sk_1,\sk_2,\ldots,\sk_r$.

The above construction of loops can be used to describe the real double Hurwitz numbers with arbitrary number of positive branch points via factorizations in the symmetric group.
Before giving the proposition, we introduce some notations first.
\begin{itemize}
    \item if $\sk_1=-1$, we set
    $I_0=\{0\}$, $I_1=\{1,\ldots,i_1\}$, $I_2=\{i_1+1,\ldots,i_2\},\ldots$, $I_{k+1}=\{i_k+1,\ldots,r\}$.
    Let $L'_{I_0}=l'_0$ and $L'_{I_1}=l'_{i_1}\cdot l'_{i_1-1}\cdots l'_1$ be the compositions of certain loops.
\item if $\sk_1=+1$, we set
$I_1=\{0,1,\ldots,i_1\}$, $I_2=\{i_1+1,\ldots,i_2\},\ldots$, $I_{k+1}=\{i_k+1,\ldots,r\}$, and $L'_{I_1}=l'_{i_1}\cdot l'_{i_1-1}\cdots l'_1\cdot l'_0$.
\end{itemize}
No matter what the sign of $\sk_1$ is, we always set $L'_{I_j}=l'_{i_{j}}\cdot l'_{i_{j}-1}\cdots l'_{i_{j-1}+1}$
which is the composition of the loops labelled by $I_j$, where $j=2,\ldots,k$.
The combinatorial properties of the conjugate action on the loops $L'_{I_1}, L'_{I_2},\ldots,L'_{I_k}$
are the same as that on the loops
$\bar l_1,\bar l_2,\ldots,\bar l_r$ constructed in Example \ref{exa:real-loops3}.
From relations $(\ref{eq:exa3})$ in Example \ref{exa:real-loops3},
we obtain that
\begin{itemize}
    \item if $\sk_1=-1$,
\begin{equation}
\label{eq:exa4}
    \conj\circ(\bar L_0\bar L_1\cdots\bar L_{j-1}\bar L_j(m))=(\bar L_0\bar L_1\cdots\bar L_{j-1}\bar L_j(m))^{-1},
\end{equation}
for $0\leq j\leq k+1$ and $m\in I_j$,
where $\bar L_n=L'_{I_n}\cdots L'_{I_1}\cdot L'_{I_0}$ for
$0\leq n\leq k$ and $\bar L_j(m)=l'_m\cdots l'_{i_{j-1}+1}\cdot L'_{I_{j-1}}\cdots L'_{I_1}\cdot L'_{I_0}$;
   \item if $\sk_1=+1$,
\begin{equation}
\label{eq:exa5}
    \conj\circ(\tilde L_1\cdots\tilde L_{j-1}\tilde L_j(m))=(\tilde L_1\cdots\tilde L_{j-1}\tilde L_j(m))^{-1},
\end{equation}
for $1\leq j\leq k+1$ and $m\in I_j$,
where $\tilde L_n=L'_{I_n}\cdots L'_{I_1}$ for
$1\leq n\leq k$ and $\tilde L_j(m)=l'_m\cdots l'_{i_{j-1}+1}\cdot L'_{I_{j-1}}\cdots L'_{I_1}$;
\end{itemize}

We denote by $\cl(\sigma)\vdash\sal_d$ the cycle type of
$\sigma\in\sal_d$.
Let $d\geq1$ and $g\geq0$ be two integers, and let $\lambda$, $\mu$ be two partitions of $d$.
\begin{definition}\label{def:real-factor}
A \textit{real factorization of type $(g,\lambda,\mu;\undl S(s))$}
is a tuple $(\gamma,\sigma_1,\tau_1,\ldots,\tau_r,\sigma_2)$ of
elements of the symmetric group $\sal_d$
satisfying:
\begin{itemize}
    \item $\sigma_2\circ\tau_r\circ\cdots\circ\tau_1\circ\sigma_1=\id$;
    \item $r=l(\lambda)+l(\mu)+2g-2$;
    \item $\cl(\sigma_1)=\lambda$, $\cl(\sigma_2)=\mu$,
    $\cl(\tau_i)=(2,1,\ldots,1)$, $i=1,\ldots,r$;
    \item the subgroup generated by $\sigma_1$,
    $\sigma_2$, $\tau_1,\ldots,\tau_r$ acts transitively
    on the set $\{1,\ldots,d\}$.
    \item $\gamma$ is an involution (i.e. $\gamma^2=\id$)
    satisfying:
    \begin{equation}\label{eq:real-factor1}
    \begin{aligned}
    \gamma\circ\sigma_1\circ\gamma&=\sigma_1^{-1} \text{and},\\
          \gamma_i\circ(\tau_i\circ\cdot\cdot\cdot\circ\tau_{1}\circ\sigma_1)\circ\gamma_i&=
    (\tau_i\circ\cdot\cdot\cdot\circ\tau_{1}\circ\sigma_1)^{-1}, \text{ for } i=1,\ldots,r,
    \end{aligned}
    \end{equation}
    where
    $\gamma_1=\left\{
    \begin{aligned}
    \gamma,~ &\text{ if } \sk_1=1;\\
    \gamma\circ\sigma_1, &\text{ if } \sk_1=-1,
    \end{aligned}
    \right.$
    and
    $\gamma_{j+1}=\left\{
    \begin{aligned}
    \gamma_j,~~~~ &\text{ if } \sk_{j+1}=\sk_j;\\
    \gamma_j\circ\tau_j\circ\cdots\circ\tau_1\circ\sigma_1, &\text{ if } \sk_{j+1}\neq \sk_j.
    \end{aligned}
    \right.$
\end{itemize}
\end{definition}
We denote by $\fl^\rb(g,\lambda,\mu;\undl S(s))$ the set of
all real factorizations of type $(g,\lambda,\mu;\undl S(s))$.
\begin{lemma}\label{lem:realDH1}
Let $g\geq0$, $d\geq1$ be two integers, and
let $\lambda$, $\mu$ be two partitions of $d$.
Suppose that $\undl S(s)=\{\sk_1,\sk_2,\ldots,\sk_r\}$ is any
sequence of signs with $s$ positive entries.
Then
$$
H^\rb_g(\lambda,\mu;s)=|\fl^\rb(g,\lambda,\mu;\undl S(s))|.
$$
\end{lemma}

\begin{proof}
The proof of this Lemma is almost the same
as the proof of \cite[Lemma $2.3$]{gpmr-2015}.
When the sequence $\undl S(s)$ is a simple sequence,
the proof of the statement is sketched in \cite[Lemma A.2]{d-2020}.
The idea is to find a bijection between the set of
real factorizations of type $(g,\lambda,\mu;\undl S(s))$ modulo the action of $\sal_d$ and the set of real covers modulo real isomorphisms.
We only explain why Construction $\ref{constr1}$
induces real factorizations of type
$(g,\lambda,\mu;\undl S(s))$ here,
and recommend the readers to \cite[Lemma $2.3$]{gpmr-2015}
for the rest of the proof.

Let $l'_0,l'_1,\ldots,l'_r$ be the $r+1$ loops with base point
$x_0$ constructed in Construction $\ref{constr1}$.
Let $\pi:C\to\cb P^1$ be a real Hurwitz cover of type $(g,\lambda,\mu,\undl p)$.
We suppose that $\pi^{-1}(x_0)=\{q_1,q_2,\ldots,q_d\}$.
Denote by $\sigma_1,\tau_1,\ldots,\tau_r$ the
monodromy actions of the loops $l'_0,l'_1,\ldots,l'_r$
on $\{q_1,q_2,\ldots,q_d\}$, respectively.
The action of the complex conjugation on
$\{q_1,q_2,\ldots,q_d\}$ is described by $\gamma$.
The combinatorial relations in equation $(\ref{eq:exa4})$ and
equation $(\ref{eq:exa5})$ imply
the relation $(\ref{eq:real-factor1})$.
In fact,
if $\sk_1=-1$, the equation $(\ref{eq:exa4})$ implies that
$$
\gamma\circ[\sigma_1\circ(\tau_{i_1}\circ\cdots\circ\tau_1\circ\sigma_1)\circ\cdots\circ(\tau_i\circ\cdots\circ\tau_1\circ\sigma_1)]\circ\gamma=[\sigma_1\circ(\tau_{i_1}\circ\cdots\circ\tau_1\circ\sigma_1)\circ\cdots\circ(\tau_i\circ\cdots\circ\tau_1\circ\sigma_1)]^{-1}.
$$
Hence we have $\gamma_i\circ(\tau_i\circ\cdots\circ\tau_1\circ\sigma_1)\circ\gamma_i=(\tau_i\circ\cdots\circ\tau_1\circ\sigma_1)^{-1}$,
where $\gamma_i$ is determined by relation $(\ref{eq:real-factor1})$.
If $\sk_1=+1$, the equation $(\ref{eq:exa5})$ implies the relation $(\ref{eq:real-factor1})$.
\end{proof}

\subsection{Tropical real double Hurwitz numbers}
Let us recall some notations about tropical cover first.
The readers may refer to \cite{cjm-2010,gpmr-2015,mr-2015,rau2019} for more details.

Unless otherwise specified,
we only consider
connected graph $\Gamma$ without $2$-valent vertices in this paper.
By $\Gamma^\circ$ we denote
the subgraph obtained by removing the $1$-valent vertices
of $\Gamma$.
The {\it genus} $g$ of $\Gamma$ is the first Betti number
$b_1(\Gamma)$ of the graph $\Gamma$.
A {\it tropical curve} $C$ is a metric graph
such that the length of an end is infinite,
and the length $\ell(e)\in\rb$ of an inner edge $e$ is finite.
An isometric homeomorphism $\varPhi:C_1^\circ\to C_2^\circ$
of two tropical curves $C_1$, $C_2$ is called an isomorphism $\varPhi:C_1\to C_2$ of the two
tropical curves.
We consider the tropical projective line $T\pb^1$
as $\rb\cup\{\pm\infty\}$ in the following.

\begin{definition}\label{def:tro-cover}
A \textit{tropical cover} $\varphi:C\to T\pb^1$ is a continuous
map satisfying:
\begin{itemize}
    \item Inner vertices condition: $\varphi(p)\in\rb$,
    for any inner vertex $p$ of $C$.
    We denote by $\undl x$ the set of images of inner vertices
    of $C$, and call $\undl x$ the inner vertices of $T\pb^1$;
    \item Leaves condition: the set of leaves of $C$ is mapped onto
    $\{\pm\infty\}$ by $\varphi$;
    \item piecewise linearity: for any
    edge $e$ of $C$, if we consider $e$ as an interval
    $[0,\ell(e)]$, there is a positive integer $\omega(e)$
    such that
    $$
    \varphi(t)=\pm\omega(e)t+\varphi(0), \forall t\in[0,\ell(e)].
    $$
    The integer $\omega(e)$ is called the \textit{weight} of $e$.
    \item Balancing condition: For any vertex $v\in C$, we choose an
    edge $e'\subset T\pb^1$ adjacent to $\varphi(v)$.
    Then the integer
    $$
    \deg(\varphi,v):=\sum_{
    \substack{e\text{ edge of } C\\
    v\in e,e'\subset\varphi(e)}}
    \omega(e)
    $$
    does not depend on the choice of $e'$.
\end{itemize}
\end{definition}

\begin{definition}\label{def:degree}
Let $\varphi:C\to T\pb^1$ be a tropical cover.
For any edge $e'$ of $T\pb^1$,
the balancing condition implies that the sum
$$
\deg(\varphi):=\sum_{
    \substack{e\text{ edge of } C\\
    e'\subset\varphi(e)}}\omega(e)
$$
is independent of $e'$, and it is called the \textit{degree} of $\varphi$.
\end{definition}

\begin{remark}\label{rem:trop-monodromy}
For a real tropical cover $\varphi:C\to T\pb^1$,
we obtain a graph with the same combinatorial
properties by forgetting the lengths of edges of
the tropical curve $C$.
Such a graph is called a monodromy graph.
In the following, we do not distinguish these
two notations.
\end{remark}

For a tropical cover $\varphi:C\to T\pb^1$,
a symmetric cycle of $\varphi$ is a pair of edges with the same weight
and adjacent to the same two vertices.
A symmetric fork of $\varphi$ is a pair of ends with the same weight and
adjacent to a same inner vertex.
We denote by $CF(\varphi)$ the set of symmetric circles and symmetric forks of $\varphi:C\to T\pb^1$.
Let $C(\varphi)$ denote the set of symmetric circles of
$\varphi:C\to T\pb^1$.

Let $d\geq1$ and $g\geq0$ be two integers. Let $\lambda$ and $\mu$ be two partitions of $d$.
Suppose $r=l(\lambda)+l(\mu)+2g-2>0$.
Fix $r$ points $\undl x=\{x_1,\ldots,x_r\}\subset\rb$
satisfying $x_1<\ldots<x_r$.

\begin{definition}
A \textit{real tropical cover $(\varphi,\rho)$ of type
$(g,\lambda,\mu,\undl x)$} is a tropical cover
$\varphi:C\to T\pb^1$ of degree $d$ together with
a choice of subset $I_\rho\subset CF(\varphi)$ such that
\begin{itemize}
    \item $C$ is a tropical curve of genus $g$;
    \item $\lambda$ (resp. $\mu$) is the tuple of weights of the
    ends adjacent to all the leaves which are mapped to $-\infty$
    (resp. $+\infty$) by $\varphi$;
    \item $\undl x$ is the set of images of the inner
    vertices of $C$;
    \item there is a choice of a colour red or blue for every component
    of the subgraph of edges of even weights in $C\setminus I_\rho^\circ$.
\end{itemize}
\end{definition}
A choice of subset $I_\rho\subset CF(\varphi)$ and
a choice of a colour red or blue for every component
of the subgraph of edges of even weights in $C\setminus I_\rho^\circ$ are called a colouring $\rho$ of $\varphi:C\to T\pb^1$.
We denote by $E(I_\rho)$ the set of inner edges of even
weights in $C\setminus I_\rho^\circ$.
An {\it isomorphism} of two real tropical covers
$(\varphi_1:C_1\to T\pb^1;\rho_1)$ and
$(\varphi_2:C_2\to T\pb^1;\rho_2)$ is
an isomorphism $\varPhi:C_1\to C_2$ of tropical curves
respecting the covering maps and the colourings.
The \textit{real multiplicity} of a real tropical cover
$(\varphi,\rho)$ is
\begin{equation}
\label{eq:mult-real-trop-cover}
\mult^\rb(\varphi,\rho)
:=2^{|E(I_\rho)|-|CF(\varphi)|}\prod_{c\in I_\rho\cap C(\varphi)}\omega(c),
\end{equation}
where $\omega(c)$ is the weight of one edge of the
symmetric cycle $c$.

\begin{remark}\label{rem:mult-comparision}
The multiplicity $\mult^\rb(\varphi,\rho)$ in equation $(\ref{eq:mult-real-trop-cover})$ is the same as that in \cite[Definition $5.1$]{mr-2015} and \cite[Definition $3.6$]{gpmr-2015}.
It is different from the multiplicity used in \cite[Definition $3.3$]{rau2019}.
Since the subset $I_\rho$ in \cite[Definition $3.3$]{rau2019} is only allowed to contain symmetric cycles and odd symmetric forks,
the multiplicity in \cite[Definition $3.3$]{rau2019} is compensated with a factor.
If there is no even symmetric fork in the tropical cover $\varphi:C\to T\pb^1$, the multiplicity $\mult^\rb(\varphi,\rho)$ in \cite[Definition $5.1$]{mr-2015} and \cite[Definition $3.6$]{gpmr-2015} coincides with that in \cite[Definition $3.3$]{rau2019}.
\end{remark}

\begin{figure}[H]
    \centering
    \begin{tikzpicture}
    \draw[line width=0.4mm] (-3,0)--(-2,0)--(-1,0.5);
    \draw[line width=0.4mm,darkblue] (-2,0)--(-1,-0.5);
    \draw[line width=0.4mm,darkblue] (-0.5,0,0)--(0.5,0)--(1.5,0.5);
    \draw[line width=0.4mm,darkblue] (0.5,0)--(1.5,-0.5);
    \draw[line width=0.4mm,cardinalred] (-3,-1.3)--(-2,-1.3);
    \draw[line width=0.4mm] (-1,-0.8)--(-2,-1.3)--(-1,-1.8);
    \draw[line width=0.4mm,darkblue] (-0.5,-1.3)--(0.5,-1.3);
    \draw[line width=0.4mm,dotted] (1.5,-0.8)--(0.5,-1.3)--(1.5,-1.8);
    \draw (-1,-2.5) node{$(1)$ Positive vertices};
    \draw[line width=0.4mm] (4,0)--(5,0)--(6,0.5);
    \draw[line width=0.4mm,cardinalred] (5,0)--(6,-0.5);
    \draw[line width=0.4mm,cardinalred] (6.5,0)--(7.5,0)--(8.5,0.5);
    \draw[line width=0.4mm,cardinalred] (7.5,0)--(8.5,-0.5);
    \draw[line width=0.4mm,darkblue] (4,-1.3)--(5,-1.3);
    \draw[line width=0.4mm] (6,-0.8)--(5,-1.3)--(6,-1.8);
    \draw[line width=0.4mm,cardinalred] (6.5,-1.3)--(7.5,-1.3);
    \draw[line width=0.4mm,dotted] (8.5,-0.8)--(7.5,-1.3)--(8.5,-1.8);
    \draw (6,-2.5) node{$(2)$ Negative vertices};
    \end{tikzpicture}
    \caption{Signed vertices: even edges are drawn in colours, odd edges in black. Dotted edges are
    the symmetric cycles or forks contained in $I_\rho$.}
    \label{fig:signed-vertices}
\end{figure}
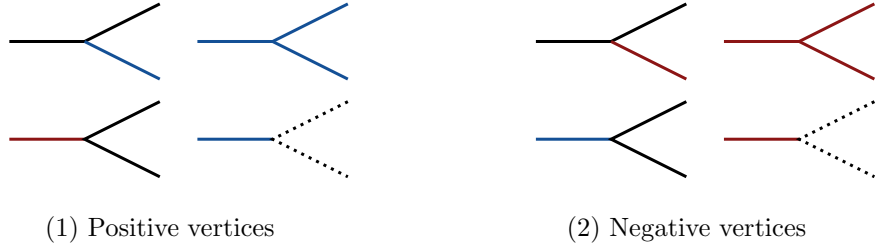

For a real tropical cover $(\varphi,\rho)$ of type
$(g,\lambda,\mu,\undl x)$,
a inner vertex $x_i\in\undl x$ is called a
{\it positive} or {\it negative} point if
it is the image of a $3$-valent vertex of $C$
depicted in Figure $\ref{fig:signed-vertices}(1)$ or
Figure $\ref{fig:signed-vertices}(2)$,
respectively, up to reflection along a vertical line.
Let $\undl x^+$ and $\undl x^-$ denote the collection of
positive and negative points respectively.
It is easy to see that a colouring $\rho(\varphi)$ of
$\varphi:C\to T\pb^1$ induces a splitting of $\undl x$
into positive and negative branch points
$\undl x^+\sqcup\undl x^-$.
Let $|\lambda|$ denote the sum of parts of $\lambda$.

\begin{theorem}[{\cite[Corollary $5.9$]{mr-2015}}]
\label{thm:real-tropical-Hurwitz}
Let $\lambda$ and $\mu$ be two partitions with $|\lambda|=|\mu|$, and let $g\geq0$ be an integer.
Suppose that the set $\undl x=\{x_1,\ldots,x_r\}\subset\rb$,
where $r=l(\lambda)+l(\mu)+2g-2>0$, has a splitting
$\undl x=\undl x^+\sqcup\undl x^-$ such that $|\undl x^+|=s$. Then
$$
H^\rb_g(\lambda,\mu;s)=\sum_{[(\varphi,\rho)]}d!\cdot\mult^\rb(\varphi,\rho),
$$
where we sum over all isomorphism classes $[(\varphi,\rho)]$
of real tropical covers of type $(g,\lambda,\mu,\undl x)$
whose positive and negative branch points reproduce the splitting $\undl x^+,\undl x^-$.
\end{theorem}

\subsection{Another proof of Markwig-Rau's correspondence theorem}
In this subsection, we illustrate real tropical covers $(\varphi,\rho)$ via factorizations in the symmetric group,
which generalizes the constructions in \cite[Section $3$]{gpmr-2015}.
By using our generalization, we give another proof of \cite[Corollary $5.9$]{mr-2015}.

From Lemma \ref{lem:realDH1}, the real double Hurwitz number
$H^\rb_g(\lambda,\mu;s)$ is determined by the number of
real factorizations of type $(g,\lambda,\mu;\undl S(s))$ for
any sequence $\undl S(s)=\{\sk_1,\sk_2,\ldots,\sk_r\}$ of signs with $s$ positive entries.
The combinatorial properties of the real factorizations
of type $(g,\lambda,\mu;\undl S(s))$ are determined by equation
$(\ref{eq:real-factor1})$. The effect of the conjugation with
$\gamma$ is concluded in \cite[Lemma $3.12$]{gpmr-2015}. For the convenience of the readers,
we recall \cite[Lemma $3.12$]{gpmr-2015} in the following.

\begin{lemma}[{\cite[Lemma $3.12$]{gpmr-2015}}]
\label{lem:action-involution}
Let $\gamma,\sigma\in\sal_d$ satisfy $\gamma\circ\sigma\circ\gamma=\sigma^{-1}$ and $\gamma^2=\id$.
Suppose that $\sigma=c_k\circ c_{k-1}\circ\cdots\circ c_1$
is the disjoint cycle decomposition of $\sigma$.
Then the conjugation with $\gamma$ can
\begin{itemize}
    \item either exchange two cycles of the same length, that is,
    $$
    \gamma\circ c_i\circ\gamma=c_j^{-1} \text{ and }
    \gamma\circ c_j\circ\gamma=c_i^{-1} \text{ for some } i\neq j;
    $$
    \item or invert a cycle $c_i$, that is $\gamma\circ c_i\circ\gamma=c_i^{-1}$. Moreover,
    \begin{itemize}
        \item if $c_i$ is of odd length,
         there is exactly one element in the cycle $c_i$
         fixed by the conjugation with $\gamma$;
         \item if $c_i$ is of even length,
         then either there are two elements which are of distance $\frac{l(c_i)}{2}$ fixed by the conjugation with $\gamma$,
         or there is no fixed point of the conjugation but a pair
         of subsequences consisting of $\frac{l(c_i)}{2}$
         consecutive elements in $c_i$ exchanged by the conjugation.
    \end{itemize}
\end{itemize}
\end{lemma}

We give two examples to explain the effect of the conjugation.
\begin{example}\label{exa:action-involution}
Let $c_1=(abcde)$, $c_2=(abcfde)$, $c_3=(bgcdhe)$ be three cycles in $\sal_8$.
$\gamma=(be)\cdot(cd)\cdot(gh)$ is an involution in $\sal_8$
satisfying $\gamma\circ c_i\circ\gamma=c_i^{-1}$,
$i=1,2,3$.
We use circle diagrams in Figure $\ref{fig:effect-conj}$ to describe cycles in $\sal_8$.
The action of $\gamma$ on cycles $c_1,c_2,c_3$ is
described by Lemma \ref{lem:action-involution}.
The cycle $c_1$ in Figure $\ref{fig:effect-conj}$ is of odd length, so it has one fixed point of $\gamma$, that is $a$.
Cycles $c_2$ and $c_3$ in Figure $\ref{fig:effect-conj}$
are of even lengths. The cycle $c_2$ has two fixed points, that are $a$ and $f$,
while the cycle $c_3$ has no fixed point.
\begin{figure}[ht]
    \centering
    \begin{tikzpicture}
    \draw[decoration={markings, mark=at position 0.18 with {\arrow{<}}},
        postaction={decorate}
        ] (-2,0) circle (1);
    \draw (-2,1) node {$a$} (-2.7,0.7) node{$e$} (-1.3,0.7) node{$b$} (-1.3,-0.7) node{$c$} (-2.7,-0.7) node{$d$};
    \draw [stealth-stealth](-2.5,0.7) -- (-1.5,0.7);
    \draw [stealth-stealth](-2.5,-0.7) -- (-1.5,-0.7);
    \draw[gray,dashed](-2,1.3) -- (-2,-1.3);
    \draw[decoration={markings, mark=at position 0.18 with {\arrow{<}}},
        postaction={decorate}
        ] (1,0) circle (1);
    \draw (1,1) node {$a$} (0.3,0.7) node{$e$} (1.7,0.7) node{$b$} (1.7,-0.7) node{$c$} (0.3,-0.7) node{$d$} (1,-1) node{$f$};
    \draw [stealth-stealth](0.5,0.7) -- (1.5,0.7);
    \draw [stealth-stealth](0.5,-0.7) -- (1.5,-0.7);
    \draw[gray,dashed](1,1.3) -- (1,-1.3);
    \draw[decoration={markings, mark=at position 0.18 with {\arrow{<}}},
        postaction={decorate}
        ] (4,0) circle (1);
    \draw (3,0) node {$h$} (3.3,0.7) node{$e$} (4.7,0.7) node{$b$} (4.7,-0.7) node{$c$} (3.3,-0.7) node{$d$} (5,0) node{$g$};
    \draw [stealth-stealth](3.5,0.7) -- (4.5,0.7);
    \draw [stealth-stealth](3.5,0) -- (4.5,0);
    \draw [stealth-stealth](3.5,-0.7) -- (4.5,-0.7);
    \draw[gray,dashed](4,1.3) -- (4,-1.3);
    \draw (-2,-1.8) node{$c_1=(abcde)$} (1,-1.8) node{$c_2=(abcfde)$} (4,-1.8) node{$c_3=(bgcdhe)$};
    \end{tikzpicture}
    \caption{Effect of conjugation with $\gamma$.}
    \label{fig:effect-conj}
\end{figure}
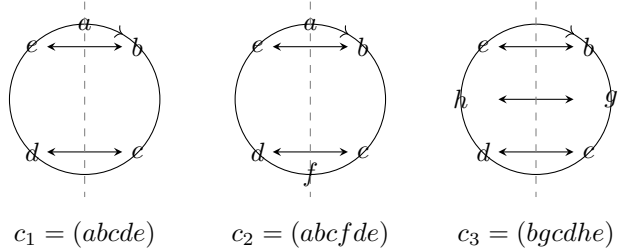
\end{example}

\begin{example}\label{exa:action-involution-var}
Let $c_1,c_2,c_3,\gamma\in\sal_8$ be the same as
Example \ref{exa:action-involution}.
It is easy to see that $\gamma_i:=\gamma\circ c_i$,
$i=1,2,3$, is an involution, and
$\gamma_i\circ c_i\circ\gamma_i=c_i^{-1}$.
Note that $\gamma_1=(ae)\cdot(bd)\cdot(gh)$, $\gamma_2=(ae)\cdot(bd)\cdot(cf)$,
$\gamma_3=(bh)\cdot(dg)$.
The effect of conjugation with $\gamma_i$ on $c_i$
is described by Figure \ref{fig:effect-conj-var}.
\begin{figure}[ht]
    \centering
    \begin{tikzpicture}
    \draw[decoration={markings, mark=at position 0.18 with {\arrow{<}}},
        postaction={decorate}
        ] (-2,0) circle (1);
    \draw (-2,1) node {$a$} (-2.9,0.5) node{$e$} (-1.1,0.5) node{$b$} (-1.3,-0.7) node{$c$} (-2.7,-0.7) node{$d$};
    \draw [stealth-stealth](-2.7,0.6) -- (-2.1,0.95);
    \draw [stealth-stealth](-2.5,-0.5) -- (-1.2,0.5);
    \draw[gray,dashed](-2.6,1.1) -- (-1.1,-1);
    \draw[decoration={markings, mark=at position 0.18 with {\arrow{<}}},
        postaction={decorate}
        ] (1,0) circle (1);
    \draw (1,1) node {$a$} (0.1,0.5) node{$e$} (1.9,0.5) node{$b$} (1.9,-0.5) node{$c$} (0.1,-0.5) node{$d$} (1,-1) node{$f$};
    \draw [stealth-stealth](0.25,0.5) -- (0.8,0.9);
    \draw [stealth-stealth](0.25,-0.5) -- (1.7,0.5);
    \draw [stealth-stealth](1.2,-0.9) -- (1.7,-0.55);
    \draw[gray,dashed](0.3,1) -- (1.7,-1);
    \draw[decoration={markings, mark=at position 0.18 with {\arrow{<}}},
        postaction={decorate}
        ] (4,0) circle (1);
    \draw (3,0) node {$h$} (3.3,0.7) node{$e$} (4,1) node{$b$} (4.7,-0.7) node{$c$} (4,-1) node{$d$} (5,0) node{$g$};
    \draw [stealth-stealth](3.1,0.1) -- (3.9,0.9);
    \draw [stealth-stealth](4.1,-0.9) -- (4.9,-0.1);
    \draw[gray,dashed](3,1) -- (5,-1);
    \draw (-2,-1.8) node{$c_1=(abcde)$} (1,-1.8) node{$c_2=(abcfde)$} (4,-1.8) node{$c_3=(bgcdhe)$};
    \end{tikzpicture}
    \caption{Effect of modified conjugation with $\gamma_i$.}
    \label{fig:effect-conj-var}
\end{figure}
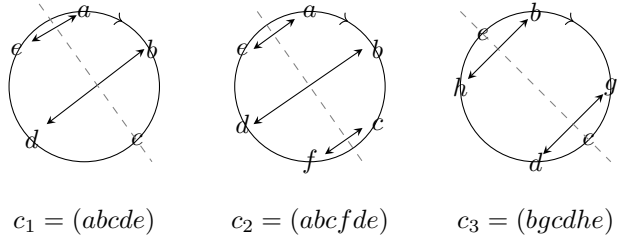
\end{example}
In the following lemma, we summarize the relations
between the action of conjugation with $\gamma$ and the action of conjugation with $\gamma_i$ described
in Example \ref{exa:action-involution-var}.

\begin{lemma}\label{lem:invol-sign-change}
Let $\gamma,\sigma$ be two cycles in $\sal_d$ satisfying $\gamma^2=\id$ and $\gamma\circ\sigma\circ\gamma=\sigma^{-1}$.
Suppose that $\sigma=c_k\circ c_{k-1}\circ\cdots\circ c_1$
is the disjoint cycle decomposition of $\sigma$.
Then
\begin{itemize}
    \item cycles exchanged by the conjugation with $\gamma$
    are also exchanged by
    the conjugation with $\gamma\circ\sigma$;
    \item any cycle $c_i$ which is inverted by the conjugation with $\gamma$
    is also inverted by the conjugation with $\gamma\circ\sigma$.
    \begin{itemize}
        \item if $c_i$ is of odd length, the only element in $c_i$ fixed by the conjugation with
        $\gamma\circ\sigma$
        is the $\lfloor\frac{l(c_i)}{2}\rfloor$-th element after the fixed point of $\gamma$
        according to the orientation of the cycle $c_i$.
        \item if $c_i$ is of even length and there are two fixed
        points of the conjugation with $\gamma$ in $c_i$, the cycle $c_i$ has no fixed point of $\gamma\circ\sigma$. Moreover, the two fixed points of $\gamma$ in $c_i$ are the two first elements in the pair of subsequences consisting of $\frac{l(c_i)}{2}$ consecutive elements in $c_i$ exchanged by the conjugation with $\gamma\circ\sigma$, respectively.
        \item if $c_i$ is of even length and there is no fixed point of $\gamma$ in $c_i$, the cycle $c_i$ contains two elements which are fixed by the conjugation with $\gamma\circ\sigma$. Moreover,
        the two fixed points of $\gamma\circ\sigma$ are the two last elements in the pair of subsequences consisting of $\frac{l(c_i)}{2}$ consecutive elements in $c_i$ exchanged by the conjugation with $\gamma$, respectively.
    \end{itemize}
\end{itemize}
\end{lemma}

\begin{proof}
Since $\gamma\circ\sigma\circ\gamma=\sigma^{-1}$
and $\gamma^2=\id$, it is easy to check that
$(\gamma\circ\sigma)\circ\sigma\circ(\gamma\circ\sigma)=\sigma^{-1}$
and $(\gamma\circ\sigma)^2=\id$.
If $\gamma\circ c_i\circ\gamma=c_j^{-1}$ and
$\gamma\circ c_j\circ\gamma=c_i^{-1}$,
one obtain that $(\gamma\circ\sigma)\circ c_i\circ(\gamma\circ\sigma)=c_j^{-1}$ and
$(\gamma\circ\sigma)\circ c_j\circ(\gamma\circ\sigma)=c_i^{-1}$.
Therefore, cycles exchanged by the conjugation with $\gamma$ are also exchanged by the conjugation with $\gamma\circ\sigma$.
Suppose that $c_i=(a_1,a_2,\ldots,a_k)$,
then $\gamma\circ c_i\circ\gamma=(\gamma(a_1),\gamma(a_2),\ldots,\gamma(a_k))$
and $(\gamma\circ\sigma)\circ c_i\circ(\gamma\circ\sigma)=(\gamma\circ\sigma(a_1),\gamma\circ\sigma(a_2),\ldots,\gamma\circ\sigma(a_k))$.
We consider the case that $c_i$ is a cycle of odd length first.
Without loss of generality, we suppose that $a_1$
is the fixed point of $\gamma$ in $c_i$.
Under the conjugation with $\gamma\circ\sigma$, $a_1$ is
mapped to $\gamma\circ\sigma(a_1)=\gamma(a_2)=a_k$,
so the fixed point of $\gamma\circ\sigma$ in $c_i$
is $a_{1+\lfloor\frac{k}{2}\rfloor}$.
Now we consider the case that $c_i$ is of even length.
Suppose that $a_1$ and $a_{1+\frac{k}{2}}$ are two
fixed points of $\gamma$ in $c_i$.
Under the conjugation with $\gamma\circ\sigma$,
$a_1$ and $a_{1+\frac{k}{2}}$ are mapped to
$\gamma\circ\sigma(a_1)=\gamma(a_2)=a_k$ and
$\gamma\circ\sigma(a_{1+\frac{k}{2}})=\gamma(a_{2+\frac{k}{2}})=a_{\frac{k}{2}-1}$,
so there is no fixed point of $\gamma\circ\sigma$ in $c_i$.
If there is no fixed point of $\gamma$ in $c_i$,
we suppose that $\gamma(a_1)=a_k$.
Under the conjugation with $\gamma\circ\sigma$,
$a_k$ and $a_{\frac{k}{2}}$ are mapped to themselves:
$\gamma\circ\sigma(a_k)=\gamma(a_1)=a_k$ and
$\gamma\circ\sigma(a_{\frac{k}{2}})=\gamma(a_{1+\frac{k}{2}})=a_{\frac{k}{2}}$.
\end{proof}

In \cite[Construction $3.13$]{gpmr-2015},
a tropical cover was constructed for a real
factorization of type $(g,\lambda,\mu;\undl S(r))$.
In the following, we construct a tropical cover
for any real factorization of type $(g,\lambda,\mu;\undl S(s))$, where $0\leq s\leq r$.
\begin{construction}
\label{const2}
Let $\undl S(s)=\{\sk_1,\sk_2,\ldots,\sk_r\}$ be any sequence of
signs with $s$ positive entries,
where $r=l(\lambda)+l(\mu)+2g-2$.
Suppose that $(\gamma,\sigma_1,\tau_1,\ldots,\tau_r,\sigma_2)$ is a real factorization
of type $(g,\lambda,\mu;\undl S(s))$.
We first draw $l(\sigma_1)$ left ends and take the cycle lengths of $\sigma_1$ as the weights of the left ends.
Two ends corresponding to two cycles which are exchanged by $\gamma$
are coloured in dotted.
An end corresponding to an
even cycle with two fixed points is coloured in red,
and ends corresponding to even cycles with no fixed point
are coloured in blue.
Since the composition $\tau_1\circ\sigma_1$ of $\tau_1$ and $\sigma_1$
either cuts a cycle of $\sigma_1$ into two cycles or
joins two cycles of $\sigma_1$,
we either draw a 3-valent vertex with one incoming edge and
two outgoing edges or draw a 3-valent vertex with
two incoming edges and one outgoing edge.
Then we mark the new edges with the lengths
of the new cycles of $\tau_1\circ\sigma_1$.
The effects of the signs $\undl S(s)$ and involution $\gamma$
are shown by the colours of the edges:
\begin{itemize}
    \item if the conjugation with $\gamma_1$ exchanges two
    cycles of $\tau_1\circ\sigma_1$, we draw the corresponding two edges in dotted;
    \item in the case that the conjugation with $\gamma_1$
    inverts a cycle of $\tau_1\circ\sigma_1$, we use colours to distinguish the action of $\gamma_1$ on even cycles:
    \begin{itemize}
        \item if $\gamma_1=\gamma$, we colour even edges
        according to the rule of $\gamma$, that is,
        we draw an edge corresponding to an even cycle with two
        fixed points of $\gamma_1$ in red, and draw an edge
        corresponding to an even cycle with no fixed point in blue;
        \item if $\gamma_1\neq\gamma$, we colour even edges
        according to the rule different from that of $\gamma$, that is,
        we draw an edge
        corresponding to an even cycle with no fixed point
        of $\gamma_1$ in red, and draw the edge corresponding
        to an even cycle with two fixed points in blue.
    \end{itemize}
    \item All other edges are drawn with normal lines.
\end{itemize}
For the compositions $\tau_i\circ\tau_{i-1}\circ\cdots\circ\tau_1\circ\sigma_1$,
$i=2,\cdots,r$, we repeat the above procedure:
\begin{itemize}
    \item draw new edges according to the cut-join operation
    of the cycles, and mark new edges with the lengths of new cycles;
    \item edges corresponding to two cycles of
    $\tau_i\circ\tau_{i-1}\circ\cdots\circ\tau_1\circ\sigma_1$
    which are exchanged by $\gamma_i$ are drawn with dotted lines;
    \item even edges of
    $\tau_i\circ\tau_{i-1}\circ\cdots\circ\tau_1\circ\sigma_1$
    are coloured according to the same rule of $\gamma_{i-1}$,
    if $\gamma_i=\gamma_{i-1}$. Otherwise, we colour even
    edges of $\tau_i\circ\tau_{i-1}\circ\cdots\circ\tau_1\circ\sigma_1$ according to the rule different from that of $\gamma_{i-1}$;
    \item all other edges are drawn with normal lines.
\end{itemize}
\end{construction}

Construction $\ref{const2}$ gives us a way to produce a real
tropical cover from a real
factorization of a particular type.
\begin{lemma}\label{lem:coloured-well}
Let $\undl S(s)$ be a sequence of signs with $s$ positive entries.
For any real factorization of type
$(g,\lambda,\mu;\undl S(s))$, Construction $\ref{const2}$ gives a real tropical cover of type
$(g,\lambda,\mu,\undl x)$ whose positive and
negative branch points produce a splitting
$\undl x=\undl x^+\sqcup\undl x^-$ in accordance
with the sequence of signs $\undl S(s)$,
\textit{i.e.} $x_i\in\undl x^+$ if and only if
$\sk_i=+1$, and $x_j\in\undl x^-$ if and only if $\sk_j=-1$.
\end{lemma}

\begin{proof}
Lemma \ref{lem:invol-sign-change} guarantees that
an even edge is drawn with only one colour in
Construction \ref{const2},
even if we use different involutions according to
the signs $\undl S(s)$ to determine colours.
Construction \ref{const2} produces a coloured
graph $C$ from a real factorization of type
$(g,\lambda,\mu;\undl S(s))$.
To see $\varphi:C\to T\pb^1$ is a real tropical cover of type
$(g,\lambda,\mu,\undl x)$ possessing a specified
splittings, we choose a subset consisting of
all the edges drawn in dotted as the subset $I_\rho$.
Note that Lemma \ref{lem:invol-sign-change}
and the analysis in the proof of
\cite[Lemma $3.12$]{gpmr-2015} imply that all
even edges in a connected component of the subgraph $C\setminus I_\rho^\circ$ are in the same colour,
and the real simple branch point $x_i\in\undl x$
has the same sign with $\sk_i\in\undl S(s)$.
Therefore, the subset $I_\rho$ and the coloured edges of $C$ form a colouring $\rho$ on $\varphi$ such that $(\varphi,\rho)$ is a real tropical cover possessing the required splitting $\undl x=\undl x^+\sqcup\undl x^-$.
\end{proof}

\begin{example}
Let $\sigma_1=(1)(234)$, $\gamma=(24)$, $\tau_1=(34)$, $\tau_2=(13)$, $\sigma_2=(24)(13)$.
Then $(\gamma,\sigma_1,\tau_1,\tau_2,\sigma_2)$ is a real factorization of type $(0,(1,3),(2,2);(+1,-1))$.
The corresponding real tropical cover is depicted in Figure $\ref{fig:coloured-graph}$.
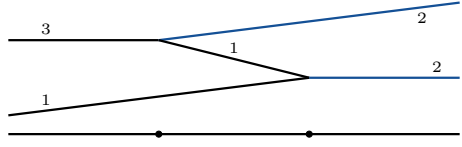
\begin{figure}
    \centering
    \begin{tikzpicture}
    \draw[line width=0.3mm] (-3,0)--(-1,0);
    \draw[line width=0.3mm,darkblue] (-1,0)--(3,0.5);
    \draw[line width=0.3mm] (-1,0)--(1,-0.5);
    \draw[line width=0.3mm,darkblue] (1,-0.5)--(3,-0.5);
    \draw[line width=0.3mm] (-3,-1)--(1,-0.5);
    \draw[line width=0.3mm] (-3,-1.25)--(3,-1.25);
    \foreach \Point in {(-1,-1.25), (1,-1.25)}
    \draw[fill=black] \Point circle (0.04);
    \draw[line width=0.3mm] (-2.5,0.15) node{\tiny{$3$}} (-2.5,-0.8) node{\tiny{$1$}} (2.5,0.3) node{\tiny{$2$}} (2.7,-0.35) node{\tiny{$2$}} (0,-0.1) node{\tiny{$1$}};
    \end{tikzpicture}
    \caption{A real tropical cover of type $(0,(1,3),(2,2);(+1,-1))$.}
    \label{fig:coloured-graph}
\end{figure}
\end{example}

\begin{remark}\label{rem:corr-trans-vertex}
Let $(\varphi:C\to T\pb^1,\rho)$ be the real tropical cover which is obtained from a given real factorization
$(\gamma,\sigma_1,\tau_1,\ldots,\tau_r,\sigma_2)$
of type $(g,\lambda,\mu;\undl S(s))$ by Construction $\ref{const2}$.
The procedure of the Construction $\ref{const2}$ shows that there is a $1$-$1$ correspondence between the transpositions in the real factorization $(\gamma,\sigma_1,\tau_1,\ldots,\tau_r,\sigma_2)$ and the inner vertices of $C$.
Moreover, for any edge in $C$ which is not in a symmetric cycle or a symmetric fork there is exactly one disjoint cycle in $\sigma_1$ or in the compositions $\tau_i\circ\cdots\circ\tau_1\circ\sigma_1$, $i=1,\ldots,r$, corresponds to it.
In the following, we will neither distinguish a transposition in $(\gamma,\sigma_1,\tau_1,\ldots,\tau_r,\sigma_2)$ from the corresponding inner vertex in $C$, nor distinguish an edge which is not in a symmetric cycle or a symmetric fork in $C$ from the corresponding disjoint cycle in $\sigma_1$ or in the compositions $\tau_i\circ\cdots\circ\tau_1\circ\sigma_1$, $i=1,\ldots,r$.
\end{remark}

Let $\sigma,\gamma\in\sal_d$ satisfy $\gamma^2=\id$ and
$\gamma\circ\sigma\circ\gamma=\sigma^{-1}$.
Once a transposition $\tau$ is composited with $\sigma$, it either joins two cycles of $\sigma$ or cuts a cycle of $\sigma$ into two.
When the effect of the involution $\gamma$ is considered,
the number of transpositions satisfying $\gamma\circ(\tau\circ\sigma)\circ\gamma=(\tau\circ\sigma)^{-1}$ was analysed in \cite[Lemma $3.15$]{gpmr-2015}.
In the following, we consider the effect of modified involution $\gamma\circ\sigma$ on the number of transpositions, that is to count the number of transpositions satisfying $(\gamma\circ\sigma)\circ(\tau\circ\sigma)\circ(\gamma\circ\sigma)=(\tau\circ\sigma)^{-1}$.
\begin{lemma}
\label{lem:numb-trans}
Let $\sigma,\gamma\in\sal_d$ satisfy $\gamma^2=\id$ and
$\gamma\circ\sigma\circ\gamma=\sigma^{-1}$.
The number of transpositions $\tau$ satisfying
$(\gamma\circ\sigma)\circ(\tau\circ\sigma)\circ(\gamma\circ\sigma)=(\tau\circ\sigma)^{-1}$
is characterized according to the cycle types
they produced from given cycle types in the following:
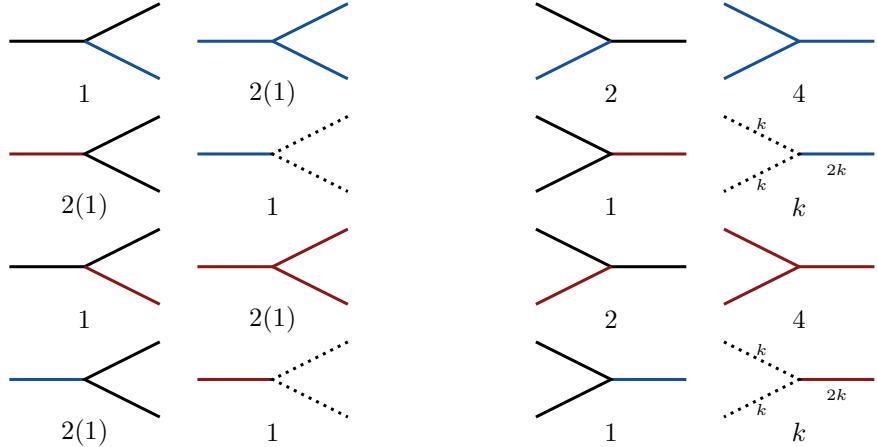
\begin{figure}[ht]
    \centering
    \begin{tikzpicture}
    \draw[line width=0.4mm] (-3,0)--(-2,0)--(-1,0.5);
    \draw[line width=0.4mm,darkblue] (-2,0)--(-1,-0.5);
    \draw (-2,-0.7) node{$1$};
    \draw[line width=0.4mm,darkblue] (-0.5,0)--(0.5,0)--(1.5,0.5);
    \draw[line width=0.4mm,darkblue] (0.5,0)--(1.5,-0.5);
    \draw (0.5,-0.7) node{$2(1)$};
    \draw[line width=0.4mm,cardinalred] (-3,-1.5)--(-2,-1.5);
    \draw[line width=0.4mm] (-1,-1)--(-2,-1.5)--(-1,-2);
    \draw (-2,-2.2) node{$2(1)$};
    \draw[line width=0.4mm,darkblue] (-0.5,-1.5)--(0.5,-1.5);
    \draw[line width=0.4mm,dotted] (1.5,-1)--(0.5,-1.5)--(1.5,-2);
    \draw (0.5,-2.2) node{$1$};
    \draw[line width=0.4mm] (4,0.5)--(5,0)--(6,0);
    \draw[line width=0.4mm,darkblue] (5,0)--(4,-0.5);
    \draw (5,-0.7) node{$2$};
    \draw[line width=0.4mm,darkblue] (6.5,0.5)--(7.5,0)--(8.5,0);
    \draw[line width=0.4mm,darkblue] (7.5,0)--(6.5,-0.5);
    \draw (7.5,-0.7) node{$4$};
    \draw[line width=0.4mm,cardinalred] (6,-1.5)--(5,-1.5);
    \draw[line width=0.4mm] (4,-1)--(5,-1.5)--(4,-2);
    \draw (5,-2.2) node{$1$};
    \draw[line width=0.4mm,darkblue] (8.5,-1.5)--(7.5,-1.5);
    \draw[line width=0.4mm,dotted] (6.5,-1)--(7.5,-1.5)--(6.5,-2);
    \draw (7.5,-2.2) node{$k$} (7,-1.9) node{\tiny{$k$}} (7,-1.1) node{\tiny{$k$}} (8,-1.7) node{\tiny{$2k$}};
    \draw[line width=0.4mm] (-3,-3)--(-2,-3)--(-1,-2.5);
    \draw[line width=0.4mm,cardinalred] (-2,-3)--(-1,-3.5);
    \draw (-2,-3.7) node{$1$};
    \draw[line width=0.4mm,cardinalred] (-0.5,-3)--(0.5,-3)--(1.5,-2.5);
    \draw[line width=0.4mm,cardinalred] (0.5,-3)--(1.5,-3.5);
    \draw (0.5,-3.7) node{$2(1)$};
    \draw[line width=0.4mm,darkblue] (-3,-4.5)--(-2,-4.5);
    \draw[line width=0.4mm] (-1,-4)--(-2,-4.5)--(-1,-5);
    \draw (-2,-5.2) node{$2(1)$};
    \draw[line width=0.4mm,cardinalred] (-0.5,-4.5)--(0.5,-4.5);
    \draw[line width=0.4mm,dotted] (1.5,-4)--(0.5,-4.5)--(1.5,-5);
    \draw (0.5,-5.2) node{$1$};
    \draw[line width=0.4mm] (4,-2.5)--(5,-3)--(6,-3);
    \draw[line width=0.4mm,cardinalred] (5,-3)--(4,-3.5);
    \draw (5,-3.7) node{$2$};
    \draw[line width=0.4mm,cardinalred] (6.5,-2.5)--(7.5,-3)--(8.5,-3);
    \draw[line width=0.4mm,cardinalred] (7.5,-3)--(6.5,-3.5);
    \draw (7.5,-3.7) node{$4$};
    \draw[line width=0.4mm,darkblue] (6,-4.5)--(5,-4.5);
    \draw[line width=0.4mm] (4,-4)--(5,-4.5)--(4,-5);
    \draw (5,-5.2) node{$1$};
    \draw[line width=0.4mm,cardinalred] (8.5,-4.5)--(7.5,-4.5);
    \draw[line width=0.4mm,dotted] (6.5,-4)--(7.5,-4.5)--(6.5,-5);
    \draw (7.5,-5.2) node{$k$} (7,-4.9) node{\tiny{$k$}} (7,-4.1) node{\tiny{$k$}} (8,-4.7) node{\tiny{$2k$}};
    \end{tikzpicture}
    \caption{The cut and join multiplicities: the numbers in brackets are used if the two edges are in a symmetric cycle or a symmetric fork.}
    \label{fig:numb-trans}
\end{figure}
\end{lemma}

\begin{proof}
Note that $\gamma_1:=\gamma\circ\sigma$ is an involution and
$\gamma_1\circ\sigma\circ\gamma_1=\sigma^{-1}$.
Therefore, the number of transpositions satisfying
$\gamma_1\circ(\tau\circ\sigma)\circ\gamma_1=(\tau\circ\sigma)^{-1}$
is the same as the number of transpositions satisfying
$\gamma\circ(\tau\circ\sigma)\circ\gamma=(\tau\circ\sigma)^{-1}$, if they produce same cycle types from given cycle types.
The case that transpositions
satisfy $\gamma\circ(\tau\circ\sigma)\circ\gamma=(\tau\circ\sigma)^{-1}$ is given in \cite[Lemma $3.15$]{gpmr-2015}.
Lemma $\ref{lem:invol-sign-change}$ and Construction $\ref{const2}$ imply that the number of transpositions satisfying
$\gamma_1\circ(\tau\circ\sigma)\circ\gamma_1=(\tau\circ\sigma)^{-1}$ depends only on the cycle types they produced from given cycle types, not dependent on the signs of the vertices corresponding to the transpositions, so we complete our proof.
\end{proof}

\begin{definition}
\label{def:admissible-number}
Let $\sigma$, $\gamma\in\sal_d$ satisfy $\gamma^2=\id$ and $\gamma\circ\sigma\circ\gamma=\sigma^{-1}$.
Let $\tau$ be a transposition satisfying $\gamma\circ(\tau\circ\sigma)\circ\gamma=(\tau\circ\sigma)^{-1}$ or
$(\gamma\circ\sigma)\circ(\tau\circ\sigma)\circ(\gamma\circ\sigma)=(\tau\circ\sigma)^{-1}$.
The entries in $\tau$ are called the \textit{admissible integers} in $\sigma$ and $\tau\circ\sigma$.
\end{definition}

\begin{remark}
Let $\sigma$, $\gamma\in\sal_d$ satisfy $\gamma^2=\id$ and $\gamma\circ\sigma\circ\gamma=\sigma^{-1}$.
If the cycle type of $\tau\circ\sigma$ and the colouring on $\tau\circ\sigma$ are known,
all the admissible integers are determined by Lemma $\ref{lem:numb-trans}$ or \cite[Lemma $3.15$]{gpmr-2015}.
\end{remark}

\begin{lemma}
\label{lem:mult-real-trop}
For a real tropical cover $(\varphi,\rho)$ of type
$(g,\lambda,\mu,\undl x)$ whose real branch points
possess a splitting $\undl x=\undl x^+\sqcup\undl x^-$ with $|\undl x^+|=s$,
there are
$$
d!\cdot\mult^\rb(\varphi,\rho)
$$
real
factorizations of type $(g,\lambda,\mu;\undl S(s)$
that produce the real tropical cover $(\varphi,\rho)$
according to Construction $\ref{const2}$,
where the signs $\undl S(s)$ are determined by the
splitting $\undl x=\undl x^+\sqcup\undl x^-$, that is, $\sk_i=+1$ if and only if $x_i\in\undl x^+$, and $\sk_j=-1$ if and only if $x_j\in\undl x^-$.
\end{lemma}

\begin{proof}
Lemma $\ref{lem:invol-sign-change}$, Construction
$\ref{const2}$, Lemma $\ref{lem:numb-trans}$ and \cite[Lemma $3.15$]{gpmr-2015}
imply that the proof of this Lemma
is the same as that of \cite[Lemma $3.16$]{gpmr-2015},
so we omit it here.
\end{proof}

\begin{proof}[Proof of Theorem $\ref{thm:real-tropical-Hurwitz}$]
By Lemma $\ref{lem:realDH1}$,
$H^\rb_g(\lambda,\mu;s)=|\fl^\rb(g,\lambda,\mu;\undl S(s))|$.
Construction $\ref{const2}$ and Lemma $\ref{lem:mult-real-trop}$ imply that
$|\fl^\rb(g,\lambda,\mu;\undl S(s))|=d!\cdot
\sum_{[(\varphi,\rho)]}\mult^\rb(\varphi,\rho)$,
where we sum over all isomorphism classes $[(\varphi,\rho)]$
of real tropical covers of type $(g,\lambda,\mu,\undl x)$
whose positive and negative branch points reproduce the splitting $\undl x^+,\undl x^-$.
Therefore, we obtain
$H^\rb_g(\lambda,\mu;s)=\sum_{[(\varphi,\rho)]}d!\cdot\mult^\rb(\varphi,\rho)$.
\end{proof}

\section{Lower bounds of real monotone double Hurwitz numbers}
\label{sec:real-mono}
In this section, we first consider the relation between conjugation with an involution and the monotonicity condition of factorizations.
According to the types of the sequences of signs that we used to describe real factorizations,
we introduce two series of numbers which can be considered as two real versions of monotone double Hurwitz numbers.
At last, we show that monotone zigzag numbers and universally monotone zigzag numbers are the lower bounds of these two series of numbers respectively.

\subsection{Monotonicity condition and involution}
In this subsection, we give some examples to illustrate the effect of conjugation with an involution on the monotonicity condition of factorizations,
then we introduce two real versions of monotone double Hurwitz numbers.
Now we recall some basic facts about monotone double Hurwitz numbers from \cite{ggpn-2013,ggpn-2013a,ggpn-2014,hahn-2019,hkl-2018}.

Let $g\geq0$ and $d\geq1$ be two integers.
Suppose that $\lambda$ and $\mu$ are two partitions of $d$.

\begin{definition}
\label{def:monot-fact}
A \textit{factorization} of type $(g,\lambda,\mu)$ is a tuple $(\sigma_1,\tau_1,\ldots,\tau_r,\sigma_2)$ of
$\sal_d$ satisfying the first four conditions in
Definition $\ref{def:real-factor}$.
Suppose that in the tuple $(\sigma_1,\tau_1,\ldots,\tau_r,\sigma_2)$ $\tau_i=(a_i,b_i)$ with
$a_i<b_i$, $i=1,\ldots,r$.
We call the factorization $(\sigma_1,\tau_1,\ldots,\tau_r,\sigma_2)$
a \textit{monotone factorization} of type $(g,\lambda,\mu)$,
if $b_i\leq b_{i+1}$ for $i=1,\ldots,r-1$.
\end{definition}
Let $\vec\fl(g,\lambda,\mu)$ denote the set of
all monotone factorizations of type $(g,\lambda,\mu)$.
Then the number
\begin{equation}\label{eq:monot-Hurwitz}
\vec H^\cb_g(\lambda,\mu)=|\vec\fl(g,\lambda,\mu)|
\end{equation}
is called \textit{the monotone double Hurwitz number} \cite{ggpn-2013,ggpn-2013a,ggpn-2014}.

\begin{remark}\label{rem:mon-def}
We use the convention in \cite{ggpn-2013,ggpn-2013a,ggpn-2014} to denote monotone double Hurwitz numbers,
since the monotonicity condition depends on a total ordering of the sheets of a branched cover.
Note that our convention is different from that in \cite{hahn-2019,hkl-2018}.
The monotone double Hurwitz number considered in \cite{hahn-2019,hkl-2018} is the number $\frac{1}{d!}\vec{H}_g(\lambda,\mu)$ in our notation.
\end{remark}


In the following, we study the effect of
involutions on monotone factorizations and consider
real analogues of monotone double Hurwitz numbers.
The symmetric group approach to real double Hurwitz
numbers (Lemma $\ref{lem:realDH1}$) suggests us some possible ways to define what
a real monotone double Hurwitz number with particular number of positive real branch points is.

\begin{definition}\label{def:signed-mixed-fact}
Let $\undl S(s)=\{\sk_1,\sk_2,\ldots,\sk_r\}$ be a sequence of signs with $s$ positive entries.
Let $(\gamma,\sigma_1,\tau_1,\ldots,\tau_r,\sigma_2)$ be a tuple of type $(g,\lambda,\mu)$ with $\tau_i=(a_i,b_i)$,
where $a_i<b_i$, $i=1,\ldots,r$.
The tuple is a \textit{real monotone factorization} of type $(g,\lambda,\mu;\undl S(s))$,
if it satisfies all the conditions listed in
Definition $\ref{def:real-factor}$ and the monotonicity condition:
$$
b_i\leq b_{i+1}, ~\forall i\in\{1,2,\ldots,r-1\}.
$$
\end{definition}
Let $\vec\fl^\rb(g,\lambda,\mu;\undl S(s))$
denote the set of all real monotone
factorizations of type $(g,\lambda,\mu;\undl S(s))$.
We consider the numbers defined as follows:
\begin{equation}\label{eq:signed-mixed-Hurwitz}
\vec{H}^\rb_g(\lambda,\mu;\undl S(s)):=|\vec\fl^\rb(g,\lambda,\mu;\undl S(s))|.
\end{equation}
The number $\vec{H}^\rb_g(\lambda,\mu;\undl S(s))$ is called the
\textit{real monotone double Hurwitz number} with $s$ positive branch points under the sequence $\undl S(s)$.
In particular, when $\undl S(s)=\{\sk_1,\sk_2,\ldots,\sk_r\}$ is a simple sequence of signs with $s$ positive entries, we use $\vec H^\rb_g(\lambda,\mu;s)$ as the abbreviation of $\vec{H}^\rb_g(\lambda,\mu;\undl S(s))$,
and call the number
$\vec{H}^\rb_g(\lambda,\mu;s)$ the
real monotone double Hurwitz number with $s$ positive branch points under simple splitting.

The effect of an involution on
monotone factorizations is quite different from
that on ordinary factorizations.
The number of real monotone factorizations
$|\vec\fl^\rb(g,\lambda,\mu;\undl S(s))|$ depends
on the sequence of signs $\undl S(s)$
(see Example $\ref{exa:real-simple-monotone}$),
while the number of real
factorizations $|\fl^\rb(g,\lambda,\mu;\undl S(s))|$
depends only on the number $s$ of positive entries in
$\undl S(s)$, and not on the sequence of signs $\undl S(s)$.

\begin{example}
\label{exa:real-simple-monotone}
Suppose that $d=3$, $g=0$ and $\lambda=\mu=(1^3)$.
From \cite[Theorem $1.1$]{ggpn-2013}, we obtain that $\vec H_0((1^3),(1^3))=8$.
We list the tuples of monotone factorizations of
type $(0,(1^3),(1^3))$ in Table $\ref{tab:simple-mon}$.
\begin{table}[ht]
    \centering
    \begin{tabular}{|c|c|}
    \hline
     $(12),(13),(23),(13)$ & $(12),(12),(13),(13)$ \\
     \hline
     $(12),(23),(13),(23)$  & $(12),(12),(23),(23)$\\
     \hline
     $(13),(23),(23),(13)$ & $(13),(13),(23),(23)$ \\
     \hline
     $(23),(13),(13),(23)$ & $(23),(23),(13),(13)$\\
     \hline
    \end{tabular}
    \caption{Factorizations of type $(0,(1^3),(1^3))$.}
    \label{tab:simple-mon}
\end{table}
The number of real monotone factorizations of type
$(0,(1^3),(1^3);(+1,+1,+1,-1))$ is $6$,
and these factorizations are listed in Table $\ref{tab:signed-simple-mon1}$.
\begin{table}[ht]
    \centering
    \begin{tabular}{|c|c|}
     \hline
     $(12),(12),(13),(13),\gamma=\id$ & $(12),(12),(23),(23),\gamma=\id$\\
     \hline
     $(13),(23),(23),(13),\gamma=(13)$ & $(13),(13),(23),(23),\gamma=\id$ \\
     \hline
     $(23),(13),(13),(23),\gamma=(23)$ & $(23),(23),(13),(13),\gamma=\id$\\
     \hline
    \end{tabular}
    \caption{Real factorizations of type $(0,(1^3),(1^3);(+1,+1,+1,-1))$.}
    \label{tab:signed-simple-mon1}
\end{table}
The number of real monotone factorizations of type
$(0,(1^3),(1^3);(+1,-1,+1,+1))$ is $4$.
We list these factorizations in Table $\ref{tab:signed-simple-mon2}$.
\begin{table}[ht]
    \centering
    \begin{tabular}{|c|c|}
     \hline
     $(12),(12),(13),(13),\gamma=(12)$ & $(12),(12),(23),(23),\gamma=(12)$\\
     \hline
     $(13),(13),(23),(23),\gamma=(13)$ & $(23),(23),(13),(13),\gamma=(23)$ \\
     \hline
    \end{tabular}
    \caption{Real factorizations of type $(0,(1^3),(1^3);(+1,-1,+1,+1))$.}
    \label{tab:signed-simple-mon2}
\end{table}
\end{example}

\begin{remark}\label{rem:rmf-1}
It follows from Example $\ref{exa:real-simple-monotone}$ that we have to consider the influence of the sequences of signs $\undl S(s)$ when we try to find some real analogues of monotone double Hurwitz numbers.
\end{remark}

Let $\sigma$ be a permutation of cycle type $\lambda$.
Note that the number of monotone factorizations $(\sigma_1,\tau_1,\ldots,\tau_r,\sigma_2)$
of type $(g,\lambda,\mu)$ with a fixed
$\sigma_1=\sigma$ and the number of real factorizations
$(\gamma,\sigma_1,\tau_1,\ldots,\tau_r,\sigma_2)$
of type $(g,\lambda,\mu;\undl S(s))$ with a fixed
$\sigma_1=\sigma$ do not depend on the choice of $\sigma$,
and they only depend on the cycle type $\lambda$ of $\sigma$
(c.f. \cite[Lemma $3.16$]{gpmr-2015}, \cite[Corollary $5.9$]{mr-2015}, \cite[Section $1.5$]{ggpn-2013}, \cite[Lemma $7$]{dk-2017}, \cite[Lemma $3.3$]{hahn-2019}).
However, the number of real monotone factorizations
$(\gamma,\sigma_1,\tau_1,\ldots,\tau_r,\sigma_2)$
of type $(g,\lambda,\mu;\undl S(s))$ with a fixed
$\sigma_1=\sigma$ does depend on the choice of $\sigma$
(see Example $\ref{exa:real-double-monotone}$).
\begin{example}
\label{exa:real-double-monotone}
The number of real factorizations of type $(0,(1,3),(2,2);(+1,+1))$ is $24$.
There are $8$ permutations of type $(1,3)$.
The number of real factorizations $(\gamma,\sigma_1,\tau_1,\tau_2,\sigma_2)$ of type $(0,(1,3),(2,2);(+1,+1))$ with a fixed $\sigma$ is $3$.
When $\sigma_1=(1)(234)$ or $\sigma_1=(4)(132)$, the real factorizations $(\gamma,\sigma_1,\tau_1,\tau_2,\sigma_2)$ of type $(0,(1,3),(2,2);(+1,+1))$ with a fixed $\sigma_1$ are listed in Table $\ref{tab:num-mon-fact1}$ and Table $\ref{tab:num-mon-fact2}$, respectively.
\begin{table}[ht]
    \centering
    \begin{tabular}{|c|c|c|c|}
     \hline
     $\sigma_1$ & $\gamma$ & $\tau_1$ & $\tau_2$\\
     \hline
     $(1)(234)$ & $(24)$ & $(34)$ &$(13)$ \\
     \hline
     $(1)(234)$ & $(34)$ & $(23)$ &$(12)$ \\
     \hline
     $(1)(234)$ & $(23)$ & $(24)$ &$(14)$ \\
     \hline
    \end{tabular}
    \caption{Real factorizations of type $(0,(1,3),(2,2);(+1,+1))$ with $\sigma_1=(1)(234)$.}
    \label{tab:num-mon-fact1}
\end{table}
\begin{table}[ht]
    \centering
    \begin{tabular}{|c|c|c|c|}
     \hline
     $\sigma_1$ & $\gamma$ & $\tau_1$ & $\tau_2$\\
     \hline
     $(4)(132)$ & $(13)$ & $(12)$ &$(24)$ \\
     \hline
     $(4)(132)$ & $(12)$ & $(23)$ &$(34)$ \\
     \hline
     $(4)(132)$ & $(23)$ & $(13)$ &$(14)$ \\
     \hline
    \end{tabular}
    \caption{Real factorizations of type $(0,(1,3),(2,2);(+1,+1))$ with $\sigma_1=(4)(132)$.}
    \label{tab:num-mon-fact2}
\end{table}
From the listed tables, we obtain that the number of real monotone factorizations of type $(0,(1,3),(2,2);(+1,+1))$ with $\sigma_1=(1)(234)$ is $1$,
and the number of real monotone factorizations with $\sigma_1=(4)(132)$ is $3$.
\end{example}

\begin{remark}\label{rem:rmf-2}
The significant difference between the number of real monotone factorizations with a fixed starting permutation and the numbers of real factorizations or monotone factorizations with a fixed starting permutation shows that
it is quite difficult to use tropical covers to investigate the properties of real analogies of monotone double Hurwitz numbers.
\end{remark}

Based on Remark $\ref{rem:rmf-1}$, Remark $\ref{rem:rmf-2}$ and Lemma $\ref{lem:realDH1}$,
we consider two series of numbers which are defined as follows:
\begin{equation}
\label{eq:def-RMH}
\begin{aligned}
\vec H^\rb_g(\lambda,\mu)&:=\inf_{0\leq s\leq r}\vec H^\rb_g(\lambda,\mu;s),\\
\vec\hl^{\rb}_g(\lambda,\mu)&:=\inf_{\undl S(s)}\vec H^\rb_g(\lambda,\mu;\undl S(s)).
\end{aligned}
\end{equation}
From \cite[Remark $1.2$]{ggpn-2013},
both $\vec H^\rb_g(\lambda,\mu)$ and $\vec\hl^{\rb}_g(\lambda,\mu)$
``count" the number of real ramified covers with labelled sheets satisfying certain ramification conditions and monotonicity conditions.
We call $\vec H^\rb_g(\lambda,\mu)$ (resp. $\vec\hl^{\rb}_g(\lambda,\mu)$) the \textit{real monotone double Hurwitz numbers relative to simple (resp. arbitrary) splittings}.
In this paper,
we are interested in investigating the asymptotic growth of the real monotone double Hurwitz numbers $\vec H^\rb_g(\lambda,\mu)$ and $\vec\hl^{\rb}_g(\lambda,\mu)$ when only simple branch points are added as the degree $d$ goes to infinity.

\subsection{Monotone zigzag numbers}
In this subsection, we characterize a type of tropical covers from zigzag covers introduced in \cite{rau2019}.
First, we recall some notations from \cite[Section $5$]{rau2019}.

Suppose that $\lambda=(\lambda_1,\ldots,\lambda_n)$ is a partition. Denote by
$\lambda^2=(\lambda_1,\lambda_1,\lambda_2,\lambda_2,\ldots,\lambda_n,\lambda_n)$.
Then any partition $\lambda$ can be uniquely decomposed into:
\begin{equation}
\label{eq:tail-decom}
\lambda=(\lambda_e,\lambda_{o,o}^2,\lambda_o),
\end{equation}
where
\begin{itemize}
    \item all parts in $\lambda_e$ are even integers;
    \item all parts in $\lambda_{o,o}$ are odd integers;
    \item all parts in $\lambda_o$ are odd integers and $\lambda_o$ does not contain repeated entries.
\end{itemize}
The decomposition $(\ref{eq:tail-decom})$ is called the \textit{tail decomposition} of $\lambda$.

\begin{remark}
The tail decomposition of a partition here is slightly different from that in \cite[Section $5$]{rau2019}.
In this paper, we do not need to distinguish the weighted $2e$ tails from the weighted $2o$ tails,
so we simplify the notation in tail decomposition.
\end{remark}

A {\it string} $S$ in a tropical curve $C$ is a
connected subgraph such that $S\cap C^\circ$ is
a closed submanifold of $C^\circ$.
We denote by $\sym(C)$ the set of symmetric circles
and symmetric odd forks of $C$.

\begin{definition}[{\cite[Definition 4.4]{rau2019}}]
\label{def:zigzag}
A \textit{zigzag cover} is a tropical cover $\varphi:C\to T\pb^1$
if there is a subset $S\subset C\setminus\sym(\varphi)$
satisfying
\begin{itemize}
    \item $S$ is either a string of odd edges or consists
    of a single inner vertex;
    \item the connected components of $C\setminus S$ are
    of the type depicted in Figure $\ref{fig:zigzag}$. In Figure $\ref{fig:zigzag}$, all the cycles and forks
    are symmetric and of odd weight.
\end{itemize}
\end{definition}
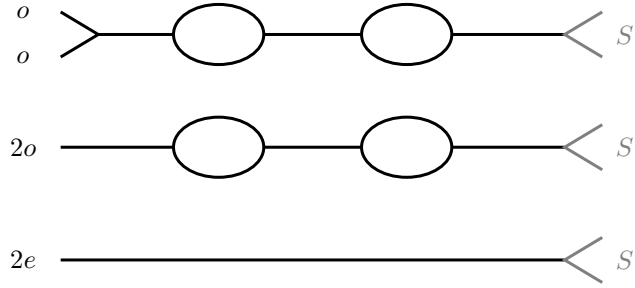
\begin{figure}[ht]
    \centering
    \begin{tikzpicture}
    \draw[line width=0.4mm] (-3,0)--(-1.5,0);
    \draw[line width=0.4mm] (-0.3,0) arc[start angle=0, end angle=360, x radius=0.6, y radius=0.4];
    \draw[line width=0.4mm] (-0.3,0)--(1,0);
    \draw[line width=0.4mm] (2.2,0) arc[start angle=0, end angle=360, x radius=0.6, y radius=0.4];
    \draw[line width=0.4mm] (2.2,0)--(3.7,0);
    \draw[line width=0.4mm,gray]  (4.2,0.3)--(3.7, 0) -- (4.2,-0.3);
    \draw[line width=0.4mm,gray] (4.5,0) node{$S$};
    \draw (-3.5,0) node{$2o$};
    \draw[line width=0.4mm]  (-3,1.8)--(-2.5, 1.5) -- (-3,1.2);
    \draw[line width=0.4mm] (-2.5,1.5)--(-1.5,1.5);
    \draw[line width=0.4mm] (-0.3,1.5) arc[start angle=0, end angle=360, x radius=0.6, y radius=0.4];
    \draw[line width=0.4mm] (-0.3,1.5)--(1,1.5);
    \draw[line width=0.4mm] (2.2,1.5) arc[start angle=0, end angle=360, x radius=0.6, y radius=0.4];
    \draw[line width=0.4mm] (2.2,1.5)--(3.7,1.5);
    \draw[line width=0.4mm,gray]  (4.2,1.8)--(3.7, 1.5) -- (4.2,1.2);
    \draw[line width=0.4mm] (-3.5,1.8) node{$o$};
    \draw[line width=0.4mm] (-3.5,1.2) node{$o$};
    \draw[line width=0.4mm,gray] (4.5,1.5) node{$S$};
    \draw[line width=0.4mm]  (-3,-1.5)--(3.7, -1.5);
    \draw[line width=0.4mm,gray]  (4.2,-1.8)--(3.7, -1.5) -- (4.2,-1.2);
    \draw[line width=0.4mm] (-3.5,-1.5) node{$2e$};
    \draw[line width=0.4mm,gray] (4.5,-1.5) node{$S$};
    \end{tikzpicture}
    \caption{Tails for zigzag covers. It does not matter
    whether $S$ turns or not here. The number of cycles
    in the first two types can be arbitrary.}
    \label{fig:zigzag}
\end{figure}

In the following, we collect some useful properties of zigzag covers from \cite{rau2019}.

\begin{proposition}[{\cite[Proposition $5.2$]{rau2019}}]
\label{prop:rau-zigzag1}
Let $\lambda,\mu$ be two partitions of $d$. If $l(\lambda_o,\mu_o)\leq2$ and $l(\lambda_{o,o},\mu_{o,o})\geq1$, there exist zigzag covers of type $(g,\lambda,\mu,\undl x)$.
\end{proposition}

\begin{proposition}[{\cite[Proposition $4.8$]{rau2019}}]
\label{prop:rau-zigzag2}
Let $\varphi:C\to T P^1$ be a zigzag cover simply branched at $\undl x$.
Choose an arbitrary
splitting $\undl x=\undl x^+\sqcup\undl x^-$ into positive and negative branch points.
Then there exists
a unique colouring $\rho$ of $\varphi$ such that the real tropical cover $(\varphi,\rho)$ has positive and
negative branch points as required.
\end{proposition}

Let $\varphi:C\to T\pb^1$ be a zigzag cover such that
there is a string $S$ of odd edges contained in $C$.
Before introducing the definition of monotone zigzag
covers, we give some useful notations.
An incoming tail depicted in Figure $\ref{fig:zigzag}$ is called an \textit{in-tail},
and an outgoing tail depicted in Figure $\ref{fig:zigzag}$ is called an \textit{out-tail}.
For a tail $t$ attached to the string $S$,
the weight of the edge in $t$
adjacent to the string $S$ is called the \textit{weight} of the tail $t$.
A vertex in $S$ at which $S$ bends is called a \textit{bent vertex}.
The vertices in $S$ other than bent vertices are called \textit{unbent vertices}.
Let $S_1,\ldots,S_n$ be pieces of $S$ obtained by
cutting $S$ at the bent vertices.
Assume that the pieces are arranged as follows:
the two ends of the string $S$ are in $S_1$ and $S_n$,
and $S_i$ intersects with $S_{i+1}$ at a bent vertex for $i=1,\ldots,n-1$.
If the piece $S_1$ contains an in-end of the string $S$,
we call the pieces $S_1,S_3,S_5,\ldots$ the \textit{in-pieces},
and call the remaining pieces $S_2,S_4,\ldots$ the \textit{out-pieces}.
If the piece $S_1$ contains an out-end of the string $S$,
we call the pieces $S_1,S_3,S_5,\ldots$ the \textit{out-pieces},
and call the remaining pieces $S_2,S_4,\ldots$ the \textit{in-pieces}.
The tails attached to bent (resp. unbent) vertices are called \textit{bent tails} (resp. \textit{unbent tails}).
We call an in-piece and all tails attached to it an \textit{in-component},
and call an out-piece and the unbent tails attached to it an \textit{out-component}.
All the bent vertices are considered as inner vertices in the corresponding in-components.

\begin{definition}\label{def:mono-zigzag}
A zigzag cover $\varphi:C\to T\pb^1$ is called
\textit{monotone} if there is a string $S\subset C$ such
that the map $\varphi$ satisfies the following conditions:
\begin{enumerate}
    \item All the unbent tails attached to in-pieces are out-tails,
    and all the unbent tails attached to out-pieces are in-tails.
    Moreover, no unbent tail contains symmetric cycles, and no unbent out-tail contains symmetric forks.
    \item Any bent tail with symmetric cycle or symmetric fork is of weight $2$.
    \item For any tail $t$, let $a_t$
    (resp. $b_t$) be the smallest
    (resp. largest) point of the images of the inner vertices in $t$
    under the map $\varphi$,
    then $[a_t,b_t]\cap[a_{t'},b_{t'}]=\emptyset$
    for any two different tails $t,t'$.
    For any component $C_i$, let $f_i$
    (resp. $g_i$) be the smallest (resp. largest) point of the images of the vertices in $C_i$ under $\varphi$,
    then $[f_i,g_i]\cap[f_j,g_j]=\emptyset$ for any two different components $C_i,C_j$.
\end{enumerate}
\end{definition}

\begin{remark}
The third condition of Definition $\ref{def:mono-zigzag}$
implies that the images of the inner vertices of a component under $\varphi$ form a consecutive sequence of points.
\end{remark}

\begin{example}
We give a monotone zigzag cover $\varphi:C\to T\pb^1$ of certain type in Figure $\ref{fig:mon-zigzag}$.
There is only one inner vertex in the out-components. More exactly, only the out-component containing $S_2$ has an inner vertex.
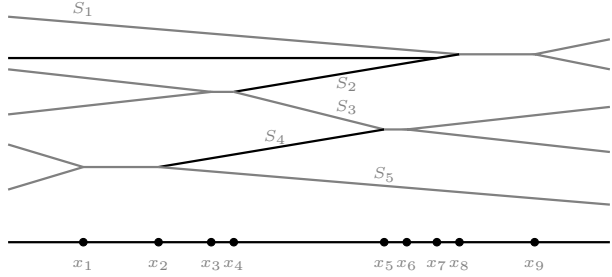
\begin{figure}[ht]
    \centering
    \begin{tikzpicture}
    \draw[line width=0.3mm,gray] (-4,-0.5)--(2,-1);
    \draw[line width=0.3mm] (2,-1)--(-1,-1.5);
    \draw[line width=0.3mm,gray] (-1,-1.5)--(1,-2);
    \draw[line width=0.3mm] (1,-2)--(-2,-2.5);
    \draw[line width=0.3mm,gray] (-2,-2.5)--(4,-3);
    \draw[line width=0.3mm,gray] (2,-1)--(3,-1);
    \draw[line width=0.3mm,gray] (4,-0.8)--(3,-1)--(4,-1.2);
    \draw[line width=0.3mm] (-4,-1.05)--(1.7,-1.05);
    \draw[line width=0.3mm,gray] (-1,-1.5)--(-1.3,-1.5);
    \draw[line width=0.3mm,gray] (-4,-1.8)--(-1.3,-1.5)--(-4,-1.2);
    \draw[line width=0.3mm,gray] (1,-2)--(1.3,-2);
    \draw[line width=0.3mm,gray] (4,-1.7)--(1.3,-2)--(4,-2.3);
    \draw[line width=0.3mm,gray] (-2,-2.5)--(-3,-2.5);
    \draw[line width=0.3mm,gray] (-4,-2.8)--(-3,-2.5)--(-4,-2.2);
    \draw[line width=0.3mm,gray] (-3,-0.4) node{\tiny{$S_1$}} (0.5,-1.4) node{\tiny{$S_2$}} (0.5,-1.7) node{\tiny{$S_3$}} (-0.45,-2.1) node{\tiny{$S_4$}} (1,-2.6) node{\tiny{$S_5$}};
    \draw[line width=0.3mm] (-4,-3.5)--(4,-3.5);
    \foreach \Point in {(-3,-3.5), (-2,-3.5),(-1.3,-3.5),(-1,-3.5),(1,-3.5),(1.3,-3.5),(1.7,-3.5),(2,-3.5),(3,-3.5)}
    \draw[fill=black] \Point circle (0.05);
    \draw[line width=0.4mm,gray] (-3,-3.8) node{\tiny{$x_1$}} (-2,-3.8) node{\tiny{$x_2$}} (-1.3,-3.8) node{\tiny{$x_3$}} (-1,-3.8) node{\tiny{$x_4$}} (1,-3.8) node{\tiny{$x_5$}} (1.3,-3.8) node{\tiny{$x_6$}} (1.7,-3.8) node{\tiny{$x_7$}}(2,-3.8) node{\tiny{$x_8$}} (3,-3.8) node{\tiny{$x_9$}};
    \end{tikzpicture}
    \caption{A monotone zigzag cover of certain type. In-components are in gray, and out-components are in black.}
    \label{fig:mon-zigzag}
\end{figure}
\end{example}

Let $\varphi:C\to T\pb^1$ be a monotone zigzag cover of type $(g,\lambda,\mu,\undl x)$.
Suppose that $\undl x^+\sqcup\undl x^-$ is a splitting of $\undl x$ into positive and negative branch points with $|\undl x^+|=s$, where $s\in\{0,\ldots,r\}$.
Denote by $\rho$ the unique colouring of the monotone zigzag cover $\varphi$
determined by the splitting $\undl x=\undl x^+\sqcup\undl x^-$ according to Proposition $\ref{prop:rau-zigzag2}$.
The colouring $\rho$ and the monotone zigzag
cover $\varphi$ form a real tropical cover with
$s$ real positive branch points.
It follows from Lemma $\ref{lem:mult-real-trop}$ that
there are $d!\cdot\mult^\rb(\varphi,\rho)$ real
factorizations of type $(g,\lambda,\mu;\undl S(s))$
producing $(\varphi,\rho)$, where $\undl S(s)$ is a sequence of signs corresponding to the splitting $\undl x=\undl x^+\sqcup\undl x^-$.
The following Lemma gives a sufficient condition to find a real monotone factorization in the $d!\cdot\mult^\rb(\varphi,\rho)$ real
factorizations producing $(\varphi,\rho)$.

\begin{lemma}\label{lem:mono-fac1}
Let $\varphi:C\to T\pb^1$ be a monotone zigzag cover of type $(g,\lambda,\mu,\undl x)$,
and let $\rho$ be the unique colouring of $\varphi$ determined by a splitting $\undl x=\undl x^+\sqcup\undl x^-$ of the branch points.
Let $(\gamma,\sigma_1,\tau_1,\ldots,\tau_r,\sigma_2)$ be a real factorization associated to $(\varphi,\rho)$ according to Construction $\ref{const2}$, where $\tau_i=(a_i,b_i)$ with $a_i<b_i$ for any $i\in\{1,\ldots,r\}$.
Assume that the real factorization
$(\gamma,\sigma_1,\tau_1,\ldots,\tau_r,\sigma_2)$
satisfies the following condition:
\begin{itemize}
    \item[$(*)$] For any $i\in\{2,\ldots,r\}$, if $\{a_i,b_i\}\subset\{a_1,b_1,\ldots,a_{i-1},b_{i-1}\}$ and $b_i\in\{a_j,b_j\}$ for some $j<i$, we have $b_j=\cdots=b_i$.
\end{itemize}
Then there is at least one real monotone factorization in the coordinate-wise conjugation class of $(\gamma,\sigma_1,\tau_1,\ldots,\tau_r,\sigma_2)$.
\end{lemma}

\begin{proof}
Let $(\gamma,\sigma_1,\tau_1,\ldots,\tau_r,\sigma_2)$ be a real factorization satisfying the condition $(*)$.
Let $\sigma_1=c_{l(\sigma_1)}\circ\cdots\circ c_1$ be the disjoint cycle decomposition of $\sigma_1$.
From the transpositions in the real factorization $(\gamma,\sigma_1,\tau_1,\ldots,\tau_r,\sigma_2)$,
we obtain a sequence of integers $b_1,\ldots,b_r$
such that if
$b_{i}\neq b_{i-1}$, $b_{i}\notin\{a_1,b_1,\ldots,a_{i-1},b_{i-1}\}$.
Up to coordinate-wise conjugation,
we obtain a real factorization
$(\bar\gamma,\bar\sigma_1,\bar\tau_1,\ldots,\bar\tau_r,\bar\sigma_2)$ of the same type with $\bar\tau_i=(\bar a_i,\bar b_i)$,
where $\bar a_i<\bar b_i$ and $\bar b_1\leq\bar b_2\leq\cdots\leq\bar b_r$.
The coordinate-wise conjugation is taken as follows.
Suppose that the sequence $b_1,b_2,\ldots,b_r$
contains $k$ distinct integers:
$$
b_1=\cdots=b_{i_1}\neq b_{i_1+1}=\cdots=b_{i_2}\neq\cdots\neq b_{i_{k-1}+1}=\cdots=b_{i_k},
$$
where $i_1<i_2<\ldots<i_k=r$.
Note that condition $(*)$ implies that $k<d$.
We rearrange the $k$ integers $b_{i_1},\ldots,b_{i_k}$
in an ascending order, that is $b_{i_1}^{(1)}<\cdots<b_{i_k}^{(k)}$.
First, the factorization $(\gamma,\sigma_1,\tau_1,\ldots,\tau_r,\sigma_2)$
is conjugated by $(b_{i_k},b_{i_k}^{(k)})$.
It is easy to see that $a_j<b_{i_k}^{(k)}$ for
any $j\in\{i_{k-1}+1,\ldots,i_k\}$.
Next, let the new factorization be conjugated
by $(b_{i_{k-1}},b_{i_{k-1}}^{(k-1)})$.
If $a_j<b_{i_{k-1}}^{(k-1)}$ for any
$j\in\{i_{k-2}+1,\ldots,i_{k-1}\}$,
we move to the next step that is
to make the achieved factorization be conjugated by
$(b_{i_{k-2}},b_{i_{k-2}}^{(k-2)})$.
Otherwise, we make the achieved factorization
be conjugated by $(a_{i_{k-1}}^*,b_{i_{k-1}}^{(k-1)})$,
where $a_{i_{k-1}}^*=\max\{a_{i_{k-2}+1},\ldots,a_{i_{k-1}}\}$.
From condition $(*)$, one obtains that $a_{i_{k-1}}^*<b_{i_{k}}^{(k)}$.
By repeating the above procedure,
we finally obtain a real monotone factorization of the same type.
\end{proof}

\begin{lemma}\label{lem:mono-zigzag1}
Assume that $\varphi:C\to T\pb^1$ is a monotone zigzag cover simply branched at $\undl x$.
Let $\rho_s$ be the unique colouring of $\varphi$
determined by the simple splitting $\undl x=\undl x^+\sqcup\undl x^-$ with $|\undl x^+|=s$.
Then for any $s\in\{0,1,\ldots,r\}$,
the number of real monotone factorizations
producing $(\varphi,\rho_s)$ according to Construction $\ref{const2}$ is non-zero.
\end{lemma}

\begin{proof}
The colouring $\rho_s$ and the monotone zigzag
cover $\varphi$ form a real tropical cover with
$s$ real positive branch points.
Suppose that $(\gamma,\sigma_1,\tau_1,\ldots,\tau_r,\sigma_2)$
is a real factorization producing the real tropical cover $(\varphi,\rho_s)$
according to Construction $\ref{const2}$.
It follows from Lemma $\ref{lem:mono-fac1}$ that if the factorization $(\gamma,\sigma_1,\tau_1,\ldots,\tau_r,\sigma_2)$ satisfies condition $(*)$,
there is a real monotone factorization in the coordinate-wise conjugation class of $(\gamma,\sigma_1,\tau_1,\ldots,\tau_r,\sigma_2)$ which produces $(\varphi,\rho_s)$ according to Construction $\ref{const2}$.

In the following, we show that in the $d!\cdot\mult^\rb(\varphi,\rho_s)$ real factorizations that produce $(\varphi,\rho_s)$ according to Construction $\ref{const2}$,
there is a factorization satisfying the condition $(*)$ in the Lemma $\ref{lem:mono-fac1}$.
For any zigzag cover $\varphi:C\to T\pb^1$, there are $6$ types inner vertices, up to reflection along the horizontal line, in the tropical curve $C$, and they are depicted in Figure $\ref{fig:vertices-in-MZ}$.
\begin{figure}[ht]
    \centering
    \begin{tikzpicture}
    \draw[line width=0.4mm] (-3,0.5)--(-2,0)--(-3,-0.5);
    \draw[line width=0.4mm] (-2,0)--(-1,0);
    \draw (-2.7,0.5) node{$o$} (-2.7,-0.5) node{$o$} (-1.2,0.2) node{$e$} (-2,-0.7) node{$(1)$};
    \draw[line width=0.4mm] (0,0)--(1,0)--(2,0.5);
    \draw[line width=0.4mm] (1,0)--(2,-0.5);
    \draw (0.2,0.2) node{$e$} (1.7,-0.5) node{$o$} (1.7,0.5) node{$o$} (1,-0.7) node{$(2)$};
    \draw[line width=0.4mm] (3,0)--(4,0);
    \draw[line width=0.4mm,gray] (5,0.5)--(4,0)--(5,-0.5);
    \draw (3.2,0.2) node{$e$} (4.7,-0.5) node{$o_2$} (4.6,0.5) node{$o_1$} (4,-0.7) node{$(3)$};
    \draw[gray] (5,0.3) node{$S$};
    \draw[line width=0.4mm,gray] (-3,-1.5)--(-2,-2)--(-1,-2);
    \draw[line width=0.4mm] (-3,-2.5)--(-2,-2);
    \draw (-2.7,-2.5) node{$e$} (-2.7,-1.5) node{$o_1$} (-1.2,-1.8) node{$o_2$} (-2,-2.7) node{$(4)$};
    \draw[gray] (-1.2,-2.2) node{$S$};
    \draw[line width=0.4mm,gray] (0,-1.5)--(1,-2)--(0,-2.5);
    \draw[line width=0.4mm] (2,-2)--(1,-2);
    \draw (0.3,-2.5) node{$o_2$} (0.3,-1.5) node{$o_1$} (1.8,-1.8) node{$e$} (1,-2.7) node{$(5)$};
    \draw[gray] (1.8,-2.2) node{$S$};
    \draw[line width=0.4mm,gray] (3,-2)--(4,-2)--(5,-1.5);
    \draw[line width=0.4mm] (4,-2)--(5,-2.5);
    \draw (3.2,-1.8) node{$o_1$} (4.7,-2.5) node{$e$} (4.6,-1.5) node{$o_2$} (4,-2.7) node{$(6)$};
    \draw[gray] (4,-1.7) node{$S$};
    \end{tikzpicture}
    \caption{Inner vertices in zigzag cover $\varphi:C\to T\pb^1$; edges in the string $S$ are drawn in gray.}
    \label{fig:vertices-in-MZ}
\end{figure}
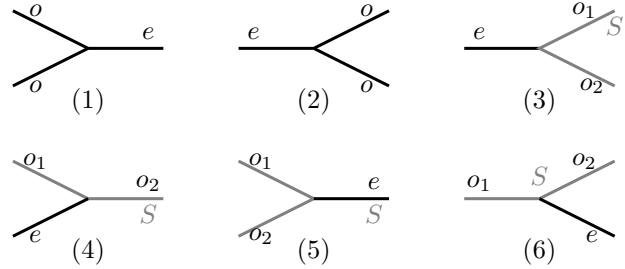
The corresponding real cut and join operations on cycles in $\sal_d$
are depicted in Figure $\ref{fig:cycle-cut-join}$.
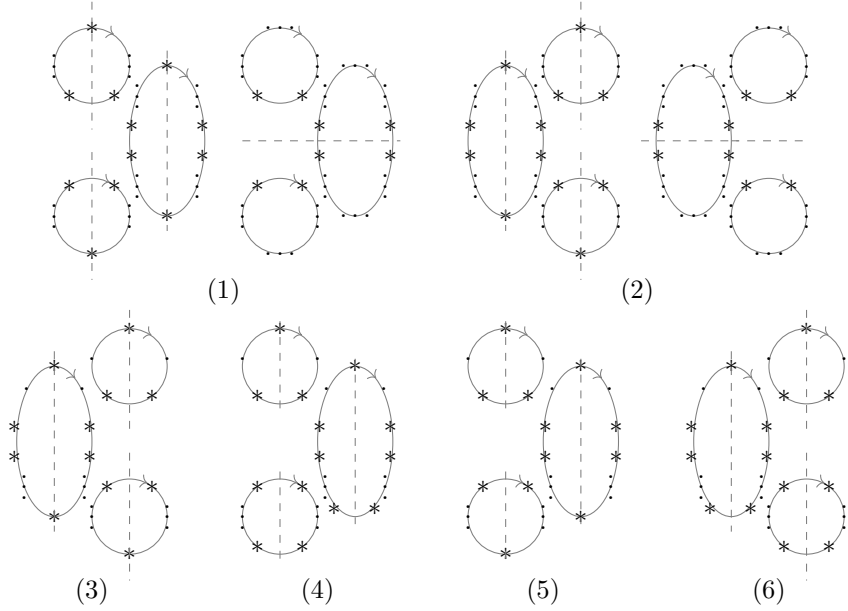
\begin{figure}[ht]
    \centering
    \begin{tikzpicture}
    \draw[decoration={markings, mark=at position 0.18 with {\arrow{<}}},
        postaction={decorate},gray
        ] (-1.5,0) circle (0.5);
    \draw (-1.5,0.5) node {$*$} (-2,0.1) node{$\vdots$} (-1,0.1) node{$\vdots$} (-1.8,-0.4) node{$*$} (-1.2,-0.4) node{$*$};
    \draw[gray,dashed](-1.5,0.85) -- (-1.5,-0.85);
    \draw[decoration={markings, mark=at position 0.2 with {\arrow{<}}},
        postaction={decorate},gray
        ] (-1.5,-2) circle (0.5);
    \draw (-1.5,-2.5) node {$*$} (-2,-1.9) node{$\vdots$} (-1,-1.9) node{$\vdots$} (-1.8,-1.6) node{$*$} (-1.2,-1.6) node{$*$};
    \draw[gray,dashed](-1.5,-1.15) -- (-1.5,-2.85);
    \draw[decoration={markings, mark=at position 0.2 with {\arrow{<}}},  postaction={decorate},gray
        ] (-0.5,-1) ellipse (0.5 and 1);
    \draw (-0.5,0) node {$*$} (-0.5,-2) node {$*$} (-0.9,-0.3) node{$\vdots$} (-0.1,-0.3) node{$\vdots$} (-0.97,-0.8) node{$*$} (-0.03,-0.8) node{$*$} (-0.9,-1.5) node{$\vdots$} (-0.1,-1.5) node{$\vdots$} (-0.97,-1.2) node{$*$} (-0.03,-1.2) node{$*$};
    \draw[gray,dashed](-0.5,0.2) -- (-0.5,-2.2);
    \draw[decoration={markings, mark=at position 0.18 with {\arrow{<}}},
        postaction={decorate},gray
        ] (1,0) circle (0.5);
    \draw (1,0.5) node {$\ldots$} (0.5,0.1) node{$\vdots$} (1.5,0.1) node{$\vdots$} (0.7,-0.4) node{$*$} (1.3,-0.4) node{$*$};
    \draw[gray,dashed](0.5,-1) -- (1.5,-1);
    \draw[decoration={markings, mark=at position 0.2 with {\arrow{<}}},
        postaction={decorate},gray
        ] (1,-2) circle (0.5);
    \draw (1,-2.5) node {$\ldots$} (0.5,-1.9) node{$\vdots$} (1.5,-1.9) node{$\vdots$} (0.7,-1.6) node{$*$} (1.3,-1.6) node{$*$};
    \draw[decoration={markings, mark=at position 0.2 with {\arrow{<}}},  postaction={decorate},gray
        ] (2,-1) ellipse (0.5 and 1);
    \draw (2,0) node {$\ldots$} (2,-2) node {$\ldots$} (1.6,-0.3) node{$\vdots$} (2.4,-0.3) node{$\vdots$} (1.53,-0.8) node{$*$} (2.47,-0.8) node{$*$} (1.6,-1.5) node{$\vdots$} (2.4,-1.5) node{$\vdots$} (1.53,-1.2) node{$*$} (2.47,-1.2) node{$*$} (0.25,-3) node{$(1)$};
    \draw[gray,dashed](1.5,-1) -- (2.6,-1);
    \draw[decoration={markings, mark=at position 0.18 with {\arrow{<}}},
        postaction={decorate},gray
        ] (5,0) circle (0.5);
    \draw (5,0.5) node {$*$} (4.5,0.1) node{$\vdots$} (5.5,0.1) node{$\vdots$} (4.7,-0.4) node{$*$} (5.3,-0.4) node{$*$};
    \draw[gray,dashed](5,0.85) -- (5,-0.85);
    \draw[decoration={markings, mark=at position 0.2 with {\arrow{<}}},
        postaction={decorate},gray
        ] (5,-2) circle (0.5);
    \draw (5,-2.5) node {$*$} (4.5,-1.9) node{$\vdots$} (5.5,-1.9) node{$\vdots$} (4.7,-1.6) node{$*$} (5.3,-1.6) node{$*$};
    \draw[gray,dashed](5,-1.15) -- (5,-2.85);
    \draw[decoration={markings, mark=at position 0.2 with {\arrow{<}}},  postaction={decorate},gray
        ] (4,-1) ellipse (0.5 and 1);
    \draw (4,0) node {$*$} (4,-2) node {$*$} (3.6,-0.3) node{$\vdots$} (4.4,-0.3) node{$\vdots$} (3.53,-0.8) node{$*$} (4.47,-0.8) node{$*$} (3.6,-1.5) node{$\vdots$} (4.4,-1.5) node{$\vdots$} (3.53,-1.2) node{$*$} (4.47,-1.2) node{$*$};
    \draw[gray,dashed](4,0.2) -- (4,-2.2);
    \draw[decoration={markings, mark=at position 0.18 with {\arrow{<}}},
        postaction={decorate},gray
        ] (7.5,0) circle (0.5);
    \draw (7.5,0.5) node {$\ldots$} (8,0.1) node{$\vdots$} (7,0.1) node{$\vdots$} (7.2,-0.4) node{$*$} (7.8,-0.4) node{$*$};
    \draw[gray,dashed](7,-1) -- (8,-1);
    \draw[decoration={markings, mark=at position 0.2 with {\arrow{<}}},
        postaction={decorate},gray
        ] (7.5,-2) circle (0.5);
    \draw (7.5,-2.5) node {$\ldots$} (8,-1.9) node{$\vdots$} (7,-1.9) node{$\vdots$} (7.2,-1.6) node{$*$} (7.8,-1.6) node{$*$};
    \draw[decoration={markings, mark=at position 0.2 with {\arrow{<}}},  postaction={decorate},gray
        ] (6.5,-1) ellipse (0.5 and 1);
    \draw (6.5,0) node {$\ldots$} (6.5,-2) node {$\ldots$} (6.1,-0.3) node{$\vdots$} (6.9,-0.3) node{$\vdots$} (6.03,-0.8) node{$*$} (6.97,-0.8) node{$*$} (6.1,-1.5) node{$\vdots$} (6.9,-1.5) node{$\vdots$} (6.03,-1.2) node{$*$} (6.97,-1.2) node{$*$} (5.75,-3) node{$(2)$};
    \draw[gray,dashed](5.8,-1) -- (7.1,-1);
    \draw[decoration={markings, mark=at position 0.18 with {\arrow{<}}},
        postaction={decorate},gray
        ] (-1,-4) circle (0.5);
    \draw (-1,-3.5) node {$*$} (-1.5,-3.9) node{$.$} (-0.5,-3.9) node{$.$} (-1.3,-4.4) node{$*$} (-0.7,-4.4) node{$*$};
    \draw[gray,dashed](-1,-3.25) -- (-1,-4.85);
    \draw[decoration={markings, mark=at position 0.2 with {\arrow{<}}},
        postaction={decorate},gray
        ] (-1,-6) circle (0.5);
    \draw (-1,-6.5) node {$*$} (-1.5,-5.9) node{$\vdots$} (-0.5,-5.9) node{$\vdots$} (-1.3,-5.6) node{$*$} (-0.7,-5.6) node{$*$};
    \draw[gray,dashed](-1,-5.15) -- (-1,-6.85);
    \draw[decoration={markings, mark=at position 0.2 with {\arrow{<}}},  postaction={decorate},gray
        ] (-2,-5) ellipse (0.5 and 1);
    \draw (-2,-4) node {$*$} (-2,-6) node {$*$} (-2.37,-4.3) node{$.$} (-1.63,-4.3) node{$.$} (-2.53,-4.8) node{$*$} (-1.53,-4.8) node{$*$} (-2.4,-5.5) node{$\vdots$} (-1.6,-5.5) node{$\vdots$} (-2.53,-5.2) node{$*$} (-1.53,-5.2) node{$*$} (-1.5,-7) node{$(3)$};
    \draw[gray,dashed](-2,-3.8) -- (-2,-6.2);
    \draw[decoration={markings, mark=at position 0.18 with {\arrow{<}}},
        postaction={decorate},gray
        ] (1,-4) circle (0.5);
    \draw (1,-3.5) node {$*$} (0.5,-3.9) node{$.$} (1.5,-3.9) node{$.$} (0.7,-4.4) node{$*$} (1.3,-4.4) node{$*$};
    \draw[gray,dashed](1,-3.4) -- (1,-4.6);
    \draw[decoration={markings, mark=at position 0.2 with {\arrow{<}}},
        postaction={decorate},gray
        ] (1,-6) circle (0.5);
    \draw (0.5,-5.9) node{$\vdots$} (1.5,-5.9) node{$\vdots$} (0.7,-5.6) node{$*$} (1.3,-5.6) node{$*$} (0.7,-6.4) node{$*$} (1.3,-6.4) node{$*$};
    \draw[gray,dashed](1,-5.4) -- (1,-6.6);
    \draw[decoration={markings, mark=at position 0.2 with {\arrow{<}}},  postaction={decorate},gray
        ] (2,-5) ellipse (0.5 and 1);
    \draw (2,-4) node {$*$} (1.62,-4.3) node{$.$} (2.38,-4.3) node{$.$} (1.53,-4.8) node{$*$} (2.47,-4.8) node{$*$} (1.6,-5.5) node{$\vdots$} (2.4,-5.5) node{$\vdots$} (1.53,-5.2) node{$*$} (2.47,-5.2) node{$*$} (1.72,-5.9) node{$*$} (2.26,-5.9) node{$*$} (1.5,-7) node{$(4)$};
    \draw[gray,dashed](2,-3.9) -- (2,-6.1);
    \draw[decoration={markings, mark=at position 0.18 with {\arrow{<}}},
        postaction={decorate},gray
        ] (4,-4) circle (0.5);
    \draw (4,-3.5) node {$*$} (3.5,-3.9) node{$.$} (4.5,-3.9) node{$.$} (3.7,-4.4) node{$*$} (4.3,-4.4) node{$*$};
    \draw[gray,dashed](4,-3.4) -- (4,-4.6);
    \draw[decoration={markings, mark=at position 0.2 with {\arrow{<}}},
        postaction={decorate},gray
        ] (4,-6) circle (0.5);
    \draw (3.5,-5.9) node{$\vdots$} (4.5,-5.9) node{$\vdots$} (3.7,-5.6) node{$*$} (4.3,-5.6) node{$*$} (4,-6.5) node{$*$};
    \draw[gray,dashed](4,-5.4) -- (4,-6.6);
    \draw[decoration={markings, mark=at position 0.2 with {\arrow{<}}},  postaction={decorate},gray
        ] (5,-5) ellipse (0.5 and 1);
    \draw (5,-4) node {$*$} (4.62,-4.3) node{$.$} (5.38,-4.3) node{$.$} (4.53,-4.8) node{$*$} (5.47,-4.8) node{$*$} (4.6,-5.5) node{$\vdots$} (5.4,-5.5) node{$\vdots$} (4.53,-5.2) node{$*$} (5.47,-5.2) node{$*$} (5,-6) node{$*$} (4.5,-7) node{$(5)$};
    \draw[gray,dashed](5,-3.9) -- (5,-6.1);
    \draw[decoration={markings, mark=at position 0.18 with {\arrow{<}}},
        postaction={decorate},gray
        ] (8,-4) circle (0.5);
    \draw (8,-3.5) node {$*$} (7.5,-3.9) node{$.$} (8.5,-3.9) node{$.$} (7.7,-4.4) node{$*$} (8.3,-4.4) node{$*$};
    \draw[gray,dashed](8,-3.25) -- (8,-4.85);
    \draw[decoration={markings, mark=at position 0.2 with {\arrow{<}}},
        postaction={decorate},gray
        ] (8,-6) circle (0.5);
    \draw (7.5,-5.9) node{$\vdots$} (8.5,-5.9) node{$\vdots$} (7.7,-5.6) node{$*$} (8.3,-5.6) node{$*$} (7.7,-6.4) node{$*$} (8.3,-6.4) node{$*$};
    \draw[gray,dashed](8,-5.15) -- (8,-6.85);
    \draw[decoration={markings, mark=at position 0.2 with {\arrow{<}}},  postaction={decorate},gray
        ] (7,-5) ellipse (0.5 and 1);
    \draw (7,-4) node {$*$} (6.63,-4.3) node{$.$} (7.37,-4.3) node{$.$} (6.47,-4.8) node{$*$} (7.47,-4.8) node{$*$} (6.6,-5.5) node{$\vdots$} (7.4,-5.5) node{$\vdots$} (6.47,-5.2) node{$*$} (7.47,-5.2) node{$*$} (6.72,-5.9) node{$*$} (7.26,-5.9) node{$*$}    (7.5,-7) node{$(6)$};
    \draw[gray,dashed](7,-3.8) -- (7,-6.2);
    \end{tikzpicture}
    \caption{Cut and join of the cycles corresponding to the inner vertices of zigzag cover $\varphi$.}
    \label{fig:cycle-cut-join}
\end{figure}
Definition $\ref{def:mono-zigzag}$ implies that the images of the inner vertices of a component under $\varphi$ form a consecutive sequence of points,
and the unbent tails attached to a piece are either all in-tails or all out-tails.
From Definition $\ref{def:mono-zigzag}$,
inner vertices in in-components may be of any type depicted in Figure $\ref{fig:vertices-in-MZ}$ except for the type $(4)$,
and inner vertices in out-components can be of type $(1)$ and type $(4)$ depicted in Figure $\ref{fig:vertices-in-MZ}$.
Suppose that $(\gamma,\sigma_1,\tau_1,\ldots,\tau_r,\sigma_2)$ is a factorization producing $(\varphi,\rho_s)$ according to Construction $\ref{const2}$.
We denote by $T_i$ (resp. $B_i$) the set of
integers (resp. larger integers) in the transpositions of this factorization which correspond to inner vertices in the component $C_i$.
An in-component intersects with an out-component at a vertex of type $(3)$ or type $(5)$ in Figure $\ref{fig:vertices-in-MZ}$.
We choose the factorization $(\gamma,\sigma_1,\tau_1,\ldots,\tau_r,\sigma_2)$ such that the larger integer in a vertex of type $(3)$ (resp. type $(5)$) is the admissible integer in the out-going (resp. in-coming) edge in in-piece.
Suppose that $C_i$ is an out-component.
The inner vertices in $C_i$ may be of type $(1)$ or type $(4)$ depicted in Figure $\ref{fig:vertices-in-MZ}$.
The vertex of type $(1)$ in $C_i$ must be a symmetric fork vertex, so the two integers in this vertex never appear in other vertices before this vertex.
Moreover, it is possible to choose the factorization such that the larger integer in the vertex of type $(4)$ in $C_i$, say $v$, comes from the even in-coming edge.
In fact, there are three possibilities for a vertex $v$ of type $(4)$.
\begin{itemize}
    \item $v$ is in an in-tail without symmetric fork or symmetric cycle.
    The admissible integer coming from the even in-coming edge never appears in other vertices before $v$,
    so it is possible to require this integer to be the larger one in $v$.
    \item $v$ is in a weight $2$ in-tail with symmetric fork or symmetric cycle.
    It follows from Lemma $\ref{lem:invol-sign-change}$ and Lemma $\ref{lem:numb-trans}$ that the admissible integer in the even in-coming edge can be chosen as the larger integer in the former vertex of $v$.
    Note that the larger integer in the former vertex does not appear before.
    \item $v$ is in a weight $>2$ in-tail with symmetric fork or symmetric cycle.
    If the signs of $v$ and its former vertex $v'$ are the same,
    the admissible integer in the even in-coming edge can be chosen as the larger integer in $v'$.
    Otherwise, it follows from Lemma $\ref{lem:invol-sign-change}$ that the two integers in $v$ have never appeared before.
\end{itemize}
Therefore, we obtain that $B_i\cap T_{i+1}=B_i\cap T_{i-1}=\emptyset$ and the sets $B_i$ and $T_i$ satisfy condition $(*)$ in Lemma $\ref{lem:mono-fac1}$.
Suppose that $C_j$ is an in-component.
The inner vertices in $C_j$ may be of any type depicted in Figure $\ref{fig:vertices-in-MZ}$ except for the type $(4)$.
Moreover, vertices of type $(1)$ or type $(2)$ can only be inner vertices in weight $2$ bent tails with symmetric fork or symmetric cycle.
It follows from Lemma $\ref{lem:invol-sign-change}$ and Lemma $\ref{lem:numb-trans}$ that the integers in a vertex $\bar v$ in $C_j$ of type $(1)$ or type $(2)$ are the same as the integers in the first vertex of the bent tail containing $\bar v$.
Note that we have already required that the larger integer in a vertex of type $(3)$ (resp. type $(5)$) is the admissible integer in the out-going (resp. in-coming) edge in in-piece.
The vertex of type $(6)$ in $C_j$ is in an unbent out-tail without symmetric fork or symmetric cycle.
We require that the larger integer in a vertex of type $(6)$ is the admissible integer in the odd out-going edge of this vertex.
In the in-component $C_j$,
To see the requirements about integers in the inner vertices in the in-component $C_j$ are compatible, we note that:
\begin{itemize}
    \item if there is a vertex of type $(1)$ in $C_j$, the odd edges in $C_j$ are all of weight $1$.
    \item if there is no vertex of type $(1)$ in $C_j$, the integers in vertex of type $(3)$ do not appear before.
    Suppose that the former vertex of the vertex $\hat v$ of type $(5)$ in $C_j$ is of type $(6)$.
    If they have the same signs,
    the larger integer of $\hat v$ is the same as the larger integer of the former vertex of $\hat v$.
    If they have different signs, Lemma $\ref{lem:invol-sign-change}$ and Lemma $\ref{lem:numb-trans}$ implies that the larger integer of $\hat v$ does not appear before.
\end{itemize}
It is obvious that
the subset of $B_j$ consisting of larger integers in vertices of type $(6)$ is disjoint with $T_{j-1}\cup T_{j+1}$.
Condition $(2)$ of Definition $\ref{def:mono-zigzag}$ and the above requirements about transpositions imply that the larger integer in $B_j$ corresponding to vertex of type $(5)$ is not in $T_{j+1}$.
Hence, $B_j\cap T_{j+1}=B_j\cap T_{j-1}=\emptyset$,
and the sets $B_j$, $T_j$ satisfy the condition $(*)$ in Lemma $\ref{lem:mono-fac1}$.
Since $B_k\cap T_{k+1}=B_k\cap T_{k-1}=\emptyset$
and $B_k$, $T_k$ satisfy the condition $(*)$ in Lemma $\ref{lem:mono-fac1}$ for any component $C_k$,
we have a factorization satisfying condition $(*)$.
Hence we complete the proof.
\end{proof}

We denote by $\vec M_{g}(\lambda,\mu)$ the set of monotone zigzag covers.
Let $\vec{N}(\varphi,s)$ denote the number of real monotone factorizations of type $(g,\lambda,\mu;s)$ associated to a real monotone zigzag cover $(\varphi,\rho_s)$, where $\rho_s$ is the unique colouring of $\varphi$ determined by a simple splitting $\undl x=\undl x^+\sqcup\undl x^-$ with $|\undl x^+|=s$.
We call the following number
\begin{equation}\label{eq:mon-zig-numb}
    \vec Z_{g}(\lambda,\mu):=\sum_{\varphi\in\vec M_{g}(\lambda,\mu)}\vec N(\varphi),
\end{equation}
the \textit{monotone zigzag number},
where $\vec{N}(\varphi):=\min_{0\leq s\leq r}\vec{N}(\varphi,s)$.

\begin{remark}
The monotone zigzag number $\vec{Z}_g(\lambda,\mu)$ is a lower bound of the number
of real monotone factorizations of type $(g,\lambda,\mu;s)$ for any $s$,
so we obtain
$$
\vec{Z}_g(\lambda,\mu)\leq\vec{H}_g^\rb(\lambda,\mu).
$$
In Section $\ref{sec:asym-real-mono}$,
we achieve the asymptotic growth of $\vec{H}_g^\rb(\lambda,\mu)$ as the degree goes to infinity by constructing a particular family of real tropical covers.
The particular family of real tropical covers are obtained by gluing a monotone zigzag cover with some monotone non-zigzag covers.
In the procedure of the construction,
the non-zero result in Lemma $\ref{lem:mono-zigzag1}$ plays a crucial role.
\end{remark}

\subsection{Universally monotone zigzag numbers}
\label{subsec:uni-mzn}

In this subsection, we give a lower bound of the real monotone double Hurwitz number  $\vec\hl^\rb_g(\lambda,\mu)$ relative to arbitrary splittings.
\begin{definition}
\label{def:mon-zigzag2}
A monotone zigzag cover $\varphi:C\to T\pb^1$ is called \textit{universally monotone} if for any in-piece of the string in $C$ there is at most one unbent out-tail attached to it.
\end{definition}

\begin{lemma}\label{lem:mono-zigzag2}
Let $\varphi:C\to T\pb^1$ be a universally monotone zigzag cover of type $(g,\lambda,\mu,\undl x)$.
Then for any splitting $\undl x=\undl x^+\sqcup\undl x^-$ of $\undl x$ into positive and negative points,
the number of real monotone factorizations
producing $(\varphi,\rho)$ according to Construction $\ref{const2}$ is non-zero,
where $\rho$ is the unique colouring of $\varphi$
determined by the splitting $\undl x=\undl x^+\sqcup\undl x^-$.
\end{lemma}

\begin{proof}
The proof of this Lemma is almost the same as that of Lemma $\ref{lem:mono-zigzag1}$,
so we only point out the difference between these two cases and omit the repeated part.
Let $\varphi$ be any universally monotone zigzag cover, and let $\undl x=\undl x^+\sqcup\undl x^-$ be any splitting of $\undl x$ into positive and negative points.
We need prove that there is a real factorization, which satisfies condition $(*)$ in Lemma $\ref{lem:mono-fac1}$, producing $(\varphi,\rho)$ according to Construction $\ref{const2}$.
By the rule for choosing transpositions in the proof of Lemma $\ref{lem:mono-zigzag1}$,
we obtain real factorizations which produce $(\varphi,\rho)$ according to Construction $\ref{const2}$.
For a factorization obtained in this way,
we use $T_i$ (resp. $B_i$) to denote the set of integers (resp. larger integers) in the transpositions corresponding to the inner vertices in the component $C_i$.
The rule for choosing integers in vertices of type $(3)$ and type $(5)$ depicted in Figure $\ref{fig:vertices-in-MZ}$ in the proof of Lemma $\ref{lem:mono-zigzag1}$ and Definition $\ref{def:mono-zigzag}$ imply that $B_i\cap T_{i+1}=B_i\cap T_{i-1}=\emptyset$.
The sets $B_i$ and $T_i$ of an out-component $C_i$ still satisfy the condition $(*)$ in Lemma $\ref{lem:mono-fac1}$.
However, if $C_i$ is an in-component which is attached with more than one unbent out-tails, the sets $B_i$ and $T_i$ may not satisfy the condition $(*)$ in Lemma $\ref{lem:mono-fac1}$ (see Figure $\ref{fig:mon-zigzag-condition}$).
\begin{figure}[ht]
    \centering
    \begin{tikzpicture}
    \draw[line width=0.3mm] (-4,0)--(-3,0);
    \draw[line width=0.3mm,gray]  (-3,0)--(0,-2.4);
    \draw[line width=0.3mm]  (0,-2.4)--(2,-2.4);
    \draw[line width=0.3mm]  (-2,-0.8)--(2,-0.8);
    \draw[line width=0.3mm]  (-1,-1.6)--(2,-1.6);
    \draw (-2.3,-0.4) node{\tiny{$o_1$}} (-2,-1.2) node{\tiny{$o_1-e_1$}} (-0.7,-2) node{\tiny{$1$}} (-0.5,-1) node{\tiny{$e_1$}} (0.5,-1.8) node{\tiny{$o_1-e_1-1$}};
    \draw (-3,0.25) node{\tiny{$+$}} (-1.8,-0.7) node{\tiny{$-$}} (-1,-1.8) node{\tiny{$+$}};
    \end{tikzpicture}
    \caption{An in-component attached with two unbent out-tails.
    The signs of the first three inner vertices are labelled by $\pm$, and the weights of the edges are also given.}
    \label{fig:mon-zigzag-condition}
\end{figure}
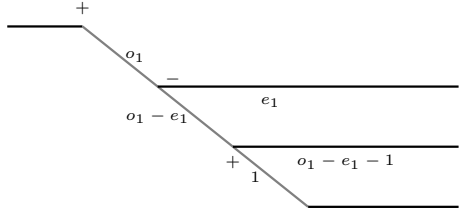
In a universally monotone zigzag cover,
any in-component is attached with at most one unbent out-tail, so the sets $B_i$ and $T_i$ corresponding to an in-component of a universally monotone zigzag cover always satisfy the condition $(*)$ in Lemma $\ref{lem:mono-fac1}$.
\end{proof}

We denote by $\vec{\ml}_{g}(\lambda,\mu)$ the set of universally monotone zigzag covers of type $(g,\lambda,\mu,\undl x)$.
Let $\vec{N}(\varphi,\undl S(s))$ denote the number of real monotone factorizations associated to a universally real monotone zigzag cover $(\varphi,\rho)$,
where $\undl S(s)$ is a sequence of signs with $s$ positive entries, and $\rho$ is the unique colouring of $\varphi$ which is determined by the splitting $\undl x=\undl x^+\sqcup\undl x^-$ associated to the sequence of signs $\undl S(s)$.
Put
$$
\vec\zl_{g}(\lambda,\mu):=\sum_{\varphi\in\vec{\ml}_{g}(\lambda,\mu)}\vec N(\varphi),
$$
where $\vec N(\varphi)=\min_{s,\undl S(s)}\vec{N}(\varphi,\undl S(s))$.
The number $\vec\zl_{g}(\lambda,\mu)$ is called the \textit{universally monotone zigzag number}.
\begin{proposition}\label{prop:lower-bound}
Fix $g\geq0$, $d\geq1$, and two partitions $\lambda$, $\mu$ with $|\lambda|=|\mu|=d$.
Then the real monotone double Hurwitz
number $\vec\hl^{\rb}_g(\lambda,\mu)$ relative to arbitrary splitting is bounded from below by the universally monotone zigzag number:
$$
\vec\zl_{g}(\lambda,\mu)\leq
\vec\hl^{\rb}_g(\lambda,\mu).
$$
\end{proposition}

\begin{proof}
It is straightforward from the definition, so we omit it.
\end{proof}

\begin{proposition}
\label{prop:nonvanish-mon-zigzag}
Let $m\geq1$ be an integer. There is a universally monotone zigzag cover $\varphi:C\to T\pb^1$ of type $(g,1^{2m+1},1^{2m+1},\undl x)$.
Moreover, any edge of $C$ is of weight $1$ or $2$, and the number $\vec N(\varphi)$ is bounded from below by $m!$.
\end{proposition}

\begin{proof}
We pick a string $S$ with two ends, then attach tails of the first type depicted in Figure $\ref{fig:zigzag}$, corresponding to $(1^{2m},1^{2m})$,
to the string $S$.
We place $g$ symmetric cycles on one of such tails.
Then we have a tropical curve $C$ (See Figure $\ref{fig:mon-zigzag1}$ for an example).
\begin{figure}[ht]
    \centering
    \begin{tikzpicture}
    \draw[line width=0.4mm,gray] (-4,-0.5)--(2,-1)--(-1,-1.5)--(1,-2)--(-2,-2.5)--(4,-3);
    \draw[line width=0.4mm] (2,-1)--(3,-1);
    \draw[line width=0.4mm] (4,-0.7)--(3,-1)--(4,-1.3);
    \draw[line width=0.4mm] (-1,-1.5)--(-1.7,-1.5);
    \draw[line width=0.4mm] (-4,-1.8)--(-1.7,-1.5)--(-4,-1.2);
    \draw[line width=0.4mm] (1,-2)--(1.7,-2);
    \draw[line width=0.4mm] (4,-1.7)--(1.7,-2)--(4,-2.3);
    \draw[line width=0.4mm] (-2,-2.5)--(-3,-2.5);
    \draw[line width=0.4mm] (-4,-2.8)--(-3,-2.5)--(-4,-2.2);
    \draw[line width=0.4mm,gray] (-3,-0.4) node{\tiny{$S$}};
    \draw[line width=0.3mm] (-4,-3.5)--(4,-3.5);
    \foreach \Point in {(-3,-3.5), (-2,-3.5),(-1.7,-3.5),(-1,-3.5),(1,-3.5),(1.7,-3.5),(2,-3.5),(3,-3.5)}
    \draw[fill=black] \Point circle (0.05);
    \draw[line width=0.4mm,gray] (-3,-3.8) node{\tiny{$x_1$}} (-2,-3.8) node{\tiny{$x_2$}} (-1.7,-3.8) node{\tiny{$x_3$}} (-1,-3.8) node{\tiny{$x_4$}} (1,-3.8) node{\tiny{$x_5$}} (1.7,-3.8) node{\tiny{$x_6$}} (2,-3.8) node{\tiny{$x_7$}} (3,-3.8) node{\tiny{$x_8$}};
    \end{tikzpicture}
    \caption{A universally monotone zigzag cover of type $(0,1^5,1^5,\undl x)$.}
    \label{fig:mon-zigzag1}
\end{figure}
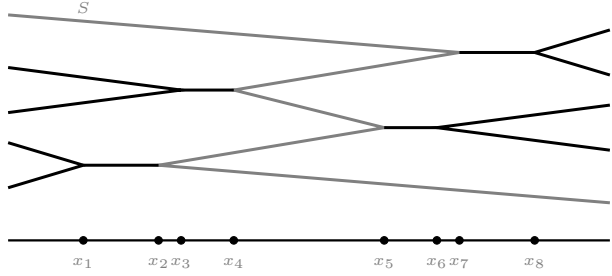
In particular, all the even edges of $C$ are of weight $2$, and odd edges of $C$ are of weight $1$.
It is obvious that there is a total order of the inner vertices of $C$ such that the tropical cover $\varphi:C\to T\pb^1$ is a universally monotone zigzag cover.
Moreover, there are at least $m!$ monotone factorizations associated to the tropical cover $\varphi:C\to T\pb^1$.
In fact, we can choose $\{m+1,m+2,\ldots,2m+1\}$ to be the set of larger integers for the transpositions in the in-components,
then there are at least $m!$ ways to arrange the smaller integers for the transpositions corresponding to the symmetric fork vertices in the bent in-tails of the in-components.
Suppose that $(\sigma_1,\tau_1,\ldots,\tau_r,\sigma_2)$
is such a monotone factorization.
Since the weight of any edge of the tropical curve $C$ is $1$ or $2$,
the transpositions corresponding to
the vertices of $C$ either join two cycles of length $1$ or
cut a cycle of length $2$ into two cycles of length $1$.
Let $\undl x^+\sqcup\undl x^-$ be the splitting of $\undl x$ into positive and negative branch points corresponding to a sequence of signs $\undl S(s)$, where $0\leq s\leq r$.
Suppose that $\rho$ is the unique colouring of $\varphi$ determined by the splitting $\undl x=\undl x^+\sqcup\undl x^-$.
From the analysis on the effect of the
conjugation with an involution in Lemma $\ref{lem:invol-sign-change}$,
it is easy to see that there is an involution
$\gamma_\rho$ such that the tuple
$(\gamma_\rho,\sigma_1,\tau_1,\ldots,\tau_r,\sigma_2)$
is a real monotone factorization producing the real universally monotone zigzag cover $(\varphi,\rho)$ according to Construction $\ref{const2}$.
Therefore, we obtain $\vec N(\varphi)\geq m!$.
\end{proof}

The following proposition shows that only zigzag covers $\varphi:C\to T\pb^1$ with edges of weight $1$ or $2$ have nonzero number $\vec N(\varphi)$.
\begin{proposition}
\label{prop:optimal-asymp}
Let $\varphi:C\to T\pb^1$ be a zigzag cover of type $(0,1^{2m+1},1^{2m+1},\undl x)$, where $m\geq1$.
Suppose that there is a weight $\omega>2$ inner edge in the string $S\subset C$.
Then the number $\vec N(\varphi)=0$.
\end{proposition}

\begin{proof}
We prove this proposition by ruling out the following $2$ possibilities.

{\bf Case $(1)$}: The tropical curve $C$ has no bent vertex.

Suppose that the maximal weight of inner edges in $S$ is $3$, and $E$ is an inner edge of weight $3$.
Note that an unbent in-tail of weight $2$ and an unbent out-tail of weight $2$ are attached to $E$ at two vertices of it.
Suppose that the vertex of $E$ which is the intersection vertex of the unbent in-tail and $S$ is mapped to $x_i\in\undl x$ by $\varphi$.
Let $\rho$ be the unique colouring of $\varphi$ determined by the simple sequence of signs $\undl S(i)$.
Then there is no real monotone factorization producing the real tropical cover $(\varphi,\rho)$.
In fact, if the integer $a$ associated to the weight $1$ incoming edge of the vertex corresponding to $x_i$ is smaller than the larger integer $b$ of the fork vertex of the in-tail attached to $E$, the larger integer of the fork vertex of the out-tail attached to $E$ is smaller than $b$.
Otherwise,
the larger integer of the intersection vertex of $E$ and the out-tail is smaller than $a$.

Suppose that the maximal weight of inner edges in $S$ is $\omega_1>3$.
Then there is no real monotone factorization producing the real tropical cover $(\varphi,\rho_1)$, where $\rho_1$ is the unique colouring of $\varphi$ determined by the simple sequence of signs $\undl S(r)$.
In fact, in this case there is a pair of unbent in-tail $t_l$ and unbent out-tail $t_r$ satisfying the following conditions:
\begin{itemize}
    \item the weight $\omega_2$ of the incoming odd edge of the intersection vertex $v$ of $t_l$ and $S$ equals the weight of the outgoing odd edge of the intersection vertex $v'$ of $t_r$ and $S$;
    \item the weights of odd edges in $S$ between $v$ and $v'$ are larger than $\omega_2$.
    \item the points in $\undl x$ corresponding to $v$ and $v'$ are not adjacent;
\end{itemize}
Note that the two transpositions corresponding to $v$ and $v'$ are the same, but the larger integer of $v$ and $v'$ is different from the larger integers of the vertices in $S$ between $v$ and $v'$.

{\bf Case $(2)$}: The tropical curve $C$ has bent vertices.

Note that the two ends of any piece of $S$ are of weight $1$, so the weight $\omega$ inner edge in $S$ is an inner edge of a piece of $S$.
Let $S_i$ be a piece containing a weight $\omega$ inner edge.
We cut the graph $C$ at the two bent vertices of $S_i$, then the subgraph of $C$ containing $S_i$ is a tropical curve with no bent vertex.
The argument in case $(1)$ implies that there is a colouring $\rho_2$ of $\varphi$ such that there is no real monotone factorization producing $(\varphi,\rho_2)$.
\end{proof}

\section{Asymptotics for real
monotone double Hurwitz numbers}
\label{sec:asym-real-mono}
In this section, we study the asymptotic growth of the real monotone double Hurwitz numbers.

\subsection{Case $\textrm{I}$: asymptotics for simple splittings}
We prove the logarithmic equivalence of the real monotone double Hurwitz numbers relative to simple splittings and the monotone double Hurwitz numbers.

\begin{proposition}
\label{prop:asym-simple}
Let $m\geq1$ be an integer.
For any $s\in\{0,1,\ldots,4m-2\}$,
there are at least $m!$ real tropical covers of type
$(0,(2,1^{2m-1}),(2,1^{2m-1}),\undl x)$
whose real branch points possess a simple splitting
$\undl x=\undl x^+\sqcup\undl x^-$ with $|\undl x^+|=s$.
Moreover, for any such real tropical cover $\varphi:C\to T\pb^1$
the number $\vec{N}(\varphi,s)\geq m!$.
\end{proposition}

\begin{proof}
For any $s\in\{0,1,\ldots,4m-2\}$, we first show that there
is a real tropical cover of type
$(0,(2,1^{2m-1}),(2,1^{2m-1}),\undl x)$ such that the real branch points possess a simple splitting
$\undl x=\undl x^+\sqcup\undl x^-$ with $|\undl x^+|=s$,
then we show that one obtain $m!$ such real tropical covers from the given cover.

We start from $4$ types of monotone components which are depicted in Figure $\ref{fig:mon-comps}$.
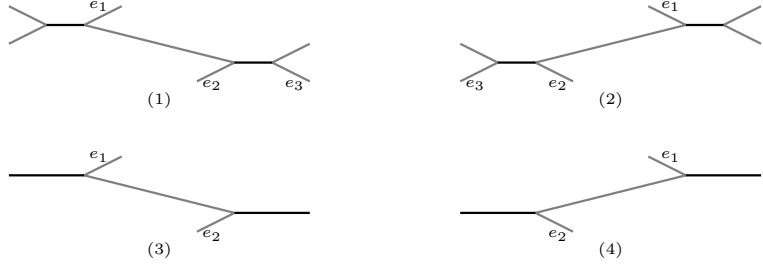
\begin{figure}[ht]
    \centering
    \begin{tikzpicture}
    \draw[line width=0.3mm,gray] (-3.5,2.25)--(-4,2)--(-2,1.5)--(-2.5,1.25);
    \draw[line width=0.3mm,gray]  (-5,2.25)--(-4.5,2)--(-5,1.75);
    \draw[line width=0.3mm]  (-4.5,2)--(-4,2);
    \draw[line width=0.3mm,gray]  (-1,1.25)--(-1.5,1.5)--(-1,1.75);
    \draw[line width=0.3mm]  (-1.5,1.5)--(-2,1.5);
    \draw[line width=0.3mm,gray] (3.5,2.25)--(4,2)--(2,1.5)--(2.5,1.25);
    \draw[line width=0.3mm,gray]  (5,2.25)--(4.5,2)--(5,1.75);
    \draw[line width=0.3mm]  (4.5,2)--(4,2);
    \draw[line width=0.3mm,gray]  (1,1.25)--(1.5,1.5)--(1,1.75);
    \draw[line width=0.3mm]  (1.5,1.5)--(2,1.5);
    \draw[line width=0.3mm,gray] (-3.5,0.25)--(-4,0)--(-2,-0.5)--(-2.5,-0.75);
    \draw[line width=0.3mm]  (-5,0)--(-4,0);
    \draw[line width=0.3mm]  (-1,-0.5)--(-2,-0.5);
    \draw[line width=0.3mm,gray] (3.5,0.25)--(4,0)--(2,-0.5)--(2.5,-0.75);
    \draw[line width=0.3mm]  (5,0)--(4,0);
    \draw[line width=0.3mm]  (1,-0.5)--(2,-0.5);
    \draw (-3,1) node{\tiny{$(1)$}} (3,1) node{\tiny{$(2)$}} (-3,-1) node{\tiny{$(3)$}} (3,-1) node{\tiny{$(4)$}};
    \draw (-3.8,2.25) node{\tiny{$e_1$}} (-2.3,1.2) node{\tiny{$e_2$}} (-1.2,1.2) node{\tiny{$e_3$}};
    \draw (3.8,2.25) node{\tiny{$e_1$}} (2.3,1.2) node{\tiny{$e_2$}} (1.2,1.2) node{\tiny{$e_3$}};
    \draw (-3.8,0.25) node{\tiny{$e_1$}} (-2.3,-0.8) node{\tiny{$e_2$}};
    \draw (3.8,0.25) node{\tiny{$e_1$}} (2.3,-0.8) node{\tiny{$e_2$}};
    \end{tikzpicture}
    \caption{Monotone components. Edges of weight $1$ are drawn in gray, and edges of weight $2$ are drawn in black.}
    \label{fig:mon-comps}
\end{figure}
We take $m-1$ monotone components $C_1,\ldots,C_{m-1}$ of type $(1)$ or $(2)$ in Figure $\ref{fig:mon-comps}$ and $1$ monotone component $C_m$ of type $(3)$ or $(4)$.
A monotone component $C_i$ is glued with component $C_{i+1}$ according to the following rule:
\begin{itemize}
    \item edge $e_2$ (resp. $e_3$) of a component of type $(1)$ is glued with edge $e_1$ of a component of type $(1)$ or $(3)$ (resp. type $(2)$ or $(4)$).
    \item edge $e_2$ (resp. $e_3$) of a component of type $(2)$ is glued with edge $e_1$ of a component of type $(2)$ or $(4)$ (resp. type $(1)$ or $(3)$).
\end{itemize}
Let $C$ be the graph obtained by gluing the monotone components $C_1,\ldots,C_m$ as above.
By shrinking or extending the lengths of the inner edges of $C$ properly, we obtain a total order of the inner vertices such that the $4$ (resp. $2$) inner vertices of a component of type $(1)$ or $(2)$ (resp. type $(3)$ or $(4)$) are $4$ (resp. $2$) consecutive vertices in this total order.
We obtain a tropical cover $\varphi:C\to T\pb^1$.

When $s$ is even, the tropical cover $\varphi:C\to T\pb^1$ has a colouring $\rho$ such that the real tropical cover $(\varphi,\rho)$ possesses a simple splitting $\undl x=\undl x^+\sqcup\undl x^-$ with $|\undl x^+|=s$.
When $s$ is odd, the above statement about the colouring may be wrong because the gluing procedure produces even inner edges without symmetric forks or symmetric cycles.
See Figure $\ref{fig:mon-tropical}$ for an example.
\begin{figure}[ht]
    \centering
    \begin{tikzpicture}
    \draw[line width=0.3mm,gray] (-3.5,2.25)--(-4,2)--(-2,1.5)--(-2.5,1.25);
    \draw[line width=0.3mm,gray]  (-5,2.25)--(-4.5,2)--(-5,1.75);
    \draw[line width=0.3mm]  (-4.5,2)--(-4,2);
    \draw[line width=0.3mm,gray]  (-1,1.75)--(-1.5,1.5)--(3,1)--(1,0.5)--(1.5,0.25);
    \draw[line width=0.3mm]  (-1.5,1.5)--(-2,1.5);
    \draw[line width=0.3mm]  (3,1)--(3.5,1);
    \draw[line width=0.3mm,gray]  (4,0.75)--(3.5,1)--(4,1.25);
    \draw[line width=0.3mm] (0.5,0.5)--(1,0.5);
    \draw[line width=0.3mm,gray]  (0,0.75)--(0.5,0.5)--(0,0.25);
    \draw[line width=0.3mm,gray] (2.5,-0.75)--(3,-1)--(1,-1.5)--(1.5,-1.75);
    \draw[line width=0.3mm,gray]  (4,-0.75)--(3.5,-1)--(4,-1.25);
    \draw[line width=0.3mm]  (3.5,-1)--(3,-1);
    \draw[line width=0.3mm,gray]  (0,-1.25)--(0.5,-1.5)--(-4,-2.25)--(-2,-2.75)--(-2.5,-3);
    \draw[line width=0.3mm]  (0.5,-1.5)--(1,-1.5);
    \draw[line width=0.3mm]  (-4.5,-2.25)--(-4,-2.25);
    \draw[line width=0.3mm,gray]  (-5,-2.5)--(-4.5,-2.25)--(-5,-2);
    \draw[line width=0.3mm]  (-2,-2.75)--(-1.5,-2.75);
    \draw[line width=0.3mm,gray]  (-1,-2.5)--(-1.5,-2.75)--(-1,-3);
    \draw (-0.5,-0.25) node{\tiny{$(1)$}} (-0.5,-3.5) node{\tiny{$(2)$}};
    \draw (-4.5,1.8) node{\tiny{$+$}} (-4,1.8) node{\tiny{$+$}} (-2,1.3) node{\tiny{$+$}};
    \draw (-1.5,1.3) node{\tiny{$-$}} (0.5,0.3) node{\tiny{$-$}} (1,0.3) node{\tiny{$-$}} (3,0.8) node{\tiny{$-$}} (3.5,0.8) node{\tiny{$-$}};
    \draw (-4.5,-2.5) node{\tiny{$+$}} (-4,-2.5) node{\tiny{$+$}} (-2,-3) node{\tiny{$+$}};
    \draw (-1.5,-3) node{\tiny{$+$}} (0.5,-1.7) node{\tiny{$+$}} (1,-1.7) node{\tiny{$-$}} (3,-1.2) node{\tiny{$-$}} (3.5,-1.2) node{\tiny{$-$}};
    \end{tikzpicture}
    \caption{Tropical covers with no colouring determined by the present simple splittings.}
    \label{fig:mon-tropical}
\end{figure}
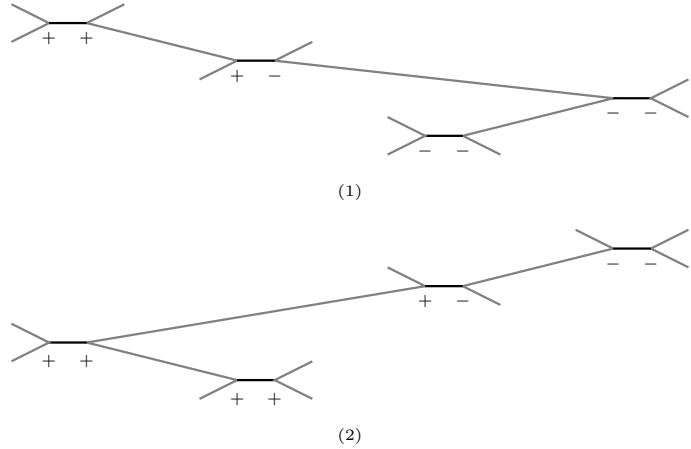
More precisely, when the inner edge connecting the $s$-th and $(s+1)$-th inner vertices is an even edge and it is not attached with a symmetric fork, the tropical cover $\varphi:C\to T\pb^1$ has no colouring such that the associated splitting $\undl x=\undl x^+\sqcup\undl x^-$ is simple and $|\undl x^+|=s$ with $s$ odd.
In this case, the $s$-th inner vertex $v_s$ is not a vertex of the component of type $(3)$ or $(4)$.
We construct a new tropical cover $\varphi':C'\to T\pb^1$ from $\varphi$ such that $\varphi'$ has an expected colouring.
Now we assume that the inner edge connecting the $s$-th and $(s+1)$-th inner vertices is an even edge and it is not attached with a symmetric fork, where $s$ is odd.
In the tropical curve $C$, either there is a bent in-tail $t_l$ with symmetric fork such that the bent vertex of $t_l$ and the bent vertex of $C_m$ with peak pointing to the left are in different sides of the line $x=x_s$ or there is a bent out-tail $t_r$ with symmetric fork such that the bent vertex of $t_r$ and the bent vertex of $C_m$ with peak pointing to the right are in different sides of the line $x=x_s$.
The tropical curve $C'$ is obtained by either exchanging $t_l$ with the in-tail of $C_m$ or exchanging $t_r$ with the out-tail of $C_m$.
See Figure $\ref{fig:mon-tropical-modified}$ for an example.
\begin{figure}[ht]
    \centering
    \begin{tikzpicture}
    \draw[line width=0.3mm,gray] (6,2.25)--(-4,2)--(-2,1.5)--(-5,1.25);
    \draw[line width=0.3mm,gray]  (-5,2.25)--(-4.5,2)--(-5,1.75);
    \draw[line width=0.3mm]  (-4.5,2)--(-4,2);
    \draw[line width=0.3mm,gray]  (6,1.75)--(-1.5,1.5)--(5,1)--(3,0.5)--(6,0.25);
    \draw[line width=0.3mm]  (-1.5,1.5)--(-2,1.5);
    \draw[line width=0.3mm]  (5,1)--(5.5,1);
    \draw[line width=0.3mm,gray]  (6,0.75)--(5.5,1)--(6,1.25);
    \draw[line width=0.3mm] (2.5,0.5)--(3,0.5);
    \draw[line width=0.3mm,gray]  (-5,0.75)--(2.5,0.5)--(0,0)--(2,-0.25)--(-5,-0.75);
    \draw[line width=0.3mm] (2,-0.25)--(6,-0.25);
    \draw[line width=0.3mm] (-5,0)--(0,0);
    \draw[line width=0.3mm] (-5,-1)--(6,-1);
    \foreach \Point in {(-4.5,-1), (-4,-1),(-2,-1),(-1.5,-1),(0,-1),(2,-1),(2.5,-1),(3,-1),(5,-1),(5.5,-1)}
    \draw[fill=black] \Point circle (0.05);
    \draw[line width=0.4mm] (-4.5,-1.25) node{\tiny{$+$}} (-4,-1.25) node{\tiny{$+$}} (-2,-1.25) node{\tiny{$+$}} (-1.5,-1.25) node{\tiny{$-$}} (0,-1.25) node{\tiny{$-$}} (2,-1.25) node{\tiny{$-$}} (2.5,-1.25) node{\tiny{$-$}} (3,-1.25) node{\tiny{$-$}} (5,-1.25) node{\tiny{$-$}} (5.5,-1.25) node{\tiny{$-$}};
    \draw (0.5,-1.5) node{\tiny{$(1)$}};
    \draw[line width=0.3mm,gray] (6,-2)--(-4,-2.25)--(-2,-2.75)--(-5,-3);
    \draw[line width=0.3mm]  (-5,-2.25)--(-4,-2.25);
    \draw[line width=0.3mm,gray]  (6,-2.5)--(-1.5,-2.75)--(5,-3.25)--(3,-3.75)--(6,-4);
    \draw[line width=0.3mm]  (-1.5,-2.75)--(-2,-2.75);
    \draw[line width=0.3mm]  (5,-3.25)--(5.5,-3.25);
    \draw[line width=0.3mm,gray]  (6,-3.5)--(5.5,-3.25)--(6,-3);
    \draw[line width=0.3mm] (2.5,-3.75)--(3,-3.75);
    \draw[line width=0.3mm,gray]  (-5,-3.5)--(2.5,-3.75)--(0,-4.25)--(2,-4.5)--(-5,-5);
    \draw[line width=0.3mm] (2,-4.5)--(6,-4.5);
    \draw[line width=0.3mm] (-0.5,-4.25)--(0,-4.25);
    \draw[line width=0.3mm,gray]  (-5,-4)--(-0.5,-4.25)--(-5,-4.5);
    \draw[line width=0.3mm] (-5,-5.25)--(6,-5.25);
    \foreach \Point in {(-4,-5.25),(-2,-5.25),(-1.5,-5.25),(-0.5,-5.25),(0,-5.25),(2,-5.25),(2.5,-5.25),(3,-5.25),(5,-5.25),(5.5,-5.25)}
    \draw[fill=black] \Point circle (0.05);
    \draw[line width=0.4mm] (-4,-5.5) node{\tiny{$+$}} (-2,-5.5) node{\tiny{$+$}} (-1.5,-5.5) node{\tiny{$+$}} (-0.5,-5.5) node{\tiny{$-$}} (0,-5.5) node{\tiny{$-$}} (2,-5.5) node{\tiny{$-$}} (2.5,-5.5) node{\tiny{$-$}} (3,-5.5) node{\tiny{$-$}} (5,-5.5) node{\tiny{$-$}} (5.5,-5.5) node{\tiny{$-$}};
    \draw (0.5,-5.75) node{\tiny{$(2)$}};
    \end{tikzpicture}
    \caption{A tropical cover depicted in $(1)$ is modified to a tropical cover depicted in $(2)$. Gray edges are of weight $1$ and black edges are of weight $2$.}
    \label{fig:mon-tropical-modified}
\end{figure}
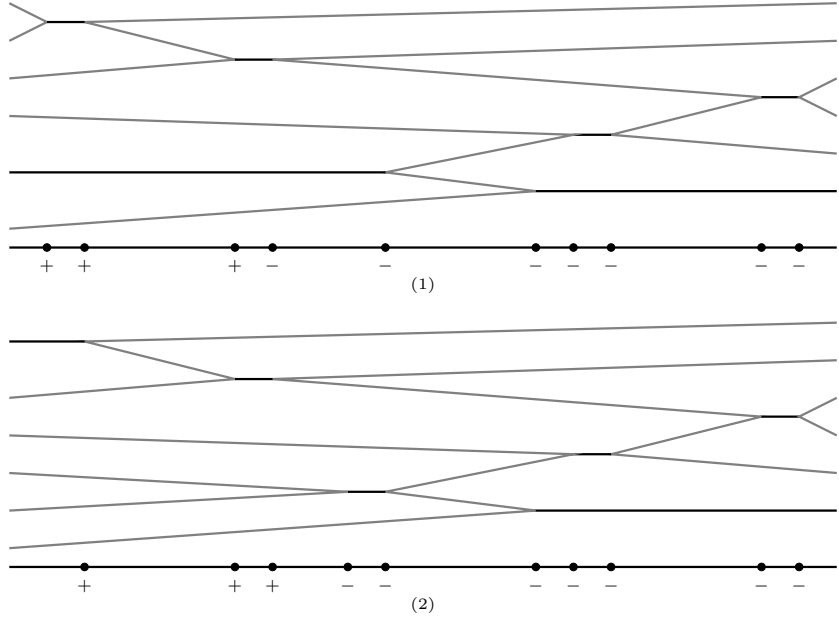
In the graph $C'$, there are two modified components having $3$ inner vertices.
We choose a total order of the inner vertices of $C'$ such that the inner vertices of any component are consecutive in the total order.
In the tropical cover $\varphi':C'\to T\pb^1$,
all the edges intersecting the line $x=x'$,
where $x_s<x'<x_{s+1}$, are odd edges,
so there is a colouring $\rho'$ on $\varphi'$ such that the real tropical cover $(\varphi',\rho')$ possesses a simple splitting $\undl x=\undl x^+\sqcup\undl x^-$ with $|\undl x^+|=s$.

Once the components $C_1,\ldots,C_m$ are fixed, the tropical covers $\varphi$ and $\varphi'$ have the same combinatorial structure, \textit{i.e.} the consecutive sequences of inner vertices of different components are arranged according to same permutation of $\{1,2,\ldots,m\}$.
Conversely, for any permutation of $\{1,2,\ldots,m\}$, there is a choice for the components $C_1,\ldots,C_m$ such that when the inner vertices of any component are arranged consecutively, the sequences of inner vertices of different components are arranged according to that permutation.
Since there are $m!$ permutations of order $m$, we have $m!$ real tropical covers for any $s$ satisfying the requirements.

Now we show that for any real tropical cover $(\varphi,\rho_s)$
of type $(0,(2,1^{2m-1}),(2,1^{2m-1}),\undl x)$ obtained above,
the number $\vec{N}(\varphi,s)\geq m!$.
For any component $C_i$ in the tropical curve $C$, there is at least one in-end of the in-tail.
We choose $\{m+2,m+3,\ldots,2m+1\}$ to be the set of larger integers of the transpositions corresponding to the inner vertices of the $m$ components.
The smaller integer of the transposition corresponding to the inner vertex of the in-tail of any component can be chosen to be any integer smaller than $m+2$,
so there are at least $m!$ choices for the transpositions producing the real tropical cover $(\varphi,\rho_s)$.
Therefore, we obtain that the number $\vec{N}(\varphi,s)\geq m!$.
\end{proof}

Denote by $\lambda_e^{max}$ the maximal even integer in $\lambda$.
We put
$$
\vec h^\rb_{g,\lambda,\mu}(m):=\vec{H}^\rb_g((\lambda,2,1^{2m}),(\mu,2,1^{2m})).
$$

\begin{theorem}
\label{thm:asym-simple-splitting}
Let $\lambda$ and $\mu$ be two partitions with $|\lambda|=|\mu|$.
Suppose that $\mu_{o,o}=\emptyset$, $l(\lambda_{o,o})>1$ and no odd integer other than $1$ appears in $\lambda_{o,o}$.
Additionally, if $\lambda$ and $\mu$ satisfy one of the following four conditions:
\begin{itemize}
    \item $l(\lambda_o)=l(\mu_o)=1$ and $|\lambda_e^{max}|>|\mu_o|$;
    \item $l(\lambda_o)=2$, $l(\mu_o)=0$ and $|\mu_e^{max}|>\max(\lambda_{o_1},\lambda_{o_2})$, where $\lambda_{o_1}$, $\lambda_{o_2}\in\lambda_o$;
    \item $l(\lambda_o)=0$, $l(\mu_o)=2$ and $|\lambda_e^{max}|>\max(\mu_{o_1},\mu_{o_2})$, where $\mu_{o_1}$, $\mu_{o_2}\in\mu_o$;
    \item $l(\lambda_o)=l(\mu_o)=0$;
\end{itemize}
there is at least one monotone zigzag cover of type $(g,\lambda,\mu,\undl x)$,
and the logarithmic asymptotics for $\vec h^\rb_{g,\lambda,\mu}(m)$ is at least $2m\log m$ as $m\to\infty$.
\end{theorem}

\begin{proof}
We prove Theorem $\ref{thm:asym-simple-splitting}$ by considering
the following four cases:

{\bf Case $(1)$:} $l(\lambda_o)=l(\mu_o)=1$ and $|\lambda_e^{max}|>|\mu_o|$.

We first show that there is a monotone zigzag cover of type $(g,\lambda,\mu,\undl x)$.
Since $|\lambda|=|\mu|=d$, the condition
$|\lambda_e^{max}|>|\mu_o|$ implies that $|\lambda|-|\lambda_e^{max}|<|\mu|-|\mu_o|$.
Suppose that $\mu_e^s\in\mu_{e}$, where $s=1,\ldots,l(\mu)-1$, and $\lambda_e^t\in(\lambda_{e},(2)^{l(\lambda_{1,1})})$, where $t=1,\ldots,l(\lambda)-l(\lambda_{1,1})-1$.
We consider sequences of integers $(k_0,k_1,k_2,\ldots,k_{N})$ defined by the following procedure.
Let $k_0=\lambda_o$ and $k_1=k_0-\mu_e^1$.
For any $i\in\{1,\ldots,N\}$, $k_{i+1}=k_i-\mu_e^{s_i}$ if $k_i>0$,
otherwise $k_{i+1}=k_i+\lambda_e^{t_i}$.
Note that every element in $\mu_e$ and
$(\lambda_{e},(2)^{l(\lambda_{1,1})})$ is used exactly once here
and $\lambda_e^{max}$ is used only when all elements in
$(\lambda_{e},(2)^{l(\lambda_{1,1})})\setminus\lambda_e^{max}$ are used up.
At the end of this procedure, we obtain that
$k_N=\mu_o$ and $k_N\neq k_{N-1}+\lambda_e^{max}$ with $k_{N-1}>0$.

Now we construct a monotone zigzag cover of type $(g,\lambda,\mu,\undl x)$.
We first choose a string $S$ with two ends,
then label the in-end with $\lambda_o$
and label the out-end with $\mu_o$.
The tails depicted in Figure $\ref{fig:zigzag}$ are attached to the string $S$ following the above sequence of integers $k_0,k_1,\ldots,k_N$.
Namely, if $k_{i+1}=k_i-\mu_e^{s_i}$ and $k_{i+1}>0$, a unbent out-tail of weight $\mu_e^{s_i}$ is attached to $S$ corresponding to $k_{i+1}$.
A bent out-tail of weight $\mu_e^{s_j}$ is attached to $S$ corresponding to $k_{j+1}<0$ when $k_{j+1}=k_j-\mu_e^{s_j}$.
Similarly, in the case $k_{u+1}=k_u+\lambda_e^{s_u}$ and $k_{u+1}<0$,
a unbent in-tail of weight $\lambda_e^{s_u}$ is attached to $S$ corresponding to $k_{u+1}$.
A bent in-tail of weight $\lambda_e^{s_v}$ is attached to $S$ corresponding to a $k_{v+1}$ when $k_{v+1}=k_v+\lambda_e^{s_v}$ and $k_{v+1}>0$. Note that if $\lambda_e^s\in(2)^{l(\lambda_{1,1})}$,
the in-tail attached to $S$ corresponding to $\lambda_e^s$ is an in-tail with symmetric fork, and we place $g$ symmetric cycles on exactly one such in-tail.
See Figure $\ref{fig:exa-mon-zigzag-1}$ for an example.
\begin{figure}[ht]
    \centering
    \begin{tikzpicture}
    \draw[line width=0.3mm,gray] (-5,2)--(3,1.5)--(-2,1)--(0,0.5)--(-4.5,0)--(4,-0.5);
    \draw[line width=0.3mm] (3,1.5)--(4,1.5);
    \draw[line width=0.3mm] (1.7,1.59)--(4,1.59);
    \draw[line width=0.3mm] (-5,1.3)--(1,1.3);
    \draw[line width=0.3mm] (-2.2,1)--(-2,1);
    \draw[line width=0.3mm] (-5,1.2)--(-2.2,1)--(-5,0.8);
    \draw[line width=0.3mm] (-1,0.75)--(4,0.75);
    \draw[line width=0.3mm] (0,0.5)--(0.5,0.5);
    \draw[line width=0.3mm] (4,0.65)--(0.5,0.5)--(4,0.35);
    \draw[line width=0.3mm] (-5,0.2)--(-2.7,0.2);
    \draw[line width=0.3mm] (-5,0)--(-4.5,0);
    \draw[line width=0.3mm] (-3,-0.09)--(4,-0.09);
    \draw[line width=0.3mm] (-5,-0.75)--(4,-0.75);
    \foreach \Point in {(-4.5,-0.75), (-3,-0.75), (-2.7,-0.75), (-2.2,-0.75), (-2,-0.75), (-1,-0.75), (0,-0.75), (0.5,-0.75), (1,-0.75), (1.7,-0.75), (3,-0.75)}
    \draw[fill=black] \Point circle (0.05);
    \draw[line width=0.4mm] (-4.5,-1)  node{\tiny{$x_1$}} (-3,-1) node{\tiny{$x_2$}} (-2.7,-1) node{\tiny{$x_3$}} (-2.2,-1) node{\tiny{$x_4$}} (-1.9,-1) node{\tiny{$x_5$}} (-1,-1) node{\tiny{$x_6$}} (0,-1) node{\tiny{$x_7$}} (0.5,-1) node{\tiny{$x_8$}} (1,-1) node{\tiny{$x_9$}} (1.7,-1) node{\tiny{$x_{10}$}} (3,-1) node{\tiny{$x_{11}$}};
    \end{tikzpicture}
    \caption{Monotone zigzag cover $\varphi:C\to T\pb^1$ when $l(\lambda_o)=l(\mu_o)=1$.}
    \label{fig:exa-mon-zigzag-1}
\end{figure}
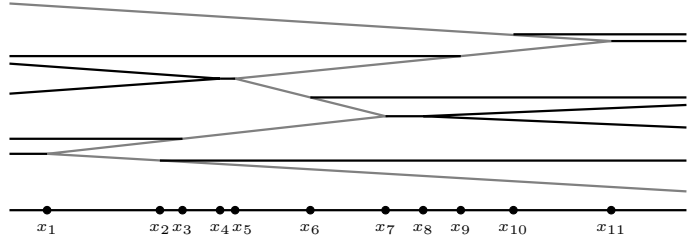
An inner edge $E$ in the string $S$ with two vertices determined by tails corresponding to $k_{i+1}$ and $k_i$ is weighted by $|k_i|$.
The edge $E$ is oriented from the vertex corresponding to $k_i$ to the vertex corresponding to $k_{i+1}$ if $k_i>0$.
Otherwise, $E$ is oriented from the vertex corresponding to $k_{i+1}$ to the vertex corresponding to $k_i$.
We obtain a graph $C$ satisfying the balancing condition, so it is a tropical curve.
Moreover, the curve $C$ satisfies conditions $(1)$ and $(2)$ of Definition $\ref{def:mono-zigzag}$.
The orientation on the edges of $C$ induces a partial order on the set of inner vertices.
By shrinking or extending the lengths of inner edges of $C$ properly,
one obtain a total order on the set of inner vertices compatible with the partial order such that $\varphi:C\to T\pb^1$ satisfies condition $(3)$ of Definition $\ref{def:mono-zigzag}$ with this total order.
At last, we get a monotone zigzag cover $\varphi:C\to T\pb^1$ of type $(g,\lambda,\mu,\undl x)$.

Now we prove the asymptotics for real monotone Hurwitz numbers relative to simple splittings when only simple branch points are added as the degree goes to infinity.
Let $\varphi:C\to T\pb^1$ be the monotone
zigzag cover of type $(g,\lambda,\mu,\undl x)$ constructed above.
Recall that $k_N=\mu_o$ and $k_N\neq k_{N-1}+\lambda_e^{max}$ with $k_{N-1}>0$,
so the last piece $S_n$ of $S$ is either attached with unbent out-tails or not attached with any unbent tail.
In any case, we attach $a$, $a\leq\frac{\lambda_e^{max}-1}{2}$, unbent in-tails of weight $2$ with symmetric forks to piece $S_{n-1}$ such that the edge $E'$ in $S_{n-1}$ with a bent vertex $v$ with peak pointing to the left has weight $1$.
Since some unbent in-tails are attached to $S_{n-1}$, the balancing condition at the bent vertex $v$ is broken.
We attach $a$ unbent out-tails of weight $2$ with symmetric forks to piece $S_{n}$ to compensate the weight, and the balancing condition is preserved at $v$ (See Figure $\ref{fig:asym-monzigzag-1}$ for a local picture of $v$).
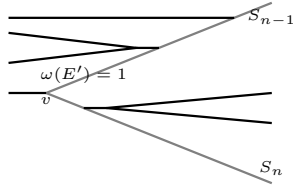
\begin{figure}[ht]
    \centering
    \begin{tikzpicture}
    \draw[line width=0.3mm,gray] (3,1.2)--(0,0)--(3,-1.2);
    \draw[line width=0.3mm] (-0.5,0)--(0,0);
    \draw[line width=0.3mm] (-0.5,1)--(2.5,1);
    \draw[line width=0.3mm] (1.2,0.6)--(1.5,0.6);
    \draw[line width=0.3mm] (-0.5,0.4)--(1.2,0.6)--(-0.5,0.8);
    \draw[line width=0.3mm] (0.5,-0.2)--(0.8,-0.2);
    \draw[line width=0.3mm] (3,0)--(0.8,-0.2)--(3,-0.4);
    \draw (0,-0.1) node{\tiny{$v$}} (0.5,0.25) node{\tiny{$\omega(E')=1$}} (3,1) node{\tiny{$S_{n-1}$}} (3,-1) node{\tiny{$S_n$}};
    \end{tikzpicture}
    \caption{A local picture of the intersection point $v$ of $S_{n-1}$ and $S_n$.}
    \label{fig:asym-monzigzag-1}
\end{figure}
Next, we cut at the middle of edge $E'$ and obtain two remaining parts $C'$ and $C''$ of $C$.
Note that $\varphi|_{C'}:C'\to T\pb^1$ is a zigzag cover with only one bent vertex $v$ in $C'$ and $\varphi|_{C''}:C''\to T\pb^1$ is a monotone zigzag covers.
Let $\tilde\varphi:\tilde C\to T\pb^1$ be a tropical cover of type $(0,(2,1^{2m-2a+1}),(2,1^{2m-2a+1}),\undl{\tilde x})$
obtained in Proposition $\ref{prop:asym-simple}$.
From the proof of Proposition $\ref{prop:asym-simple}$, the cover $\tilde\varphi:\tilde C\to T\pb^1$ can be chosen such that the monotone components $C_1$ and $C_{m-a}$ are of type $(1)$ and type $(3)$ depicted in Figure $\ref{fig:mon-comps}$, respectively.
See Figure $\ref{fig:mon-tropical-modified}(1)$ for an example of $\tilde\varphi:\tilde C\to T\pb^1$.
We glue the remaining half-edges of $E'$ in $C'$ and $C''$ to the ends of the string of the monotone components $C_1$ and $C_{m-a}$ of the graph $\tilde C$.
We get a tropical cover $\bar\varphi$ of type $(g,(\lambda,2,1^{2m}),(\mu,2,1^{2m}),\undl{\bar x})$ (See Figure $\ref{fig:asym-glued-zigzag-1}$ for an example).
\begin{figure}[ht]
    \centering
    \begin{tikzpicture}
    \draw[line width=0.3mm,gray] (2.1,2)--(3.5,1.8)--(2.6,1.6)--(3,1.4)--(2.1,1);
    \draw[line width=0.3mm] (3.5,1.8)--(4,1.8);
    \draw[line width=0.3mm] (2.1,1.6)--(2.6,1.6);
    \draw[line width=0.3mm] (3,1.4)--(4,1.4);
    \draw[line width=0.3mm] (2.8,1.5)--(4,1.5);
    \draw[line width=0.3mm] (2.25,1.18)--(2.5,1.18);
    \draw[line width=0.3mm] (2.1,1.08)--(2.25,1.18)--(2.1,1.28);
    \draw[line width=0.3mm,gray] (-4,2)--(1.5,2);
    \draw[line width=0.3mm,gray,dotted] (1.5,2)--(1.9,2);
    \draw[line width=0.3mm] (-4,1.4)--(1.5,1.4);
    \draw[line width=0.3mm,dotted] (1.5,1.4)--(1.9,1.4);
    \draw[line width=0.3mm,gray] (1.9,1)--(-1.2,0.7)--(-0.8,0.6)--(-1.9,0.3);
    \draw[line width=0.3mm] (-1.6,0.7)--(-1.2,0.7);
    \draw[line width=0.3mm] (-1.9,0.8)--(-1.6,0.7)--(-1.9,0.6);
    \draw[line width=0.3mm] (-0.8,0.6)--(-0.6,0.6);
    \draw[line width=0.3mm,gray] (1.9,0.7)--(-0.6,0.6)--(1.2,0.2)--(0.8,0)--(1.9,-0.2);
    \draw[line width=0.3mm] (1.2,0.2)--(1.6,0.2);
    \draw[line width=0.3mm] (1.9,0.3)--(1.6,0.2)--(1.9,0.1);
    \draw[line width=0.3mm] (0.6,0)--(0.8,0);
    \draw[line width=0.3mm,gray] (-1.9,0.1)--(0.6,0)--(-0.4,-0.3)--(0.4,-0.6)--(-1.9,-0.9);
    \draw[line width=0.3mm] (-1.9,-0.3)--(-0.4,-0.3);
    \draw[line width=0.3mm] (0.4,-0.6)--(1.9,-0.6);
    \draw[line width=0.3mm,gray] (-4,0.8)--(-2.6,0.8);
    \draw[line width=0.3mm,gray,dotted] (-2.6,0.8)--(-2.1,0.8);
    \draw[line width=0.3mm] (-4,-0.6)--(-2.6,-0.6);
    \draw[line width=0.3mm,dotted] (-2.6,-0.6)--(-2.1,-0.6);
    \draw[line width=0.3mm,gray] (4,0.8)--(2.6,0.8);
    \draw[line width=0.3mm,gray,dotted] (2.6,0.8)--(2.1,0.8);
    \draw[line width=0.3mm] (4,-0.6)--(2.6,-0.6);
    \draw[line width=0.3mm,dotted] (2.6,-0.6)--(2.1,-0.6);
    \draw[line width=0.3mm,gray] (-2.1,-0.9)--(-3.5,-1.5)--(-2.1,-2);
    \draw[line width=0.3mm] (-4,-1.5)--(-3.5,-1.5);
    \draw[line width=0.3mm] (-2.8,-1.75)--(-2.4,-1.75);
    \draw[line width=0.3mm] (-2.1,-1.65)--(-2.4,-1.75)--(-2.1,-1.85);
    \draw[line width=0.3mm,gray] (4,-2)--(-1.5,-2);
    \draw[line width=0.3mm,gray,dotted] (-1.5,-2)--(-1.9,-2);
    \draw[line width=0.3mm] (4,-1.4)--(-1.5,-1.4);
    \draw[line width=0.3mm,dotted] (-1.5,-1.4)--(-1.9,-1.4);
    \draw[line width=0.3mm] (-3.5,-1.7) node{\tiny{$v$}} (-2.2,-0.8) node{\tiny{$1$}} (2.2,0.9) node{\tiny{$1$}} (-1.8,-0.8) node{\tiny{$1$}} (1.8,0.9) node{\tiny{$1$}};
    \draw[line width=0.3mm] (-3,-2) node{\tiny{$C'$}} (0,-0.8) node{\tiny{$\tilde C$}} (3,1.2) node{\tiny{$C''$}};
    \draw[line width=0.3mm,gray,dashed] (-2,2.1)--(-2,-2.4);
    \draw[line width=0.3mm,gray,dashed] (2,2.1)--(2,-2.4);
    \draw[line width=0.3mm] (-4,-2.3)--(4,-2.3);
    \foreach \Point in {(-3.5,-2.3), (-2.8,-2.3), (-2.4,-2.3), (-1.6,-2.3), (-1.2,-2.3), (-0.8,-2.3), (-0.6,-2.3), (-0.4,-2.3), (0.4,-2.3), (0.6,-2.3), (0.8,-2.3), (1.2,-2.3), (1.6,-2.3), (2.25,-2.3), (2.5,-2.3), (2.6,-2.3), (2.8,-2.3), (3,-2.3), (3.5,-2.3)}
    \draw[fill=black] \Point circle (0.05);
    \draw[line width=0.4mm] (-3,-2.5)  node{\tiny{$\undl x'$}} (0,-2.5) node{\tiny{$\tilde{\undl x}$}} (3,-2.5) node{\tiny{$\undl x''$}};
    \end{tikzpicture}
    \caption{An example of the glued tropical cover $\bar\varphi: \bar C\to T\pb$.}
    \label{fig:asym-glued-zigzag-1}
\end{figure}
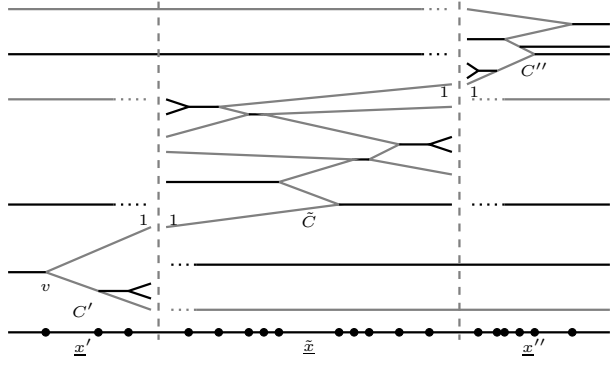
Let $\undl{\bar x}=\undl {\bar x}^+\sqcup\undl{\bar x}^-$ be a simple splitting with $|\undl{\bar x}^+|=s$.
This splitting induces a simple splitting on each of the three sets $\undl x'$, $\undl{\tilde x}$, $\undl x''$.
Suppose that the induced splitting on $\undl{\tilde x}$ is $\undl{\tilde x}=\undl{\tilde x}^+\sqcup\undl{\tilde x}^-$ with $|\undl{\tilde x}^+|=\tilde s$.
From Proposition $\ref{prop:asym-simple}$,
the tropical cover $\tilde\varphi:\tilde C\to T\pb^1$ can be chosen to be real under the corresponding splitting on $\undl{\tilde x}$.
Let $z_1$ (resp. $z_2$) be any integer satisfying $1\leq z_1\leq\lambda_e^{max}$ (resp. $\lambda_e^{max}<z_2\leq\lambda_e^{max}+2m-2a+2$).
We use cycles of $\sal_d$ with elements in $\{1,2,\ldots,\lambda_e^{max}\}$, $\{z_1,\lambda_e^{max}+1,\lambda_e^{max}+2,\ldots,\lambda_e^{max}+2m-2a+2\}$ and
$\{z_2,\lambda_e^{max}+2m-2a+3,\lambda_e^{max}+2m-2a+2,\ldots,|\lambda|+2m+2\}$ to produce real tropical covers $\varphi|_{C'}$, $\tilde\varphi$ and $\varphi|_{C''}$ with the induced splittings according to Construction
$\ref{const2}$, respectively.
The integers $z_1$ and $z_2$ are used to produce the bridge edges of the three parts $C'$, $\tilde C$ and $C''$ respectively.
Once we choose two integers, say $z_1$ and $z_2$, as the integers corresponding to the out-end of the string in $C'$ and the out-end of the string in $C_{m-a}$ in $\tilde C$ respectively, the composition of the three sets of cycles of $\sal_d$ in the above produces the real tropical cover $\bar\varphi$ with the splitting $\undl{\bar x}=\undl{\bar x}^+\sqcup\undl{\bar x}^-$.
From Lemma $\ref{lem:mono-zigzag1}$ and Proposition $\ref{prop:asym-simple}$, the number $\vec{N}(\varphi|_{C''})>0$ and $\vec N(\tilde\varphi,\tilde s)\geq(m-\lambda_e^{max})!$.
Note that $C'$ contains only one bent vertex and all the unbent tails attached to $C'$ are out-tails.
By a similar argument as in the proof of Lemma $\ref{lem:mono-zigzag1}$, one obtains that the number $\vec{N}(\varphi|_{C'})>0$.
Hence, we obtain $\vec N(\bar\varphi,s)\geq(m-\lambda_e^{max})!$.
It follows from Proposition $\ref{prop:asym-simple}$ and the above construction that there are at least $(m-\lambda_e^{max})!$ real tropical covers $(\bar\varphi,\bar\rho_s)$ with $\vec N(\bar\varphi,s)\geq(m-\lambda_e^{max})!$.
The number of real monotone Hurwitz numbers relative to simple splittings is bounded from below by:
$$
\begin{aligned}
\vec h^\rb_{g,\lambda,\mu}(m)&\geq\sum_{(\bar\varphi,\bar\rho_s)}N(\bar\varphi,s)\\
&\geq(m-\lambda_e^{max})!^2.
\end{aligned}
$$
Since $\log[(m-\lambda_e^{max})!^2]\sim2\log m$ as $m\to\infty$,
we obtain the expected logarithmic asymptotics for $\vec h^\rb_{g,\lambda,\mu}(m)$.

{\bf Case $(2)$:} $l(\lambda_o)=2$, $l(\mu_o)=0$ and $\mu_e^{max}>\max(\lambda_{o_1},\lambda_{o_2})$.

Since the idea of the proof for this case is the same as in case $(1)$,
we only give a sketch here.
As in the case $(1)$,
we denote by $\mu_e^{max}$ the maximal even integer in $\mu$.
The condition
$\mu_e^{max}>\max(\lambda_{o_1},\lambda_{o_2})$ implies that $|\lambda|-|\lambda_{o_2}|>|\mu|-|\mu_e^{max}|$.
Let $\mu_e^s\in\mu_{e}$ and $\lambda_e^t\in(\lambda_{e},(2)^{l(\lambda_{1,1})})$ be the same as in case $(1)$.
The sequence of integers  $(k_0,k_1,k_2,\ldots,k_{N})$ that we consider here is defined as follows.
Let $k_0=\lambda_{o_1}$ and $k_1=k_0-\mu_e^1$.
For any $i\in\{1,\ldots,N\}$, $k_{i+1}=k_i-\mu_e^{s_i}$ if $k_i>0$,
otherwise $k_{i+1}=k_i+\lambda_e^{t_i}$.
Note that every element in $\mu_{e}$ and
$(\lambda_{e},(2)^{l(\lambda_{1,1})})$ is used exactly once here
and $\mu_e^{max}$ is used only when all elements in
$\mu_{e}\setminus\mu_e^{max}$ are used up.
In this case, we obtain that
$k_N=-\lambda_{o_2}$ and $k_N\neq k_{N-1}-\mu_e^{max}$ with $k_{N-1}<0$.

A monotone zigzag cover $\varphi:C\to T\pb^1$ of type $(g,\lambda,\mu,\undl x)$ is constructed in the same way as that in case $(1)$.
Note that the string in $C$ has two in-ends,
and the last piece $S_n$ is not attached with unbent out-tails.
We attach $a$, $a\leq\frac{\lambda_{o_2}-1}{2}$, unbent in-tails of weight $2$ with symmetric fork and $a$ unbent out-tails of weight $2$ with symmetric fork to the weight $\lambda_{o_2}$ in-end of the string such that the inner edge $E'$ connecting the weight $2$ in-tail and out-tail is of weight $1$ (See Figure $\ref{fig:asym-monzigzag-2}$ for an example).
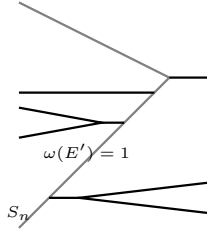
\begin{figure}[ht]
    \centering
    \begin{tikzpicture}
    \draw[line width=0.3mm,gray] (0,1)--(2,0)--(0,-2);
    \draw[line width=0.3mm] (2,0)--(2.5,0);
    \draw[line width=0.3mm] (1.8,-0.2)--(0,-0.2);
    \draw[line width=0.3mm] (1.4,-0.6)--(1.1,-0.6);
    \draw[line width=0.3mm] (0,-0.4)--(1.1,-0.6)--(0,-0.8);
    \draw[line width=0.3mm] (0.4,-1.6)--(0.8,-1.6);
    \draw[line width=0.3mm] (2.5,-1.4)--(0.8,-1.6)--(2.5,-1.8);
    \draw (0.9,-1) node{\tiny{$\omega(E')=1$}} (0,-1.8) node{\tiny{$S_n$}};
    \end{tikzpicture}
    \caption{A local picture of the edge $E'$.}
    \label{fig:asym-monzigzag-2}
\end{figure}
As in case $(1)$, we cut at the middle of edge $E'$ and glue the half edges with the two ends of a type $(0,(2,1^{2m-2a+1}),(2,1^{2m-2a+1}),\undl{\tilde x})$ tropical cover obtained in Proposition $\ref{prop:asym-simple}$,
then we have a tropical cover $\bar\varphi:\bar C\to T\pb^1$ of type $(g,(\lambda,2,1^{2m}),(\mu,2,1^{2m}),\undl{\bar x})$.
By a similar argument as in the case $(1)$, we have at least $(m-\mu_e^{max})!$ real tropical covers $(\bar\varphi,\bar\rho_s)$ with $\vec N(\bar\varphi)\geq(m-\mu_e^{max})!$,
and hence completes the proof.

{\bf Case $(3)$:} $l(\lambda_o)=0$, $l(\mu_o)=2$ and $\lambda_e^{max}>\max(\mu_{o_1},\mu_{o_2})$.

As in the case $(2)$, we sketch the proof of Theorem $\ref{thm:asym-simple-splitting}$ for the case $l(\lambda_o)=0$ and $l(\mu_o)=2$ as follows.
Because $|\lambda|=|\mu|=d$ and $\lambda_e^{max}>\max(\mu_{o_1},\mu_{o_2})$,
we have $|\mu|-|\mu_{o_2}|>|\lambda|-|\lambda_e^{max}|$.
Let $\mu_e^s\in\mu_{e}$ and $\lambda_e^t\in(\lambda_{e},(2)^{l(\lambda_{1,1})})$ be the same as case $(1)$.
The sequence of integers  $(k_0,k_1,k_2,\ldots,k_{N})$ is defined as follows.
Let $k_0=-\mu_{o_1}$ and $k_1=k_0+\lambda_e^1$.
For any $i\in\{1,\ldots,N\}$, $k_{i+1}=k_i-\mu_e^{s_i}$ if $k_i>0$,
otherwise $k_{i+1}=k_i+\lambda_e^{t_i}$.
In this procedure, every element in $\mu_{e}$ and
$(\lambda_{e},(2)^{l(\lambda_{1,1})})$ is used exactly once.
Moreover, the integer $\lambda_e^{max}$ is used only when all elements in
$(\lambda_{e},(2)^{l(\lambda_{1,1})})\setminus\lambda_e^{max}$ are used up.
Since $\lambda_e^{max}>\max(\mu_{o_1},\mu_{o_2})$ and $\lambda_e^{max}$ is the integer we used last to add to certain $k_i$,
$k_N=\mu_{o_2}$ and $k_N\neq k_{N-1}+\lambda_e^{max}$ with $k_{N-1}>0$.
As in the case $(2)$, we construct a monotone zigzag cover $\varphi:C\to T\pb^1$ of type $(g,\lambda,\mu,\undl x)$ and
the string in $C$ has two out-ends.
The last piece $S_n$ is not attached with unbent in-tails.
We attach some out-tails and in-tails to the weight $\mu_{o_2}$ out-end of the string in $C$ such that the inner edge $E'$ connecting the out-tail and in-tail is of weight $1$ (See Figure $\ref{fig:asym-monzigzag-3}$ for an example).
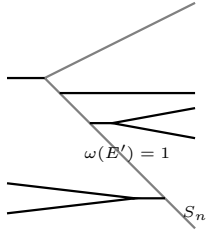
\begin{figure}[ht]
    \centering
    \begin{tikzpicture}
    \draw[line width=0.3mm,gray] (0,1)--(-2,0)--(0,-2);
    \draw[line width=0.3mm] (-2,0)--(-2.5,0);
    \draw[line width=0.3mm] (-1.8,-0.2)--(0,-0.2);
    \draw[line width=0.3mm] (-1.4,-0.6)--(-1.1,-0.6);
    \draw[line width=0.3mm] (0,-0.4)--(-1.1,-0.6)--(0,-0.8);
    \draw[line width=0.3mm] (-0.4,-1.6)--(-0.8,-1.6);
    \draw[line width=0.3mm] (-2.5,-1.4)--(-0.8,-1.6)--(-2.5,-1.8);
    \draw (-0.9,-1) node{\tiny{$\omega(E')=1$}} (0,-1.8) node{\tiny{$S_n$}};
    \end{tikzpicture}
    \caption{$C$ and $\tilde C$ type $(0,((1)^{2m},\mu_{o_2}),((1)^{2m},\mu_{o_2}),\undl{\tilde x})$}
    \label{fig:asym-monzigzag-3}
\end{figure}
We proceed a surgery at the edge $E'$ similar to that in the case $(1)$
(See Figure $\ref{fig:asym-glued-zigzag-1}$ for more details).
Then we get a real tropical cover $(\bar\varphi,\bar\rho)$ of type $(g,(\lambda,2,1^{2m}),(\mu,2,1^{2m}),\undl{\bar x})$.
By a similar argument as in the case $(1)$, we have at least $(m-\lambda_e^{max})!$ real tropical covers $(\bar\varphi,\bar\rho)$ with $\vec N(\bar\varphi,s)\geq(m-\lambda_e^{max})!$.
Therefore, we obtain the logarithmic growth of $\vec h^\rb_{g,\lambda,\mu}(m)$ as $m\to\infty$.

{\bf Case $(4)$:} $l(\lambda_o,\mu_o)=0$.

The construction in this case is similar to that in the case $(2)$,
that is $l(\lambda_o)=2$ and $l(\mu_o)=0$,
so we only sketch how to define the sequence of integers $k_0,k_1,\ldots,k_N$.
Since there is no odd integer that appears an odd times in $\lambda$ and $\mu$,
we use a pair of ones in $\lambda_{1,1}$ to replace
the two odd integers $\lambda_{o_1},\lambda_{o_2}$ in the case $(2)$.
Suppose that $\mu_e^s\in\mu_{e}$, where $s=1,\ldots,l(\mu)$, and $\lambda_e^t\in(\lambda_{e},(2)^{l(\lambda_{1,1})-1})$, where $t=1,\ldots,l(\lambda)-l(\lambda_{1,1})-1$.
Let $k_0=1$ and $k_1=k_0-\mu_e^1$.
For any $i\in\{1,\ldots,N\}$, $k_{i+1}=k_i-\mu_e^{s_i}$ if $k_i>0$,
otherwise $k_{i+1}=k_i+\lambda_e^{t_i}$.
Every element in $\mu_{e}$ and
$(\lambda_{e},(2)^{l(\lambda_{1,1})-1})$ is used exactly once here.
Note that $k_N=-1$ and $k_N\neq k_{N-1}-\lambda_e^{s}$ with $k_{N-1}<0$.
The rest of the proof is the same as that for the case $(2)$, so we omit it.
\end{proof}

\begin{corollary}
\label{cor:asym-simple-splitting}
Under the same conditions of Theorem $\ref{thm:asym-simple-splitting}$,
real monotone double Hurwitz numbers $\vec h^\rb_{g,\lambda,\mu}(m)$ relative to simple splittings are logarithmically equivalent to the monotone double Hurwitz numbers $\vec h^\cb_{g,\lambda,\mu}(m)$:
$$
\log\vec h^\rb_{g,\lambda,\mu}(m)\sim2m\log m\sim\log\vec h^\cb_{g,\lambda,\mu}(m),
\text{ as }m\to\infty,
$$
where $\vec h^\cb_{g,\lambda,\mu}(m)=\vec H^\cb_g((\lambda,2,1^{2m}),(\mu,2,1^{2m}))$.
\end{corollary}

\begin{proof}
It follows from \cite[Theorem $1.1$]{ggpn-2017} or \cite[Theorem $4.4$]{ggpn-2014} that $\vec h^\cb_{g,\lambda,\mu}(m)\sim2m\log m$ as $m\to\infty$.
It is obviously that $\vec h^\rb_{g,\lambda,\mu}(m)\leq\vec h^\cb_{g,\lambda,\mu}(m)$,
so Corollary $\ref{cor:asym-simple-splitting}$ follows from Theorem $\ref{thm:asym-simple-splitting}$.
\end{proof}

\subsection{Case $\textrm{II}$: asymptotics for arbitrary splitting}
We study the logarithmic asymptotics for real monotone double Hurwitz numbers relative to arbitrary splittings.

\begin{theorem}\label{thm:asymp-arbitrary-splitting}
Let $\lambda$ and $\mu$ be two partitions with $|\lambda|=|\mu|$.
Suppose that $\mu_{o,o}=\emptyset$, $l(\lambda_{o,o})>1$ and no odd integer other than $1$ appears in $\lambda_{o,o}$.
Assume that $\mu_i+\mu_j>\lambda_e^{max}$ for any two parts $\mu_i$, $\mu_j$ of $\mu_{e}$.
Additionally, if $\lambda$ and $\mu$ satisfy one of the following four conditions:
\begin{enumerate}
    \item[$(1)$] $\lambda_o=\mu_o=(1)$;
    \item[$(2)$] $l(\lambda_o)=2$, $l(\mu_o)=0$ and $\lambda_{o_1}\neq1,\lambda_{o_2}=1$, where $\lambda_{o_1},\lambda_{o_2}\in\lambda_o$;
    \item[$(3)$] $l(\lambda_o)=0$, $l(\mu_o)=2$ and $\mu_{o_1}\neq1,\mu_{o_2}=1$, where $\mu_{o_1},\mu_{o_2}\in\mu_o$;
    \item[$(4)$] $l(\lambda_o)=l(\mu_o)=0$;
\end{enumerate}
there is at least one universally monotone zigzag cover of type $(g,\lambda,\mu, \undl x)$
and the universally monotone zigzag number $\vec\zl_{g}((\lambda, 1^{2m}),(\mu, 1^{2m}))\geq m!$.
\end{theorem}

\begin{proof}
Since $\mu_i+\mu_j>\lambda_e^{max}$ for any two parts $\mu_i$, $\mu_j$ of $\mu_{e}$,
the construction of monotone zigzag cover of type $(g,\lambda,\mu,\undl x)$ in the proof of Theorem $\ref{thm:asym-simple-splitting}$ produces a universally monotone zigzag cover of type $(g,\lambda,\mu,\undl x)$.
Therefore, we omit the construction of universally monotone zigzag cover of type $(g,\lambda,\mu,\undl x)$ here.
Now we show that the universally monotone zigzag number $\vec\zl_{g}((\lambda, 1^{2m}),(\mu, 1^{2m}))$ is bounded from below as follows.

{\bf Case $(1)$:} $\lambda_o=\mu_o=(1)$.

Let $\varphi:C\to T\pb^1$ be a universally monotone
zigzag cover of type $(g,\lambda,\mu,\undl x)$ constructed in the above.
Note that in this case the out-end of the string $S$ of the universally monotone zigzag cover is of weight $1$.
We glue the out-end of the string $S$ with the in-end of the string of the universally monotone zigzag cover $\tilde\varphi:\tilde C\to T\pb^1$ of type $(0,(1^{2m+1}),(1^{2m+1}),\undl{\tilde x})$ (See Figure $\ref{fig:mon-zigzag1}$ for an example).
By shrinking or extending the length of the glued inner edge properly,
we get a universally monotone zigzag cover $\bar\varphi$ of type $(g,(\lambda,(1)^{2m}),(\mu,(1)^{2m}),\undl{\bar x})$.
Let $\undl{\bar x}=\undl{\bar x}^+\sqcup\undl{\bar x}^-$ be an arbitrary splitting with $|\undl{\bar x}^+|=s$.
This splitting induces a splitting of each of the two sets $\undl x$, $\undl{\tilde x}$.
Let $z$ be an undetermined integer satisfying $2m+1\leq z\leq|\lambda|+2m$.
We use cycles in $\sal_d$ with elements in $\{1,2,\ldots,2m,z\}$, $\{2m+1,2m+2,\ldots,|\lambda|+2m\}$ to produce universally real monotone zigzag covers $\tilde\varphi$, $\varphi$ with the induced splittings according to Construction    $\ref{const2}$, respectively.
The integer $z$ is used to produce the in-end of the string of $\tilde C$.
Once we choose $z$ as the integer corresponding to the out-end of the string in $C$, the composition of the two sets of cycles in the above produces universally real monotone zigzag cover $\bar\varphi$ with the splitting $\undl{\bar x}=\undl{\bar x}^+\sqcup\undl{\bar x}^-$.
From Lemma $\ref{lem:mono-zigzag2}$ and Proposition $\ref{prop:nonvanish-mon-zigzag}$, the numbers $\vec{N}(\varphi)>0$ and $\vec N(\tilde\varphi)\geq m!$,
so the number $\vec N(\bar\varphi)\geq m!$.
Hence, we obtain the lower bound of $\vec\zl_{g}((\lambda, 1^{2m}),(\mu, 1^{2m}))\geq m!$.

{\bf Case $(2)$:} $l(\lambda_o)=2$, $l(\mu_o)=0$ and $\lambda_{o_1}\neq1,\lambda_{o_2}=1$.

Let $\varphi:C\to T\pb^1$ be a universally monotone zigzag cover of type $(g,\lambda,\mu,\undl x)$.
Note that the string in $C$ has two in-ends,
and the in-end of the last piece $S_n$ is of weight $1$.
We glue the weight $1$ in-end of the string in $C$ with the out-end of the string of a universally monotone zigzag cover $\tilde\varphi:\tilde C\to T\pb^1$ of type $(0,(1^{2m+1}),(1^{2m+1}),\undl{\tilde x})$ (See Figure $\ref{fig:asym-uni-monzigzag2}$ for an example).
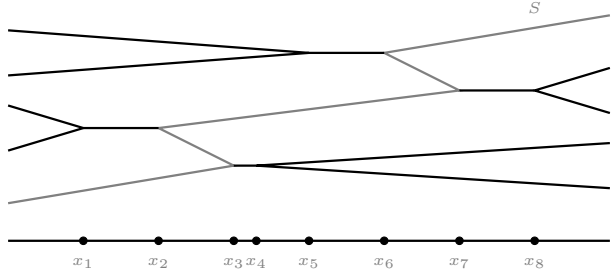
\begin{figure}[ht]
    \centering
    \begin{tikzpicture}
    \draw[line width=0.3mm,gray] (4,-0.5)--(1,-1)--(2,-1.5)--(-2,-2)--(-1,-2.5)--(-4,-3);
    \draw[line width=0.3mm] (-2,-2)--(-3,-2);
    \draw[line width=0.3mm] (-4,-1.7)--(-3,-2)--(-4,-2.3);
    \draw[line width=0.3mm] (-1,-2.5)--(-0.7,-2.5);
    \draw[line width=0.3mm] (4,-2.8)--(-0.7,-2.5)--(4,-2.2);
    \draw[line width=0.3mm] (0,-1)--(1,-1);
    \draw[line width=0.3mm] (-4,-0.7)--(0,-1)--(-4,-1.3);
    \draw[line width=0.3mm] (2,-1.5)--(3,-1.5);
    \draw[line width=0.3mm] (4,-1.8)--(3,-1.5)--(4,-1.2);
    \draw[line width=0.3mm,gray] (3,-0.4) node{\tiny{$S$}};
    \draw[line width=0.3mm] (4,-3.5)--(-4,-3.5);
    \foreach \Point in {(3,-3.5), (2,-3.5),(0,-3.5),(1,-3.5),(-1,-3.5),(-0.7,-3.5),(-2,-3.5),(-3,-3.5)}
    \draw[fill=black] \Point circle (0.05);
    \draw[line width=0.4mm,gray] (3,-3.8) node{\tiny{$x_8$}} (2,-3.8) node{\tiny{$x_7$}} (0,-3.8) node{\tiny{$x_5$}} (1,-3.8) node{\tiny{$x_6$}} (-1,-3.8) node{\tiny{$x_3$}} (-0.7,-3.8) node{\tiny{$x_4$}} (-2,-3.8) node{\tiny{$x_2$}} (-3,-3.8) node{\tiny{$x_1$}};
    \end{tikzpicture}
    \caption{A universally monotone zigzag cover of type $(0,1^5,1^5,\undl x)$.}
    \label{fig:asym-uni-monzigzag2}
\end{figure}
By shrinking or extending the length of the glued inner edge properly,
we get an universally monotone zigzag cover $\bar\varphi$ of type $(g,(\lambda,(1)^{2m}),(\mu,(1)^{2m}),\undl{\bar x})$.
By a similar argument as that in the case $(1)$, we have $\vec N(\bar\varphi)\geq m!$ and $\vec\zl_{g}((\lambda, 1^{2m}),(\mu, 1^{2m}))\geq m!$.

{\bf Case $(3)$:} $l(\lambda_o)=0$, $l(\mu_o)=2$ and $\mu_{o_1}\neq1,\mu_{o_2}=1$.

Let $\varphi:C\to T\pb^1$ be a universally monotone zigzag cover of type $(g,\lambda,\mu,\undl x)$.
The string in $C$ has two out-ends.
We glue the weight $1$ out-end of the string in $C$ with the in-end of the string of a universally monotone zigzag cover $\tilde\varphi:\tilde C\to T\pb^1$ of type $(0,(1^{2m+1}),(1^{2m+1}),\undl{\tilde x})$ (See Figure $\ref{fig:mon-zigzag1}$ for an example).
By shrinking or extending the length of the glued inner edge properly,
we get a universally monotone zigzag cover $\bar\varphi$ of type $(g,(\lambda,(1)^{2m}),(\mu,(1)^{2m}),\undl{\bar x})$.
By a similar argument as that in the case $(1)$, we have $\vec N(\bar\varphi)\geq m!$ and $\vec\zl_{g}((\lambda, 1^{2m}),(\mu, 1^{2m}))\geq m!$.

{\bf Case $(4)$:} $l(\lambda_o)=l(\mu_o)=0$.

The proof of this case is the same as the argument in case $(2)$, so we omit it here.
\end{proof}

\begin{corollary}
\label{cor:asymp-arbitrary-splitting}
Under the same conditions of Theorem $\ref{thm:asymp-arbitrary-splitting}$,
the logarithmic asymptotics for $\vec\zl_{g,\lambda,\mu}(m)$ is at least $m\log m$ as $m\to\infty$, where $\vec\zl_{g,\lambda,\mu}(m)=\vec\zl_{g}((\lambda, 1^{2m}),(\mu, 1^{2m}))$.
\end{corollary}

\begin{proof}
Since $\log m!\sim m\log m$ as $m\to\infty$,
Corollary $\ref{cor:asymp-arbitrary-splitting}$ follows from Theorem $\ref{thm:asymp-arbitrary-splitting}$.
\end{proof}

\section{Application to real mixed double Hurwitz numbers}
\label{sec:5}

\subsection{Mixed zigzag numbers}
The construction of universally monotone zigzag covers can also be used to study the logarithmic asymptotics for the lower bounds of mixed double Hurwitz numbers.

\begin{definition}\label{def:mixed-fact}
Let $\undl S(s)=\{\sk_1,\sk_2,\ldots,\sk_r\}$ be a sequence of signs with $s$ positive entries.
Suppose that $(\gamma,\sigma_1,\tau_1,\ldots,\tau_r,\sigma_2)$ is a tuple of type $(g,\lambda,\mu)$ with $\tau_i=(a_i,b_i)$,
where $a_i<b_i$, $i=1,\ldots,r$.
This tuple is a \textit{real $k$-mixed factorization} of type $(g,\lambda,\mu;\undl S(s),k)$,
if it satisfies all the conditions listed in
Definition $\ref{def:real-factor}$ and the following condition:
$$
b_i\leq b_{i+1}, ~\forall i\in\{1,2,\ldots,k-1\}.
$$
\end{definition}
We denote by $\vec\fl^\rb(g,\lambda,\mu;\undl S(s),k)$
the set of real $k$-mixed
factorizations of type $(g,\lambda,\mu;\undl S(s),k)$.
Let
\begin{equation}\label{eq:mixed-Hurwitz}
\vec{H}^\rb_g(\lambda,\mu;\undl S(s),k):=|\vec\fl^\rb(g,\lambda,\mu;\undl S(s),k)|.
\end{equation}
The number $\vec{H}^\rb_g(\lambda,\mu;\undl S(s),k)$ is called the
\textit{real $k$-mixed double Hurwitz number} with $s$ positive branch points under the sequence $\undl S(s)$.

Let $\varphi:C\to T\pb^1$ be a tropical cover,
and let $V$ be a subset of inner vertices of $C$.
Suppose that $C'\subset C$ is a subgraph of $C$ containing $V$.
The subgraph $C'$ is called {\it a minimal subgraph containing $V$}, if there does not exist a proper subgraph $C''$ of $C'$ such that $V$ is also a subset of inner vertices of $C''$.

\begin{definition}
\label{def:k-mixed-zig}
A zigzag cover $\varphi:C\to T\pb^1$ of type $(g,\lambda,\mu,\undl x)$ is called
\textit{$k$-mixed} if there exists a string $S\subset C$ such that $\varphi$ satisfies the following two conditions:
\begin{itemize}
    \item the restriction of $\varphi$ on the
    minimal subgraph $C'$ containing $\varphi^{-1}(\undl x^k)$ is a universally monotone zigzag cover of certain type, where $\undl x^k=\{x_1,\ldots,x_k\}$.
    \item if the string $S\not\subset C'$, the subgraph $S\cap(C\setminus C')$ is connected.
\end{itemize}
\end{definition}
We denote by $\vec\ml_{g}(\lambda,\mu;k)$ the set of $k$-mixed zigzag covers of type $(g,\lambda,\mu,\undl x)$.
Let $\vec{N}_k(\varphi,\undl S(s))$ denote the number of real $k$-mixed factorizations associated to a real $k$-mixed zigzag cover $(\varphi,\rho)$, where $\undl S(s)$ is a sequence of signs with $s$ positive entries, and $\rho$ is the unique colouring of $\varphi$ which is determined by the splitting $\undl x=\undl x^+\sqcup\undl x^-$ associated to the sequence of signs $\undl S(s)$.
We set
\begin{equation}
    \vec{\zl}_g(\lambda,\mu;k):=
    \sum_{\varphi\in\vec{\ml}_{g}(\lambda,\mu;k)}\vec{N}_k(\varphi),
\end{equation}
where $\vec{N}_k(\varphi)=\min_{s,\undl S(s)}\vec{N}_k(\varphi,\undl S(s))$.
The number $\vec{\zl}_g(\lambda,\mu;k)$ is called
\textit{the $k$-mixed zigzag number}.

\begin{remark}
The $k$-mixed zigzag numbers
$\vec{\zl}_g(\lambda,\mu;k)$ do not depend on the number $s$ of the positive real branch points. We refer the readers to \cite[Remark $5.3$]{rau2019} or \cite[Remark $3.9$]{d-2020} for more details.
\end{remark}

\begin{proposition}\label{prop:lower-bound-mixed}
Fix $g\geq0$, $d\geq1$, and two partitions $\lambda$, $\mu$ with $|\lambda|=|\mu|=d$.
Then the real $k$-mixed double Hurwitz
number $\vec{H}^\rb_g(\lambda,\mu;\undl S(s),k)$ is bounded from below by the $k$-mixed zigzag number:
$$
\vec\zl_{g}(\lambda,\mu;k)\leq
\vec{H}^\rb_g(\lambda,\mu;\undl S(s),k).
$$
\end{proposition}

\begin{proof}
It is straightforward from the definition, so we omit it.
\end{proof}

\begin{lemma}
\label{lem:k-mix-index}
Let $\varphi:C\to T\pb^1$ be a $k$-mixed zigzag cover of type $(g,\lambda,\mu,\undl x)$, and $C'$ be the minimal subgraph of $C$ containing $\varphi^{-1}(\undl x^k)$,
where $\undl x^k=\{x_1,x_2,\ldots,x_k\}$.
Suppose that $\varphi|_{C'}$ is a degree $d_\varphi$ universally monotone zigzag cover.
Then the number
$$
\vec{N}_k(\varphi)\geq (d-d_\varphi)!.
$$
\end{lemma}

\begin{proof}
Let $\rho_s$ be a colouring of the $k$-mixed zigzag cover $\varphi:C\to T\pb^1$
possessing an arbitrary splitting $\undl x=\undl x^+\sqcup\undl x^-$ with $|\undl x^+|=s$.
From Lemma $\ref{lem:mult-real-trop}$,
there are $d!\cdot\mult^\rb(\varphi,\rho_s)$ real
factorizations of type $(g,\lambda,\mu;\undl S(s))$
associated to $(\varphi,\rho_s)$.
Let $C'$ be the minimal subgraph having
$\varphi^{-1}(\undl x^k)$ as inner vertices.
Then $\varphi|_{C'}:C'\to T\pb^1$ is a universally monotone zigzag cover of degree $d_\varphi$.
From Lemma $\ref{lem:mono-zigzag1}$,
there is a real monotone factorization
$(\gamma_1,\sigma_1',\tau_1,\ldots,\tau_k,\sigma_2')$
producing ($\varphi|_{C'},\rho_s|_{C'}$).

If the string $S\subset C'$,
the permutation $\sigma_1'=\sigma_1$ and $d_\varphi=d$.
The statement of Lemma $\ref{lem:k-mix-index}$ follows from Lemma $\ref{lem:mono-zigzag1}$ for the case $k=r$.
In the case $k<r$, a similar argument as in the proof of Theorem $\ref{thm:asym-simple-splitting}$ shows that
there is at least one real factorization
$(\gamma_1,(\sigma_2')^{-1},\tau_{k+1},\ldots,\tau_r,\sigma_2)$
such that $(\gamma_1,\sigma_1,\tau_1,\ldots,\tau_r,\sigma_2)$ produces the real $k$-mixed zigzag cover $(\varphi,\rho_s)$. Hence, we have $\vec N_k(\varphi)\geq 1$.

If the string $S\not\subset C'$,
we know that $S\cap(C\setminus C')$ is connected.
Let $C''$ be the minimal connected subgraph of $C$ containing $\varphi^{-1}(\undl x\setminus\undl x^k)$.
Suppose that the degree of the restricted zigzag cover $\varphi|_{C''}:C''\to T\pb^1$ is $d_1$.
From Lemma $\ref{lem:mult-real-trop}$,
there are $d_1!\cdot\mult^\rb(\varphi|_{C''},\rho_s|_{C''})$ real
factorizations
associated to $(\varphi|_{C''},\rho_s|_{C''})$.
Suppose that the weight of the in-end of the string $S|_{C''}$, which connects to $S|_{C'}$, is $l$.
We want to calculate the number of real factorizations
$(\gamma_2,\sigma_1'',\tau_{k+1},\ldots,\tau_{r},\sigma_2'')$
associated to $(\varphi|_{C''},\rho_s|_{C''})$
such that the real factorizations $(\gamma_1,\sigma_1',\tau_1,\ldots,\tau_k,\sigma_2')$
and $(\gamma_2,\sigma_1'',\tau_{k+1},\ldots,\tau_r,\sigma_2'')$
form real factorizations which
produce the real $k$-mixed zigzag cover $(\varphi,\rho_s)$.
Once $(\gamma_1,\sigma_1',\tau_1,\ldots,\tau_k,\sigma_2')$ is fixed,
the permutation associated to the bridge edge connecting
$S|_{C'}$ with $S|_{C''}$ is also fixed by $\sigma_2'$.
This contributes a factor $\frac{l\cdot(d_1-l)!}{d_1!}$ to the total number of required real
factorizations
associated to $(\varphi|_{C''},\rho_s|_{C''})$,
so there are at least $l\cdot(d_1-l)!\cdot\mult^\rb(\varphi|_{C''},\rho_s|_{C''})$ real factorizations producing the real $k$-mixed zigzag cover $(\varphi,\rho_s)$.
Since $\varphi|_{C''}:C''\to T\pb^1$ is a zigzag cover,
it follows from \cite[Proposition $4.7$]{rau2019} and Remark $\ref{rem:mult-comparision}$ that the multiplicity $\mult^\rb(\varphi|_{C''},\rho_s|_{C''})$ is an odd number.
From Definition $\ref{def:k-mixed-zig}$,
we get $d_1-l=d-d_\varphi$.
Therefore, $\vec N_k(\varphi)\geq(d-d_\varphi)!$.
\end{proof}

\subsection{Asymptotics for real mixed double Hurwitz numbers}
In this section, we prove the logarithmic equivalence of real mixed double Hurwitz numbers and complex double Hurwitz numbers under certain conditions.

\begin{theorem}
\label{thm:asym-mix-zigzag1}
Let $\lambda$ and $\mu$ be two partitions of $d$ with $l(\lambda_o)=l(\mu_o)=1$.
Suppose that there are two partitions $\lambda'$ and $\mu'$ with $|\lambda'|=|\mu'|\leq d$ such that
\begin{itemize}
    \item $\lambda_o'=\lambda_o$, $\lambda'\setminus\lambda_o'\subset\lambda_{e}$ and $\mu'\setminus\mu_o'\subset\mu_{e}$;
    \item $\lambda'_{o,o}=\mu'_{o,o}=\emptyset$, $l(\lambda'_o)=l(\mu'_o)=1$ and $l(\lambda')+l(\mu')=k+2$;
    \item $\mu'_i+\mu'_j>\max(\lambda_o,\lambda'_{max})$ for any two entries $\mu'_i$, $\mu'_j$ of $\mu'_{e}$,
    and $\lambda'_{max}>\mu_o'$, where $\lambda'_{max}$ is the maximal integer in $\lambda'_e$.
\end{itemize}
Then for any $m\geq1$, there is at least one $k$-mixed zigzag cover $\varphi:C\to T\pb^1$ of type $(g,(\lambda,1^{2m}), (\mu,1^{2m}), \undl x)$.
Moreover, there is an integer $m_0$ such that for any $m\geq m_0$ the $k$-mixed zigzag number
$$
\vec\zl_{g,\lambda,\mu;k}(m)\geq (m-m_0)!^4\cdot(2m)!,
$$
where $\vec\zl_{g,\lambda,\mu;k}(m)=\vec{\zl}_{g}((\lambda,1^{2m}),(\mu,1^{2m});k)$.
\end{theorem}

\begin{proof}
Since $\mu'_i+\mu'_j>\max(\lambda_o,\lambda'_{max})$ for any two entries $\mu'_i$, $\mu'_j$ of $\mu'_{e}$ and $\lambda'_{max}>\mu_o'$,
the construction of monotone zigzag cover of type $(0,\lambda',\mu',\undl x')$ in the case $(1)$ of the proof of Theorem $\ref{thm:asym-simple-splitting}$ produces a universally monotone zigzag cover of type $(0,\lambda',\mu',\undl x')$.
Let $\varphi':C'\to T\pb^1$ be a universally monotone zigzag cover of type $(0,\lambda',\mu',\undl x')$.
From \cite[Proposition $5.2$]{rau2019}, there is a zigzag cover $\varphi'':C''\to T\pb^1$ of type $(g,(\lambda\setminus\lambda',\mu'_o,1^{1m}),(\mu\setminus(\mu'\setminus\mu'_o),1^{2m}),\undl{x''})$.
The out-end of the string $S'\subset C'$ and the in-end of the string $S''\subset C''$ are of the same weight $\mu'_o$.
We glue the out-end of the string $S'\subset C'$ with the in-end of the string $S''\subset C''$,
then we obtain a graph $C$. See Figure $\ref{fig:exa-mix-zigzag-1}$ for an example.
\begin{figure}[ht]
    \centering
    \begin{tikzpicture}
    \draw[line width=0.3mm,gray] (-5,2.5)--(-1,2)--(-3,1.5)--(-2,1.2)--(-4,0.8)--(3,-0.2)--(1,-0.7)--(2,-0.9)--(1.5,-1.2)--(5,-1.5);
    \draw[line width=0.3mm] (-1,2)--(5,2);
    \draw[line width=0.3mm] (-5,1.5)--(-3,1.5);
    \draw[line width=0.3mm] (-2,1.2)--(5,1.2);
    \draw[line width=0.3mm] (-5,0.8)--(-4,0.8);
    \draw[line width=0.3mm] (5,-0.2)--(3,-0.2);
    \draw[line width=0.3mm] (-0.3,-0.3)--(2.6,-0.3);
    \draw[line width=0.3mm] (-5,-0.6)--(-0.3,-0.3)--(-5,0);
    \draw[line width=0.3mm] (0.5,-0.7)--(1,-0.7);
    \draw[line width=0.4mm] (0.5,-0.7) arc[start angle=0, end angle=360, x radius=0.2, y radius=0.1];
    \draw[line width=0.3mm] (-0.1,-0.7)--(0.1,-0.7);
    \draw[line width=0.4mm] (-0.1,-0.7) arc[start angle=0, end angle=360, x radius=0.2, y radius=0.1];
    \draw[line width=0.3mm] (-5,-0.7)--(-0.5,-0.7);
    \draw[line width=0.3mm] (2,-0.9)--(4.5,-0.9);
    \draw[line width=0.3mm] (5,-0.7)--(4.5,-0.9)--(5,-1.1);
    \draw[line width=0.3mm] (-5,-1.2)--(1.5,-1.2);
    \draw[line width=0.3mm] (-5,-1.75)--(5,-1.75);
    \foreach \Point in {(-1,-1.75), (-3,-1.75), (-2,-1.75), (-4,-1.75), (3,-1.75), (1,-1.75), (2,-1.75), (1.5,-1.75), (-0.5,-1.75), (-0.3,-1.75),(-0.1,-1.75),(0.1,-1.75),(0.5,-1.75),(2.6,-1.75),(4.5,-1.75)}
    \draw[fill=black] \Point circle (0.05);
    \draw[line width=0.3mm,gray] (-3,2.4) node{\tiny{$S$}} (-0.7,0.5) node{\tiny{The glued edge}};
    \end{tikzpicture}
    \caption{A $k$-mixed zigzag cover $\varphi:C\to T\pb^1$, where $k=4$.}
    \label{fig:exa-mix-zigzag-1}
\end{figure}
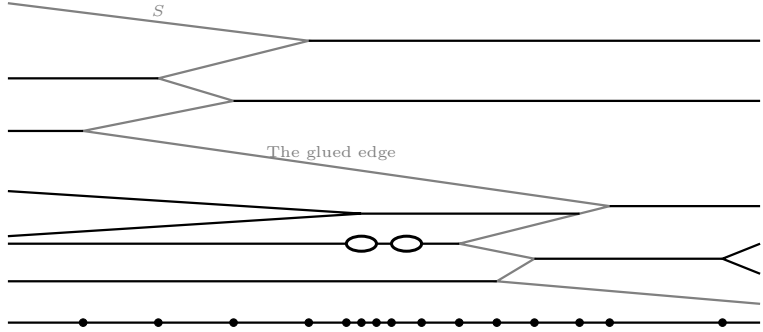
After shrinking or extending the lengths of the edges of $C$ properly,
there is a total order of the inner vertices of $C$ such that the inner vertices of $C'$ are the first $k$ ones. Then we obtain a $k$-mixed zigzag cover $\varphi:C\to T\pb^1$ of type $(g,(\lambda,1^{2m}),(\mu,1^{2m}),\undl x)$.

By Lemma $\ref{lem:k-mix-index}$, the number $\vec{N}_k(\varphi)\geq(d+2m-d_\varphi)!$,
where $d_\varphi=\deg(\varphi')=|\lambda'|$.
Note that $|\lambda'|\leq|\lambda|=d$,
so $\vec{N}_k(\varphi)\geq(2m)!$.
From \cite[Proposition $5.9$]{rau2019}, for sufficiently large $m$,
there are at least $(m-m_0)!^4$ zigzag covers of type $(g,(\lambda\setminus\lambda',\mu'_o,1^{1m}),(\mu\setminus(\mu'\setminus\mu'_o),1^{2m}),\undl{x''})$, where $m_0$ is a fixed integer.
Hence, there are at least $(m-m_0)!^4$ $k$-mixed zigzag covers of type $(g,(\lambda,1^{2m}),(\mu,1^{2m}),\undl x)$ with $\vec{N}_k(\varphi)\geq(2m)!$ for sufficiently large $m$.
Therefore, the $k$-mixed zigzag number $\vec\zl_{g,\lambda,\mu;k}(m)\geq(m-m_0)!^4\cdot(2m)!$ for any $m\geq m_0$.
\end{proof}

\begin{corollary}
\label{cor:asym-mix-zigzag1}
Under the same conditions of Theorem $\ref{thm:asym-mix-zigzag1}$,
the $k$-mixed zigzag number
$\vec\zl_{g,\lambda,\mu;k}(m)$
is logarithmically equivalent to the complex double Hurwitz number:
$$
\vec\zl_{g,\lambda,\mu;k}(m)
\sim 6m\log m\sim h^\cb_{g,\lambda,\mu}(m),
\text{ as }m\to\infty,
$$
where $h^\cb_{g,\lambda,\mu}(m)= H^\cb_g((\lambda,1^{2m}),(\mu,1^{2m}))$.
\end{corollary}

\begin{proof}
Let $H^\cb_g(m)=H^\cb_g((1)^m,(1)^m)$.
It follows from \cite[Equation $5$]{dyz-2017} and the proof of \cite[Theorem $5.10$]{rau2019} that $\log h^\cb_{g,\lambda,\mu}(m)\leq
\log H^\cb_g(m)\sim 6m\log m$ as $m\to\infty$.
Since $\log[(m-m_0)!^4\cdot(2m)!]\sim6m\log m$,
the statement follows from Theorem $\ref{thm:asym-mix-zigzag1}$.
\end{proof}

\begin{remark}
It follows from Proposition $\ref{prop:lower-bound-mixed}$ and Corollary $\ref{cor:asym-mix-zigzag1}$ that the real $k$-mixed double Hurwitz numbers are logarithmically equivalent to the complex double Hurwitz numbers.
Our construction can also be used to analyze the asymptotic growth of $k$-mixed zigzag numbers when adding both monotone and ordinary simple branch points or only adding monotone simple branch points.
When considering the asymptotics for these cases, the arguments are similar to what we have done in Section $\ref{sec:asym-real-mono}$,
so we restrict ourselves to the case that only ordinary simple branch points are added.
\end{remark}


\begin{thebibliography}{10}

\bibitem{bbm-2011}
B.~Bertrand, E.~Brugall\'e, and G.~Mikhalkin.
\newblock {Tropical open Hurwitz numbers}.
\newblock {\em Rend. Semin. Mat. Univ. Padova}, 125:157--171, 2011.

\bibitem{cadoret-2005}
A.~Cadoret.
\newblock {Counting real Galois covers of the projective line}.
\newblock {\em Pacific J. Math.}, 219(1):53--81, 2005.

\bibitem{cjm-2010}
R.~Cavalieri, P.~Johnson, and H.~Markwig.
\newblock {Tropical Hurwitz numbers}.
\newblock {\em J. Algebraic Combin.}, 32(2):241--265, 2010.

\bibitem{cjm-2011}
R.~Cavalieri, P.~Johnson, and H.~Markwig.
\newblock {Wall crossings for double Hurwitz numbers}.
\newblock {\em Adv. Math.}, 228(4):1894--1937, 2011.

\bibitem{cm-2016}
R.~Cavalieri and E.~Miles.
\newblock {\em {Riemann surfaces and algebraic curves: A first course in
  Hurwitz theory}}, volume~87 of {\em London Mathematical Society Student
  Texts}.
\newblock Cambridge University Press, Cambridge, 2016.

\bibitem{d-2020}
Y.~Ding.
\newblock {On the lower bounds for real double Hurwitz numbers}.
\newblock {\em arXiv:2010.00899}, 2020.

\bibitem{dk-2017}
N.~Do and M.~Karev.
\newblock {Monotone orbifold Hurwitz numbers}.
\newblock {\em J. Math. Sci. (N.Y.)}, 226(5):568--587, 2017.

\bibitem{dl-2022}
N.~Do and D.~Lewa\'nski.
\newblock {On the Goulden-Jackson-Vakil conjecture for double Hurwitz numbers}.
\newblock {\em Adv. Math.}, 403:31, 2022.

\bibitem{dyz-2017}
B.~Dubrovin, D.~Yang, and D.~Zagier.
\newblock {Classical Hurwitz numbers and related combinatorics}.
\newblock {\em Mosc. Math. J.}, 17(4):601--633, 2017.

\bibitem{elsv-2001}
T.~Ekedahl, S.~Lando, M.~Shapiro, and A.~Vainshtein.
\newblock {Hurwitz numbers and intersections on moduli spaces of curves}.
\newblock {\em Invent. Math.}, 146(2):297--327, 2001.

\bibitem{er-2019}
B.~El~Hilany and J.~Rau.
\newblock {Signed counts of real simple rational functions}.
\newblock {\em J Algebraic Combin.},
  https://doi.org/10.1007/s10801-019-00906-6, 2019.

\bibitem{gz2018}
P.~Georgieva and A.~Zinger.
\newblock {Real Gromov-Witten theory in all genera and real enumerative
  geometry: construction}.
\newblock {\em Ann. of Math. (2)}, 188(3):685--752, 2018.

\bibitem{ggpn-2013}
I.~P. Goulden, M.~Guay-Paquet, and J.~Novak.
\newblock {Monotone Hurwitz numbers in genus zero}.
\newblock {\em Canad. J. Math.}, 65(5):1020--1042, 2013.

\bibitem{ggpn-2013a}
I.~P. Goulden, M.~Guay-Paquet, and J.~Novak.
\newblock {Polynomiality of monotone Hurwitz numbers in higher genera}.
\newblock {\em Adv. Math.}, 238:1--23, 2013.

\bibitem{ggpn-2014}
I.~P. Goulden, M.~Guay-Paquet, and J.~Novak.
\newblock {Monotone Hurwitz numbers and the HCIZ integral}.
\newblock {\em Ann. Math. Blaise Pascal}, 21(1):71--89, 2014.

\bibitem{ggpn-2017}
I.~P. Goulden, M.~Guay-Paquet, and J.~Novak.
\newblock {On the convergence of monotone Hurwitz generating functions}.
\newblock {\em Ann. Comb.}, 21(1):73--81, 2017.

\bibitem{gjv-2005}
I.~P. Goulden, D.~M. Jackson, and R.~Vakil.
\newblock {Towards the geometry of double Hurwitz numbers}.
\newblock {\em Adv. Math.}, 198(1):43--92, 2005.

\bibitem{gpmr-2015}
M.~Guay-Paquet, H.~Markwig, and J.~Rau.
\newblock {The combinatorics of real double Hurwitz numbers with real positive
  branch points}.
\newblock {\em Int. Math. Res. Not. IMRN}, 2016(1):258--293, 2016.

\bibitem{hahn-2019}
M.~A. Hahn.
\newblock {A monodromy graph approach to the piecewise polynomiality of simple,
  monotone and Grothendieck dessins d'enfants double Hurwitz numbers}.
\newblock {\em Graphs Combin.}, 35(3):729--766, 2019.

\bibitem{hkl-2018}
M.~A. Hahn, R.~Kramer, and D.~Lewa\'{n}ski.
\newblock {Wall-crossing formulae and strong piecewise polynomiality for mixed
  Grothendieck dessins d'enfant, monotone, and double simple Hurwitz numbers}.
\newblock {\em Adv. Math.}, 336:38--69, 2018.

\bibitem{hl-2020}
M.~A. Hahn and D.~Lewa\'{n}ski.
\newblock {Wall-crossing and recursion formulae for tropical Jucys covers}.
\newblock {\em Trans. Amer. Math. Soc.}, 373(7):4685--4711, 2020.

\bibitem{hl-2022}
M.~A. Hahn and D.~Lewanski.
\newblock {Tropical Jucys covers}.
\newblock {\em Math. Z.}, 301(2):1719--1738, 2022.

\bibitem{hurwitz-1891}
A.~Hurwitz.
\newblock {Ueber Riemann'sche Fl\"{a}chen mit gegebenen Verzweigungspunkten}.
\newblock {\em Math. Ann.}, 39(1):1--60, 1891.

\bibitem{iks2003}
I.~Itenberg, V.~Kharlamov, and E.~Shustin.
\newblock {Welschinger invariant and enumeration of real rational curves}.
\newblock {\em Int. Math. Res. Not.}, 49:2639--2653, 2003.

\bibitem{iks2004}
I.~Itenberg, V.~Kharlamov, and E.~Shustin.
\newblock {Logarithmic equivalence of the Welschinger and the Gromov-Witten
  invariants}.
\newblock {\em Russian Math. Surveys}, 59(6):1093--1116, 2004.

\bibitem{iks2007}
I.~Itenberg, V.~Kharlamov, and E.~Shustin.
\newblock {New cases of logarithmic equivalence of Welschinger and
  Gromov-Witten invariants}.
\newblock {\em Proc. Steklov Inst. Math.}, 258(1):65--73, 2007.

\bibitem{iks2013b}
I.~Itenberg, V.~Kharlamov, and E.~Shustin.
\newblock {Welschinger invariants of real del Pezzo surfaces of degree $
  \geqslant3$}.
\newblock {\em Math. Ann.}, 355(3):849--878, 2013.

\bibitem{iz-2018}
I.~Itenberg and D.~Zvonkine.
\newblock {Hurwitz numbers for real polynomials}.
\newblock {\em Comment. Math. Helv.}, 93(3):441--474, 2018.

\bibitem{johnson-2015}
P.~Johnson.
\newblock {Double Hurwitz numbers via the infinite wedge}.
\newblock {\em Trans. Amer. Math. Soc.}, 367(9):6415--6440, 2015.

\bibitem{ks-2015}
V.~Kharlamov and R.~R$\breve{a}$sdeaconu.
\newblock {Counting real rational curves on K3 surfaces}.
\newblock {\em Int. Math. Res. Not. IMRN}, 2015(14):5436--5455, 2015.

\bibitem{kls-2019}
R.~Kramer, D.~Lewanski, and S.~Shadrin.
\newblock {Quasi-polynomiality of monotone orbifold Hurwitz numbers and
  Grothendieck's dessins d'enfants}.
\newblock {\em Doc. Math.}, 24:857--898, 2019.

\bibitem{lzz-2000}
A.-M. Li, G.~Zhao, and Q.~Zheng.
\newblock {The number of ramified covering of a Riemann surface by Riemann
  surface}.
\newblock {\em Comm. Math. Phys.}, 213(3):685--696, 2000.

\bibitem{mr-2015}
H.~Markwig and J.~Rau.
\newblock {Tropical real Hurwitz numbers}.
\newblock {\em Math. Z.}, 281(1-2):501--522, 2015.

\bibitem{mikhalkin-2005}
G.~Mikhalkin.
\newblock {Enumerative tropical algebraic geometry in $\mathbb{R}^2$}.
\newblock {\em J. Amer. Math. Soc.}, 18(2):313--377, 2005.

\bibitem{okounkov-2000}
A.~Okounkov.
\newblock {Toda equations for Hurwitz numbers}.
\newblock {\em Math. Res. Lett.}, 7(4):447--453, 2000.

\bibitem{op-2006}
A.~Okounkov and R.~Pandharipande.
\newblock {Gromov-Witten theory, Hurwitz theory, and completed cycles}.
\newblock {\em Ann. of Math. (2)}, 163(2):517--560, 2006.

\bibitem{rau2019}
J.~Rau.
\newblock {Lower bounds and asymptotics of real double Hurwitz numbers}.
\newblock {\em Math. Ann.}, 375(1-2):895--915, 2019.

\bibitem{ssv-2008}
S.~Shadrin, M.~Shapiro, and A.~Vainshtein.
\newblock {Chamber behavior of double Hurwitz numbers in genus 0}.
\newblock {\em Adv. Math.}, 217(1):79--96, 2008.

\bibitem{shustin2015}
E.~Shustin.
\newblock {On higher genus Welschinger invariants of del Pezzo surfaces}.
\newblock {\em Int. Math. Res. Not. IMRN}, 16:6907--6940, 2015.

\bibitem{wel2005a}
J.-Y. Welschinger.
\newblock Invariants of real symplectic $4$-manifolds and lower bounds in real
  enumerative geometry.
\newblock {\em Invent. Math.}, 162(1):195--234, 2005.

\bibitem{wel2005b}
J.~Y. Welschinger.
\newblock {Spinor states of real rational curves in real algebraic convex
  3-manifolds and enumerative invariants}.
\newblock {\em Duke Math. J.}, 127(1):89--121, 2005.

\end{thebibliography}

\end{document}